\documentclass{amsart}

\usepackage[all,hyperref]{bi-discrete}
\usepackage[backend=biber,maxbibnames=5,maxalphanames=5,style=alphabetic,bibencoding=utf8,giveninits,url=false,isbn=false]{biblatex}
\usepackage{tikz}

\usetikzlibrary{positioning, arrows.meta, shapes.geometric, shadows, calc}
\AtBeginBibliography{\small}
\SetSymbolFont{stmry}{bold}{U}{stmry}{m}{n}

\newcommand*\osc{\mathop{}\!\mathrm{osc}_k}
\newcommand*{\oscp}[1]{\mathop{}\!\mathrm{osc}_{k,#1}}
\newcommand*{\osckp}[2]{\mathop{}\!\mathrm{osc}_{#1,#2}}
\newcommand*\PiL{\Pi^\mathrm{L}_k}
\newcommand*\PiSZ{\Pi^\mathrm{SZ}_k}
\newcommand*\PiP{\Pi^\mathrm{P}_{B^\Omega_R}}
\newcommand\tilu{\widetilde{u}}
\newcommand\tiluh{\widetilde{u}_h}
\newcommand{\diff}[1]{\mathop{}\!\mathrm{d}#1}
\newcommand{\barp}{\bar{p}}
\newcommand{\vertex}{\mathtt{z}}
\newcommand{\setX}{\mathrm{X}}
\newcommand{\setF}{\mathrm{F}}
\newcommand{\setA}{\mathrm{A}}
\renewcommand{\setW}{\mathrm{W}}
\renewcommand{\setL}{\mathrm{L}}
\renewcommand{\setC}{\mathrm{C}}
\renewcommand{\setB}{\mathrm{B}}

\author{Lars Diening}
\address[L. Diening]{Fakultät für Mathematik, Universität Bielefeld, 33615 Bielefeld, Germany.}
\email{lars.diening@uni-bielefeld.de}
\urladdr{http://www.bi-discrete.com}

\author{Julian Rolfes}
\address[J. Rolfes]{Fakultät für Mathematik, Universität Bielefeld, 33615 Bielefeld, Germany.}
\email{jrolfes@math.uni-bielefeld.de}
\urladdr{http://www.bi-discrete.com}

\author{Abner J. Salgado}
\address[A.J. Salgado]{Department of Mathematics, University of Tennessee, Knoxville, TN 37996, USA}
\email{asalgad1@utk.edu}
\urladdr{https://math.utk.edu/people/abner-salgado/}

\title{Regularity transfer for the Ritz projection}

\begin{document}



\begin{abstract}
  Let $\Omega\subset\RRn$ be a convex polytope in arbitrary dimension $n\geq 1$, and consider the Ritz projection on quasi-uniform triangulations of $\Omega$. We prove a new pointwise estimate in terms of the sharp maximal function and use it to establish stability in weighted Lebesgue, Orlicz, and Lorentz spaces. We then introduce discrete higher-order smoothness spaces based on oscillations and a smooth companion construction, and prove discrete analogues of classical embedding and compactness results. Finally, we establish stability of the Ritz projection in these discrete spaces, together with localized estimates. These results provide a framework for transferring higher-order regularity through the Ritz projection. As part of our analysis, we also study approximation properties of the Scott--Zhang projection in Besov and Triebel--Lizorkin spaces.
\end{abstract}

\maketitle

\tableofcontents

\section{Introduction}
\label{sec:introduction}

This work is dedicated to the study of stability and approximation properties of the Ritz projection in non-energy norms. The Ritz projection maps the solution of the Laplace problem to its discrete counterpart. More precisely, let $\Omega \subset \RRn$ with $n \geq 1$ be a bounded and convex polytopal domain, $\mathcal{T}_h$ be a quasi-uniform, shape regular mesh, and $V_h \subset \setW^{1,\infty}_0(\Omega)$ denote the Lagrange finite element space consisting of piecewise polynomials of degree $k \geq 1$. Then the Ritz projection $R_h u$ of a function $u \in \setW^{1,1}_0(\Omega)$ is the unique element of $V_h$ that satisfies
\begin{equation}
\label{eq:RitzProjection}
  \int_\Omega \nabla R_h u \cdot \nabla \varphi_h \diff{x} = \int_\Omega \nabla u \cdot \nabla \varphi_h \diff{x}, \qquad \forall \varphi_h \in V_h.
\end{equation}
We are interested in determining for which function spaces $\setX(\Omega)$ the estimate
\begin{align}
  \label{eq:intro-ritz-stability}
  \norm{ \nabla R_h u }_{\setX(\Omega)} \lesssim \norm{ \nabla u }_{\setX(\Omega)}
\end{align}
is valid. This stability is equivalent to the best approximation property
\[
  \norm{ \nabla ( u - R_h u) }_{\setX(\Omega)} \lesssim \inf_{v_h \in V_h} \norm{ \nabla ( u - v_h) }_{\setX(\Omega)}.
\]
This follows immediately by a simple application of the triangle inequality.
The mere definition of the Ritz projection implies its stability and optimal approximation properties for $\setX(\Omega) = \setL^2(\Omega)$.

The study of these properties in other norms, however, is not a trivial endeavor. Under varying sets of assumptions this question has been addressed in many works. The case $\setX(\Omega) = \setL^\rho(\Omega)$, $\rho \in [1, \infty]$, is the content of several papers. The case $\rho=\infty$ was studied in \cite{MR645661} for $n=2$. This was extended to $n=3$ in \cite[Chapter 8]{MR2373954} and \cite{MR2495783}. Those papers also contain an interpolation and duality argument to extend to $\rho\in(1,\infty]$.  The case $\rho=1$ has been treated in~\cite{MR5022107}.
For mildly graded meshes, see \cite{MR2869035}.

Another important example is provided by weighted Lebesgue spaces $\setX(\Omega) = \setL^\rho_\omega(\Omega)$ for $\rho \in (1,\infty)$. For $\omega \in \mathcal{A}_{\rho/2}$ with $\rho\geq 2$, this was positively answered in \cite{MR4060457} for $n = 2,3$. The general case $\omega \in \mathcal{A}_\rho$ with $\rho\in(1,\infty)$ is relevant for the Poisson problem with singular forces like the Dirac delta function and was solved in~\cite{MR4748207}. For $\rho=1$, this was again done in~\cite{MR5022107}.

The crucial component of our previous work \cite{MR4748207} is the pointwise estimate
\begin{equation}
\label{eq:PtWiseEstimate}
  \abs{ \nabla R_hu} \lesssim M[\nabla u] 
\end{equation}
for $n = 2,3$, where $M[\cdot]$ denotes the Hardy--Littlewood maximal function.  This pointwise estimate immediately implies the stability in any Banach function space $\setX(\Omega)$, in the sense of~\cite[Definition 2.1]{Lerner2017}, on which the maximal function is norm bounded.  Banach function spaces have the ideal property: if $\abs{v}\leq\abs{w}$ almost everywhere, then $\norm{v}_{\setX(\Omega)}\leq\norm{w}_{\setX(\Omega)}$. This includes for example the Lebesgue spaces $\setL^\rho(\Omega)$ and the weighted Lebesgue spaces $\setL^\rho_\omega(\Omega)$ for $1<\rho \leq \infty$ and $\omega \in \mathcal{A}_\rho$. 

Although the estimate~\eqref{eq:PtWiseEstimate} with its consequences is very strong, it is restricted to $n=2,3$ and cannot be applied to $\setX(\Omega) = \setBMO(\Omega)$. The paper \cite{MR1677478} shows boundedness of~$R_h$ for $\setX(\Omega) = \setBMO(\Omega)$ and $n \geq 1$. It also covers the case of Morrey spaces.

Building on the oscillation estimates of~\cite{MR1677478},
we establish a pointwise comparison between the oscillations
of the discrete and continuous gradients.
More precisely, we prove
\begin{align}
  \label{eq:intro-pointwise-sharp-aux}
  M^\sharp_\tau[\nabla R_h u]
  &\lesssim M^\sharp_\tau[\nabla u],
  \qquad \tau\in(1,2].
\end{align}
Unlike the pointwise estimate~\eqref{eq:PtWiseEstimate},
this estimate holds in every space dimension $n\geq 1$.
It yields stability of the Ritz projection in
Lebesgue spaces $\setL^\rho(\Omega)$ and weighted
Lebesgue spaces $\setL^\rho_\omega(\Omega)$ for
$\rho\in(1,\infty)$ and $\omega\in\mathcal A_\rho$,
as well as in $\setBMO(\Omega)$.
Further applications to Banach function spaces are
discussed in Section~\ref{sec:banach}.

An additional advantage of~\eqref{eq:intro-pointwise-sharp-aux}
is that it allows us to control the discrete gradient
directly in terms of the data.
Indeed, if $u$ solves
$-\Delta u=\divergence F$ in $\Omega$ with zero boundary
data, then the continuous sharp maximal estimate
from Appendix~\ref{sec:cont-sharp-maxim} gives
\begin{align}
  M^\sharp_\tau[\nabla R_h u]
  &\lesssim M^\sharp_\tau[\nabla u]
  \lesssim M^\sharp_\tau[F]
  \lesssim M_\tau[F].
\end{align}
Such a data estimate cannot be obtained from
\eqref{eq:PtWiseEstimate}: even in two dimensions,
there are bounded vector fields $F$ for which
$\nabla u\in\setBMO(\Omega)\setminus\setL^\infty(\Omega)$.
Consequently, $M[\nabla u]$ cannot be controlled
pointwise by $M_\tau[F]$.
An explicit example is given in
Remark~\ref{rem:counterexample-linfty}.

Furthermore, we extend the stability of~$R_h$ to spaces with additional smoothness like differentiability. In particular, we will consider Campanato spaces, Sobolev spaces, Triebel--Lizorkin spaces and Besov spaces.  In doing so we go beyond the intrinsic Lipschitz regularity of discrete functions.  To this end, we introduce discrete Triebel--Lizorkin and Besov spaces that are insensitive to scales below the mesh size $h$.  These spaces will have similar embedding and compactness properties as their non-discrete counterparts.  Additionally, we introduce a so-called \emph{smooth companion}, which is a smooth approximation of the finite element function and possesses similar approximation properties.  The stability estimates we establish in Section~\ref{sec:StabRitz} then have the form
\begin{align}
  \label{eq:intro-BF-stab}
  \norm{\nabla R_h u}_{\setA^s_{\rho,q,h}(\Omega)}
  &\lesssim
  \norm{\nabla u}_{\setA^s_{\rho,q}(\Omega)}.
\end{align}
Here, $\setA \in \{\setB,\setF\}$ may denote Besov or Triebel--Lizorkin spaces and $\setA^s_{\rho,q,h}(\Omega)$ our discrete smoothness versions.
Related discrete elliptic regularity estimates are available in mesh-dependent Sobo\-lev norms. For $\setW^{2,\rho}$-estimates measured in broken norms of the gradients for continuous Lagrange finite element functions, see~\cite{FengHenningsNeilan2017}. For discontinuous Galerkin methods we refer to~\cite{SmearsSuli2013} for discrete $\setW^{2,2}$-estimates and to~\cite{FengNeilanSchnake2018} for discrete $\setW^{2,\rho}$-estimates for interior penalty methods.

\section{Geometry and discretization}
\label{sec:Prelims}

In the following, the relation $A \lesssim B$ means that there exists a nonessential constant $c > 0$ such that $A \leq cB$. The value of this constant may change at each occurrence. The relation $A \eqsim B$ is shorthand for $A \lesssim B \lesssim A$. Equality by definition will be denoted by $A \coloneqq B$. As expected, $B \eqqcolon A$ means $A \coloneqq B$. For $a \in \RR$ we set $a_+ \coloneqq \max\{a,0\}$.

We will always assume that $\Omega \subset \RRn$ ($n \geq 1$) is a convex polytopal domain. For $x_0 \in \RRn$ and $R>0$ we denote by $B_R(x_0)$ the ball of radius $R$ centered at $x_0$. If the center is irrelevant we will then write $B_R$, whereas if the center and radius are irrelevant we simply write $B$. We also define $B_R^\Omega(x_0) \coloneqq B_R(x_0) \cap \Omega$, and use the notational shorthands $2B_R(x_0) \coloneqq B_{2R}(x_0)$ and $2B_R \coloneqq B_{2R}$. We call $D \subseteq \Omega$ regular (in the sense of~\cite[Definition 1.1]{MR2247585}), if for every $x_0 \in \overline{D}$ there holds
\begin{equation}
\label{eq:BiteOffEstimate}
  \frac{ |B_t(x_0) | }{ |B_t(x_0) \cap D| } \lesssim 1, \qquad \forall t \in (0,\diameter D),
\end{equation}
where the implicit constant is independent of $x$ and $t$, but it may depend on $D$.
Since $\Omega$ is convex and polytopal, it is regular. If $D$ satisfies the interior cone condition or is Lipschitz, then it is regular.

If $D\subset \RRn$ denotes a bounded domain and $\setX(D)$ is a space of functions over $D$, we will denote by $\norm{\cdot}_{\setX(D)}$ the standard norm in this space, while $\abs{\cdot}_{\setX(D)}$ means the seminorm. We will adhere to standard notation regarding most function spaces, or provide an explanation when we deviate from it.

Since $\Omega$ is a polytope it can be meshed exactly. By $\{ \mathcal{T}_h \}_{h>0}$ we denote a quasi-uniform, in the sense of \cite[Definition 22.20]{MR4242224} family of triangulations of $\Omega$, where $h>0$ is the mesh size, such that $\diameter T \leq h$ for every $T \in \mathcal{T}_h$. In particular, for all $h>0$ and $T \in \mathcal{T}_h$, we have $\abs{T} \eqsim h^n$. For $T \in \mathcal{T}_h$ we define
\[
  \omega(T) \coloneqq \{ T' \in \mathcal{T}_h \ : \ T' \cap T \ne \emptyset \}
\]
to be the patch, or neighborhood, of $T$. We assume that the elements $T \in \mathcal{T}_h$ are closed. If $\wp_h \subset \mathcal{T}_h$, we denote
\[
  \omega(\wp_h) \coloneqq \bigcup_{T \in \wp_h} \omega(T).
\]
With these definitions we can introduce higher order patches recursively. If $m \in \setN$, then $\omega^{m+1}(T) \coloneqq \omega( \omega^m(T))$ and $\omega^{m+1}(\wp_h) \coloneqq \omega( \omega^m(\wp _h) )$. In an abuse of notation, we may use these symbols to denote collections of elements, or their union.
 
Since some of our statements will be of local nature, we need to define discrete subdomains.
\begin{definition}[Discrete subdomain]
  \label{def:discrete-domain}
  A subdomain $D \subseteq \Omega$ is a \emph{discrete subdomain} on the mesh $\mathcal{T}_h$ if it satisfies the following properties:
  \begin{enumerate}[label={\alph{*})},leftmargin=*]
    \item \textbf{Submesh:} The closure of $D$ is the union of elements of $\mathcal{T}_h$.
    \item \textbf{Regular:} $D$ is regular in the sense of~\eqref{eq:BiteOffEstimate}.
    \item \textbf{Minimal size:} There holds $\diameter D \geq 4h$.
  \end{enumerate}
\end{definition}
Note that $\Omega$ is always a discrete subdomain.

Now Let $m \in \setN_0$ and $r \in \setN_0$. The space of all polynomials (in $n$ variables) of degree at most $m$ is denoted by $\mathcal{P}_m$. We define, for each $h>0$,
\[
  \mathcal{L}_m^r(\mathcal{T}_h) \coloneqq \Big\{ w_h \in \setC^{r-1}(\overline{\Omega}) \ : \ w_{h|T} \in \mathcal{P}_m \text{ for all } T \in \mathcal{T}_h \Big\},
\]
where, for simplicity of notation only, we set $\setC^{-1}(\overline{\Omega}) \coloneqq \setL^\infty(\Omega)$. Consequently, elements of $\mathcal{L}_m^0(\mathcal{T}_h)$ are piecewise polynomial, but are not globally continuous and can have jumps on the faces of the simplices that comprise $\mathcal{T}_h$.

The finite element space that we shall work with is as follows. We fix $k \in \setN$ and define, again for each $h>0$,
\[
  V_h \coloneqq \mathcal{L}_k^1(\mathcal{T}_h) \cap \setW^{1,2}_0(\Omega).
\]
In our proofs we will need two standard projections onto finite element spaces. The first one is the classical Lagrange interpolant $\PiL \colon \setC^0(\overline{\Omega}) \to \mathcal{L}^1_k(\mathcal{T}_h)$, defined by interpolation of a continuous function at the nodes.
The second projection we shall make use of is the so-called Scott--Zhang interpolation operator $\PiSZ \colon \setW^{1,1}(\Omega) \to \mathcal{L}^1_k(\mathcal{T}_h)$, defined in \cite{MR1011446}. We use the classical version that preserves discrete boundary data. We recall that this operator is $\setW^{1,1}$-stable and it satisfies the following approximation estimate: for every $T\in \mathcal{T}_h$ and $l = 0,1$,
\begin{equation}
\label{eq:SZApprox}
  \norm{\nabla^l(v - \PiSZ v)}_{\setL^2(T)} \lesssim h^{2-l} \abs{\nabla^2 v}_{\setL^2(\omega(T))}, \qquad \forall v \in \setW^{1,1}(\Omega) \cap\setW^{2,2}(\omega(T)).
\end{equation}
We finally mention that both of these projections preserve homogeneous boundary values.

\section{Oscillation estimates}
\label{sec:Oscillation}

In this section we begin the proof of the regularity transfer properties of discrete solutions. The main result of this section is the restricted oscillation estimate in Theorem~\ref{thm:osc-estimate}, which for a ball $B \subset \Omega$ controls the restricted oscillation of $\nabla R_h u$, where $R_h u$ is defined in \eqref{eq:RitzProjection}, in terms of the restricted oscillation of $\nabla u$ in a slightly larger ball plus a lower order term.
Our oscillation results are higher order improvements of the $\setBMO$ estimates in~\cite{MR1677478}. Our technique is similar and is based on \Caccioppoli{} and comparison arguments.

The oscillation here is defined as follows. Let $p \in [1,\infty]$ and $v \in \setL^p(\Omega)$. For $m \in \setN$, the local oscillation of $v$ of order $m$ at point $x \in \Omega$ and scale $t>0$ is
\begin{equation}
\label{eq:LocOscillation}
\mathrm{osc}_{m,p} v(x,t) \coloneqq 
  \begin{dcases}
    \bigg(\inf_{q \in \mathcal{P}_{m-1}} \dashint_{B_t^\Omega(x)} \abs{v-q}^p \diff{y} \bigg)^{1/p}, & p<\infty,
    \\
    \inf_{q \in \mathcal{P}_{m-1}} \norm{v-q}_{\setL^\infty(B_t^\Omega(x))} , & p = \infty.
  \end{dcases}
\end{equation}
For $D \subseteq \Omega$ regular we also define the restricted oscillation of $v$ of order $m$ at point $x \in D$ and scale $t>0$ by
\begin{align*}
  \mathrm{osc}_{m,p}^D v(x,t) \coloneqq 
    \begin{dcases}
    \bigg(\inf_{q \in \mathcal{P}_{m-1}} \dashint_{D \cap B_t(x)} \abs{v-q}^p \diff{y} \bigg)^{1/p}, & p<\infty,
        \\
        \inf_{q \in \mathcal{P}_{m-1}} \norm{v-q}_{\setL^\infty(D \cap B_t(x))} , & p = \infty.
    \end{dcases}
\end{align*}

With this notation we may state the main result of this section.

\begin{theorem}[Oscillation reduction]
\label{thm:osc-estimate}
  For $\beta \in (0,k)$, there are $\theta \in (0,1)$ and  $\Lambda >1$ such that, if $R > \Lambda h$, $u \in \setW^{1,2}_0(\Omega)$, the following estimates hold.
  \begin{enumerate}[label={\alph{*})},leftmargin=*]
    \item If $x \in \Omega$ is such that $B_R(x) \Subset \Omega$, then
    \[
      \oscp{2}^{B_R} \nabla R_h u(x,\theta R) \leq \theta^\beta \oscp{2}^{B_R} \nabla R_h u (x,R) + c_\beta \oscp{2}^{B_R} \nabla u (x,R).
    \]
    
  \item There exists $\sigma \in (0,1]$ depending on the geometry of $\Omega$ such that provided $\beta \in (0,\sigma)$, for every $x \in \overline\Omega$ there holds
    \[
    \osckp{1}{2}^{B_R} \nabla R_h u(x,\theta R) \leq \theta^\beta \osckp{1}{2}^{B_R} \nabla R_h u (x,R) + c_\beta \osckp{1}{2}^{B_R} \nabla u (x,R).
    \]
  \end{enumerate}
\end{theorem}

\begin{remark}[Alternative version]
  In Theorem~\ref{thm:osc-estimate-ii} we extend Theorem~\ref{thm:osc-estimate} and replace $\oscp{2}$ by $\oscp{\tau}$ for $\tau \in (1,2)$.
\end{remark}

As a first step we compare the local best approximation, in the $\setW^{1,2}$-seminorm, via polynomials of degree $m \in \setN$, against the best approximation of a gradient, in the $\setL^2$-norm, via polynomials of degree $m-1$. For $m=1$ these two quantities are equivalent, since every vector $\xi \in \RRn$ is the gradient of a linear polynomial. For $m>1$ this is more involved, and quantified in the following result.

\begin{lemma}[Equivalent oscillation]
  \label{lem:inf_equiv}
  Let $D \subset \Omega$ be star-shaped with respect to a ball $B$ and $v \in \setW^{1,\rho}(\Omega)$ with $1\leq \rho \leq \infty$. For every $m \in \setN$ we have
  \begin{equation}
  \label{eq:inf_equiv}
    \inf_{p \in \mathcal{P}_m} \norm{\nabla (v-p)}_{\setL^\rho(D)}
    \eqsim
    \inf_{q \in (\mathcal{P}_{m-1})^n} \norm{\nabla v-q}_{\setL^\rho(D)},
  \end{equation}
  where the implicit constant depends on $m$ and $\abs{D}/\abs{B}$.
\end{lemma}
\begin{proof}
  It is clear that the right hand side is smaller than the left one. Thus we focus on the other direction.  To prove \eqref{eq:inf_equiv}, we use the averaged Taylor polynomial of order $m$, see \cite{MR559195} or \cite[Section 4.1]{MR2373954}.
  The averaged Taylor operators, extended to integrable functions by integration by parts in the averaging variable, satisfy, for $j=0,\ldots,m$,
  \begin{align*}
    \frT_B^j(p) &= p
    && \forall p\in\mathcal P_j,\\
    \norm{\frT_B^j(w)}_{\setL^\rho(D)}
                &\lesssim \norm{w}_{\setL^\rho(B)}
    && \forall w\in\setL^\rho(D),\\
    \nabla\frT_B^m(v) &= \frT_B^{m-1}(\nabla v)
    && \forall v\in\setW^{1,1}(D).
  \end{align*}
  When applying these, for $q \in (\mathcal{P}_{m-1})^n$, we can calculate
  \begin{align*}
    \inf_{p \in \mathcal{P}_m} \norm{\nabla (v-p)}_{\setL^\rho(D)}
    &\leq \norm{\nabla (v-\frT_B^m(v))}_{\setL^\rho(D)} 
    = \norm{\nabla v -\frT_B^{m-1}(\nabla v)}_{\setL^\rho(D)} \\
    &= \norm{\nabla v -q - \frT_B^{m-1}(\nabla v-q)}_{\setL^\rho(D)} 
    \lesssim \norm{\nabla v -q}_{\setL^\rho(D)}.
  \end{align*}
  Since $q$ was chosen arbitrarily, this gives the desired estimate.
\end{proof}
We now state an auxiliary oscillation estimate which, together with Lemma~\ref{lem:inf_equiv}, will allow us to obtain Theorem~\ref{thm:osc-estimate}.
\begin{lemma}[Oscillation estimate]
\label{lem:campanato}
  There is $\Lambda >0$ such that if $R > \Lambda h$ and $\rho \in (h,\tfrac12 R)$ then, for every $u \in \setW^{1,2}_0(\Omega)$, the following estimates hold:
  \begin{enumerate}[label={\alph{*})},leftmargin=*]
  \item Let $m \in \setN$ with $m \le k$. For every ball $B_R \Subset \Omega$ and all concentric balls $B_\rho \subset B_R$,
    \begin{align*}
      \inf_{p \in \mathcal{P}_m} \dashint_{B_\rho} \abs{\nabla (R_h u - p)}^2\diff{x}
      &\leq C_1 \bigg( \frac{\rho^{2m}}{R^{2m}} + \frac{h}{R}
      \frac{R^n}{\rho^n}\bigg) \inf_{p \in \mathcal{P}_m} \dashint_{B_R}
      \abs{\nabla (R_h u-p)}^2\diff{x} \\
      &+ C_2 \frac{R^n}{\rho^n} \inf_{p \in \mathcal{P}_m}\dashint_{B_R} \abs{\nabla (u-p)}^2\diff{x}.
    \end{align*}
    
  \item There exists $\sigma \in (0,1]$ only depending on the geometry of $\Omega$ such that for every ball $B_R$ with a center in $\Omega$, and all concentric balls $B_\rho \subset B_R$,
    \begin{align*}
      \inf_{\xi \in \RRn} \dashint_{B^\Omega_\rho} \abs{\nabla R_h u - \xi}^2\diff{x}
      &\leq C_1 \bigg( \frac{\rho^{2\sigma}}{R^{2\sigma}} + \frac{h}{R}
      \frac{R^n}{\rho^n}\bigg) \inf_{\xi \in \RRn} \dashint_{B^\Omega_R}
      \abs{\nabla R_h u-\xi}^2\diff{x} \\
      &+ C_2 \frac{R^n}{\rho^n} \inf_{\xi \in \RRn}\dashint_{B^\Omega_R} \abs{\nabla u-\xi}^2\diff{x}.
    \end{align*}
  \end{enumerate}
\end{lemma}

Let us now show how, from Lemma~\ref{lem:campanato}, the proof of Theorem~\ref{thm:osc-estimate} easily follows.

\begin{proof}[Proof of Theorem~\ref{thm:osc-estimate}]
  This is a consequence of Lemma \ref{lem:campanato}. In view of Lemma \ref{lem:inf_equiv}, there holds
  \[
    \big(\oscp{2}^{B_R} \nabla R_h u(x,R)\big)^2 \eqsim \inf_{p \in \mathcal{P}_k} \dashint_{B^\Omega_R} \abs{\nabla(R_h u-p)}^2 \diff{x},
  \]
  meaning this is just a matter of adjusting constants.

  Consider an interior ball. If $\theta \coloneqq \frac{\rho}{R}$ is small enough, there holds
  \[
    C_1 \Big( \frac{\rho}R \Big)^{2k} = C_1 \theta^{2k}  \leq \frac12 \theta^{2\beta}.
  \]
  Having fixed $\theta$ we, if necessary, further increase $\Lambda$, thus making $\frac{h}{R} < \tfrac1\Lambda$ small enough, so as to guarantee
  \[
    C_1  \frac{h}R  \bigg( \frac{R}\rho \bigg)^n = C_1 \frac{h}R  \theta^{-n} \leq C_1 \Lambda^{-1} \theta^{-n} \leq \frac12 \theta^{2\beta}.
  \]

  For balls that intersect the boundary the argument is similar, provided we restrict $\beta \in (0,\sigma)$.
\end{proof}

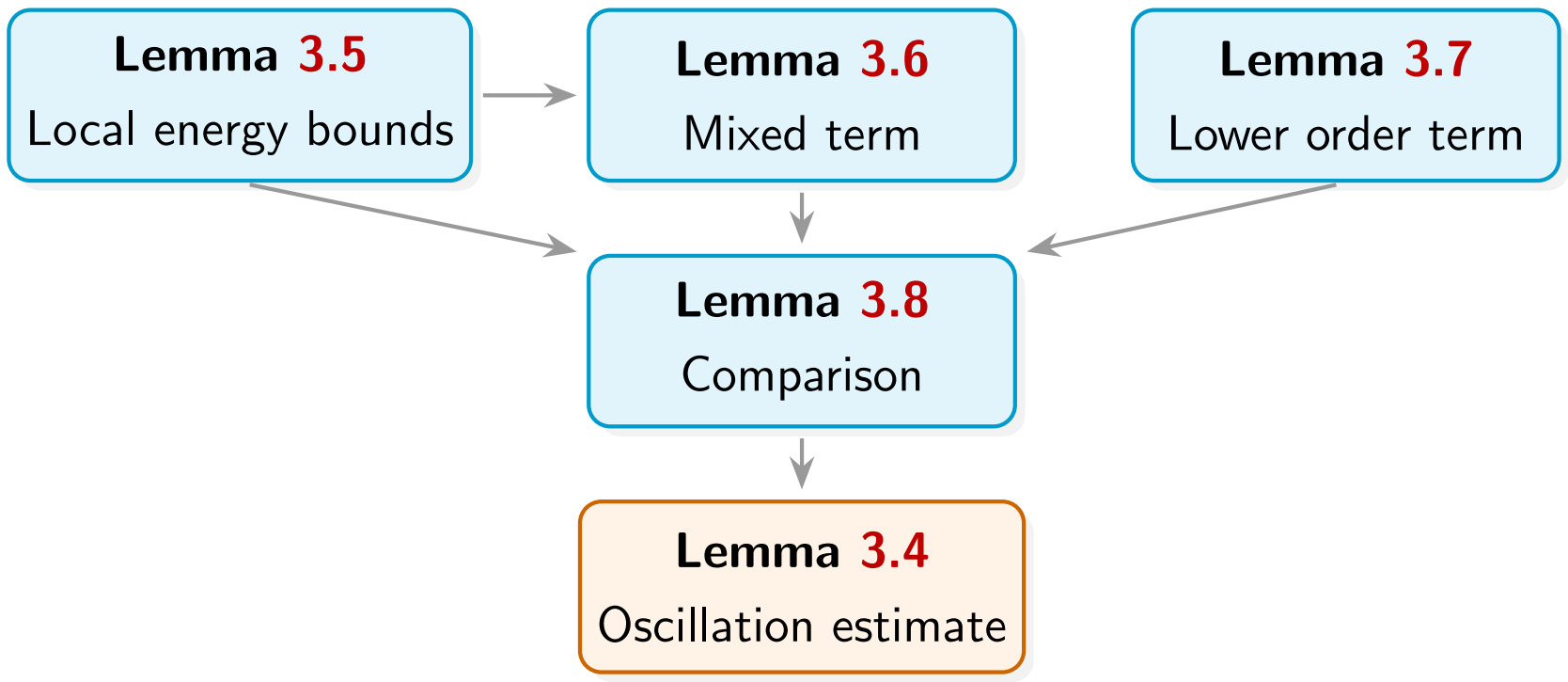
\begin{figure}[htbp]
  \begin{tikzpicture}[
      % --- Style Definitions ---
      basebox/.style={
          rectangle,
          rounded corners=1ex,
          draw,
          thick,
          align=center,
          minimum width=3cm,
          minimum height=1.2cm,
          font=\sffamily,
          drop shadow={opacity=0.1, shadow xshift=2pt, shadow yshift=-2pt}
      },
      theorem/.style={
          basebox,
          fill=orange!10,
          draw=orange!80!black
      },
      lemma/.style={
          basebox,
          fill=cyan!10,
          draw=cyan!80!black
      },
      myarrow/.style={
          -Stealth,
          thick,
          draw=gray!80,
          shorten >=2pt,
          shorten <=2pt
      }
  ]

      \node[lemma] (lem6) {\textbf{Lemma \ref{lem:comparison}}\\[1mm] Comparison};

      \node[lemma, above=0.5cm of lem6] (lem4) {\textbf{Lemma \ref{lem:mixed-term}}\\[1mm] Mixed term};
      \node[lemma, left=0.8cm of lem4] (lem3) {\textbf{Lemma \ref{lem:grad-psi-estimate}}\\[1mm] Local energy bounds};
      \node[lemma, right=0.8cm of lem4] (lem5) {\textbf{Lemma \ref{lem:lower-order-estimate}}\\[1mm] Lower order term};

      \node[theorem, below=0.5cm of lem6] (lem2) {\textbf{Lemma \ref{lem:campanato}}\\[1mm] Oscillation estimate};

      \draw[myarrow] (lem3) -- (lem4);

      \draw[myarrow] (lem3.south) -- (lem6.north west);
      \draw[myarrow] (lem4.south) -- (lem6.north);
      \draw[myarrow] (lem5.south) -- (lem6.north east);

      \draw[myarrow] (lem6) -- (lem2);

  \end{tikzpicture}
  \caption{Visual roadmap of the proof for Lemma \ref{lem:campanato}.}
  \label{fig:proof-roadmap}
\end{figure}

The remainder of this section will be dedicated to the proof of Lemma~\ref{lem:campanato}. For clarity of exposition, we organize this in a series of auxiliary results. Figure~\ref{fig:proof-roadmap} shows the connections between the statements of this section.

We begin with a simple reduction. Let, from now on, $x_0 \in \overline\Omega$ and $R>0$ be fixed. For every $r>0$ we let $B_r \coloneqq B_r(x_0)$. If $B_R \not\subset \Omega$ then we set $q_R \equiv 0$. Otherwise, that is if $B_R \Subset \Omega$, we fix $m \in \setN$ with $m \leq k$, and we let $q_R \in \mathcal{P}_m$ be arbitrary. Then we define
\[
  \tilu \coloneqq u-q_R, \qquad \qquad \tiluh \coloneqq R_h u-q_R.
\]
Notice that these definitions depend on $x_0$ and $R$. In order not to overburden the notation we do not indicate this, with the understanding that 
functions are related by
\begin{align}
  \label{eq:tiluh-def-equation}
\begin{alignedat}{2}
  \int_\Omega \nabla \tiluh \cdot\nabla \phi_h \diff{x} &= \int_\Omega \nabla \tilu \cdot \nabla \phi_h \diff{x}
                                                   \qquad && \forall \phi_h \in V_h
    \\
  \tiluh &= \tilu = -q_R && \text{ on } \partial\Omega.
\end{alignedat}
\end{align}

Having defined $\tilu$, we fix $\barp \in \mathcal{P}_m$ that satisfies
\begin{equation}
\label{eq:pMinimizesLocalH1seminorm}
  \barp \in \argmin_{p \in \mathcal{P}_m} \int_{B^\Omega_R} \abs{\nabla(\tilu-p)}^2 \diff{x},
\end{equation}
this polynomial is unique up to constants. 

Assume now that $B_R \Subset \Omega$. Since $q_R$ was chosen arbitrarily, to prove the interior estimate of Lemma~\ref{lem:campanato} it is then enough to prove that
\begin{equation}
\label{eq:campanato-tilde-estimate}
  \begin{aligned}
    \inf_{q \in \mathcal{P}_m} \dashint_{B_\rho} \abs{\nabla (\tiluh - q)}^2\diff{x}
    &\lesssim \bigg( \frac{\rho^{2m}}{R^{2m}} + \frac{h}{R}
    \frac{R^n}{\rho^n}\bigg) \dashint_{B_R} \abs{\nabla \tiluh}^2\diff{x}
    \\
    &+ \frac{R^n}{\rho^n} \dashint_{B_R}
    \abs{\nabla (\tilu-\barp)}^2\diff{x}.
  \end{aligned}
\end{equation}
We are now ready to begin the proof.
  
\subsection{Localization and comparison}
  
We begin by having a closer look at the polynomial $\barp$. By an orthogonal decomposition of $\setW^{1,2}(B^\Omega_R)$ we can represent $\barp = p_0 + p_H$, where $p_H$ is harmonic and $p_0$ vanishes on $\partial B^\Omega_R$. In this case, we have $\Delta p_0 = \Delta \barp$. Furthermore, owing to \cite[Theorem 5.5]{MR1805196}, if $B_R^\Omega = B_R$ (an interior ball), there holds $p_0,p_H \in \mathcal{P}_m$.

As a means to establishing the claimed oscillation estimates we compare $\tiluh$ to a harmonic function. This way, many of the properties we shall need can be reduced to those of harmonic functions, which are well-known. Let $w \in \setW^{1,2}(B^\Omega_R)$ be the solution to
\begin{equation}
  \label{eq:comparison-vw}
  \begin{dcases}
    -\Delta w = 0, & \text{in } B^\Omega_R,
      \\
    w =\tiluh, & \text{on } \partial B^\Omega_R.
  \end{dcases}
\end{equation}

Our goal now is to compare $\tiluh$ with the harmonic function $w$, defined in \eqref{eq:comparison-vw}. For that we need to localize the functions $\tiluh$ and $\tilu$. To this end we introduce a discrete cut-off function $\zeta_h \in \mathcal{L}^1_1(\mathcal{T}_h)$ such that $0 \leq \zeta_h \leq 1$,
\begin{equation}
\label{eq:DefOfzetahCutOff}
  \zeta_h \equiv 0 \text{ in } \overline{\Omega} \setminus B_{\frac34R}^\Omega,
  \qquad \qquad \zeta_h \equiv 1 \text{ in } B_{\frac12R}^\Omega,
  \qquad \qquad 
  \norm{ \nabla \zeta_h }_{\setL^\infty(\Omega)} \lesssim \frac1R.
\end{equation}
For this to be possible, a sufficiently large $\Lambda$ must be chosen. Furthermore, we choose $\Lambda$ large enough such that also
\[
  \bigcup_{T \cap \support \zeta_h \neq \emptyset} \omega(T) \subset B_{\frac78R}.
\]
With this at hand, we can define the functions
\begin{equation}
    \label{eq:def-psi}
  \begin{aligned}
    \psi &\coloneqq \zeta_h^2 (\tiluh - p_0 - w), \\
    \psi_h &\coloneqq \PiL (\zeta_h^2 (\tiluh - p_0 - \PiSZ w)).
  \end{aligned}
\end{equation}

We will need bounds on the energy of $\psi$, $\psi_h$, and their difference.

\begin{lemma}[Local energy bounds]
\label{lem:grad-psi-estimate}
  For $\psi$ and $\psi_h$, defined in \eqref{eq:def-psi}, we have
  \[
    \int_{B^\Omega_R} \abs{\nabla \psi}^2 \diff{x} \lesssim \int_{B^\Omega_R} \zeta_h^2 \abs{\nabla(\tiluh-p_0-w)}^2 \diff{x} + \frac1{R^2} \int_{B^\Omega_R} \abs{\tiluh-p_0-w}^2 \diff{x},
  \]
  as well as the local estimate on the difference
  \begin{align*}
    \int_T \abs{\nabla(\psi-\psi_h)}^2 \diff{x}
    \lesssim \sup_{x \in T} \zeta_h^2(x) \Big( h^2 \int_{\omega(T)} \abs{\nabla^2 w}^2 \diff{x} + \frac1{R^2}\int_{T} \abs{\tiluh-p_0-w}^2 \diff{x}\Big),
  \end{align*}
  and for every $\xi \in \RRn$ the global estimate
  \begin{align*}
    \int_{B^\Omega_R} \abs{\nabla(\psi-\psi_h)}^2 \diff{x}
    \lesssim \frac{h^2}{R^2} \int_{B^\Omega_R} \abs{\nabla \tiluh -\xi}^2 \diff{x} + \frac1{R^2} \int_{B^\Omega_R} \abs{\tiluh-p_0-w}^2 \diff{x}.
  \end{align*}
\end{lemma}
\begin{proof}
  The first inequality is a direct consequence of the identity
  \begin{align*}
    \nabla \psi = \zeta_h^2\nabla(\tiluh-p_0-w)+2\zeta_h\nabla\zeta_h(\tiluh-p_0-w)
  \end{align*}
  and the fact that $\norm{\nabla\zeta_h}_{\setL^\infty(\Omega)} \lesssim R^{-1}$.

  The last estimate follows from summing the second one over all elements $T$ such that $T \cap B_R^\Omega \neq \emptyset$, using the support properties of $\zeta_h$, and applying Lemma~\ref{lem:properties-w}.
  It remains to prove the local energy estimate on the difference $\psi - \psi_h$. By expanding $\nabla\psi$ as before and using, again, that $\abs{\nabla\zeta_h} \lesssim R^{-1}$ there holds
  \begin{align*}
    \lefteqn{\int_T \abs{\nabla(\psi-\psi_h)}^2 \diff{x}} \quad & \\
    &\lesssim \int_T \abs{\zeta_h\nabla\zeta_h(\tiluh-p_0-w)}^2 \diff{x}
    + \int_T \abs{\zeta_h^2\nabla(\tiluh-p_0-w)-\nabla\psi_h}^2 \diff{x} \\
    &\lesssim \frac{ \sup_{x \in T} \zeta_h^2(x) }{R^2} \int_T \abs{\tiluh-p_0-w}^2 \diff{x}
    + \int_T \abs{\zeta_h^2\nabla(\tiluh-p_0-w)-\nabla\psi_h}^2 \diff{x}.
  \end{align*}
  For the second term, we use the Scott--Zhang projection to estimate
  \begin{align*}
    \lefteqn{\int_T \abs{\zeta_h^2\nabla(\tiluh-p_0-w)-\nabla\psi_h}^2 \diff{x}}
    \quad & \\
          &\lesssim \int_T \abs{\zeta_h^2\nabla(\tiluh-p_0-\PiSZ w)-\nabla\psi_h}^2 \diff{x}
         + \sup_{x \in T} \zeta_h^2(x)\int_T \abs{\nabla(w-\PiSZ w)}^2\diff{x}.
  \end{align*}
  By interpolation properties of the Scott--Zhang projection we have
  \begin{align*}
    \sup_{x \in T} \zeta_h^2(x) \int_T \abs{\nabla(w-\PiSZ w)}^2\diff{x}
    &\lesssim h^2\sup_{x \in T} \zeta_h^2(x) \int_{\omega(T)} \abs{\nabla^2 w}^2\diff{x}.
  \end{align*}
  Finally, let $\{ \vertex_i \}_{i\in I}$ be the Lagrange nodes of $T$ associated with the Lagrange interpolant, and $\{\oldphi_i\}_{i\in I} \subset \mathcal{P}_k$ the corresponding basis functions. Plugging in the definition of $\psi_h$ we get
  \begin{align*}
    \MoveEqLeft \int_T \abs{\zeta_h^2\nabla(\tiluh-p_0-\PiSZ w)-\nabla\psi_h}^2 \diff{x} = \\
    &= \int_T \bigabs{\sum_{i \in I} \big(\zeta_h^2(x)-\zeta_h^2(\vertex_i)\big) (\tiluh-p_0-\PiSZ w)(\vertex_i)\nabla\oldphi_i(x)}^2 \diff{x}.
  \end{align*}
  By shape regularity as well as scaling arguments there holds $\norm{\nabla \oldphi_i}_{\setL^\infty(T)} \lesssim h^{-1}$. We can also obtain the estimate
  \[
    \abs{\zeta_h^2(x)-\zeta_h^2(\vertex_i)} \lesssim \frac{h}{R} \sup_{x \in T} \zeta_h(x).
  \]
  Thus, using that $\{\vertex_i\}_{i\in I}$ is a unisolvent set for $\mathcal{P}_k$,
  \begin{align*}
    \MoveEqLeft \int_T \abs{\zeta_h^2\nabla(\tiluh-p_0-\PiSZ w)-\nabla\psi_h}^2 \diff{x} \\
    &\lesssim \sup_{x \in T} \zeta_h^2(x) \frac{\abs{T}}{R^2} \sum_{i \in I} \abs{(\tiluh-p_0-\PiSZ w)(\vertex_i)}^2 \\
    &\lesssim \sup_{x \in T} \zeta_h^2(x) \frac1{R^2}\int_T \abs{\tiluh-p_0-\PiSZ w}^2 \diff{x} \\
    &\lesssim \sup_{x \in T} \zeta_h^2(x) \frac1{R^2}\bigg(\int_T \abs{\tiluh-p_0-w}^2 \diff{x}
    + \int_T \abs{w-\PiSZ(w)}^2 \diff{x} \bigg) .
  \end{align*}
  The last term can again be bounded by interpolation estimates.
\end{proof}

We now have all the tools to give two auxiliary estimates that will play an important role in the proof of Lemma~\ref{lem:campanato}.

\begin{lemma}[Mixed term]
\label{lem:mixed-term}
With the previously established definitions we have, for any $\epsilon \in (0,1)$ and all $\xi \in \RRn$,
  \begin{align*}
    \MoveEqLeft \int_{B^\Omega_R} \nabla(\tiluh-p_0-w) \cdot \nabla(\psi-\psi_h) \diff{x} \\
    &\lesssim \epsilon \int_{B^\Omega_R} \zeta_h^2\abs{\nabla(\tiluh-p_0-w)}^2 \diff{x}
    + c_\epsilon \frac{h^2}{R^2} \int_{B^\Omega_R} \abs{\nabla \tiluh -\xi}^2 \diff{x} \\
    & + c_\epsilon \frac{1}{R^2} \int_{B^\Omega_R} \abs{\tiluh-p_0-w}^2 \diff{x}.
  \end{align*}
\end{lemma}
\begin{proof}
  For conciseness, within this proof, we set $\chi \coloneqq \tiluh - p_0 - w$. We split the integral into those $T$ for which $T \cap B_R^\Omega \neq \emptyset$, on each integral use the Cauchy--Schwarz inequality, then the local estimate of Lemma~\ref{lem:grad-psi-estimate}, and finally Young's inequality to arrive at
  \begin{align*}
    \MoveEqLeft \int_{B_R^\Omega} \nabla \chi \cdot \nabla (\psi- \psi_h) \diff{x} \lesssim \\
    &\sum_{ \substack{ T \in \mathcal{T}_h \\ T \cap \support \zeta_h \neq \emptyset}}
    \bigg[
      \epsilon \sup_{x \in T} \zeta_h^2(x) \int_T \abs{ \nabla \chi}^2 \diff{x}
      +
    c_\epsilon h^2 \int_{\omega(T)} \abs{ \nabla^2 w }^2 \diff{x} + \frac{c_\epsilon}{R^2}\int_{T} \abs{\chi}^2 \diff{x} \bigg].
  \end{align*}
  We may now apply inverse estimates to obtain
  \begin{align*}
    \MoveEqLeft \sup_{x \in T} \zeta_h^2(x) \int_T \abs{ \nabla \chi}^2 \diff{x} \lesssim \\
    &\sup_{x \in T} \zeta_h^2(x) \int_T \abs{\nabla(\tiluh-p_0-\PiSZ w)}^2 \diff{x} + \int_T \abs{\nabla(\PiSZ w-w)}^2 \diff{x} \\
    &\lesssim \int_T \zeta_h^2 \abs{\nabla(\tiluh-p_0-\PiSZ w)}^2 \diff{x}
    + \int_T \abs{\nabla(\PiSZ w-w)}^2 \diff{x} \\
    &\lesssim \int_T \zeta_h^2 \abs{\nabla \chi}^2 \diff{x}
    + \int_T \abs{\nabla(\PiSZ w-w)}^2 \diff{x}.
  \end{align*}
  Using interpolation estimates and gathering all the previously obtained estimates therefore gives us, after summation over $T$
  \begin{align*}
    \int_{B^\Omega_R} \nabla\chi \cdot \nabla(\psi-\psi_h) \diff{x} &\lesssim
    \epsilon \int_{B_R^\Omega} \zeta_h^2 \abs{\nabla \chi}^2 \diff{x} 
    + c_\epsilon \frac1{R^2}\int_{B_R^\Omega} \abs{\chi}^2 \diff{x} \\
    & + 
    c_\epsilon h^2 \sum_{ \substack{ T \in \mathcal{T}_h \\ T \cap \support \zeta_h \neq \emptyset}}
    \int_{\omega(T)} \abs{\nabla^2 w}^2 \diff{x}.
  \end{align*}
  Applying Lemma~\ref{lem:properties-w} and summing over all $T$ finishes the proof.
\end{proof}

Our last preliminary estimate deals with the lower order term.

\begin{lemma}[Lower order term]
\label{lem:lower-order-estimate}
In the established setting, for every $\xi \in \RRn$, we have
  \begin{align*}
    \frac1{R^2} \int_{B^\Omega_R} \abs{\tiluh-p_0-w}^2 \diff{x} \lesssim \frac{h}{R} \int_{B^\Omega_R} \abs{\nabla \tiluh -\xi}^2 \diff{x} + \int_{B^\Omega_R} \abs{\nabla (\tilu-\barp)}^2 \diff{x}.
  \end{align*}
\end{lemma}
\begin{proof}
  Fix $\xi \in \RRn$. We proceed by duality. The dual problem is then to find $z \in \setW^{1,2}_0(B^\Omega_R)$ such that
  \begin{equation}
  \label{eq:DefOfDualProblem}
    \int_{B_R^\Omega} \nabla \phi \cdot \nabla z \diff{x} = \int_{B_R^\Omega} ( \tiluh - p_0 - w ) \phi \diff{x}, \quad \forall \phi \in \setW^{1,2}_0(B^\Omega_R).
  \end{equation}
  We will also need a finite element approximation of $z$ that also has zero boundary values on $B^\Omega_R$, i.e., its support must be contained in $B_R^\Omega$. To achieve this, we define a cut-off function $\tau \in \setC^2(B^\Omega_R)$ that satisfies
  \[
    \indicator_{B_{R-2h}} \leq \tau \leq \indicator_{B_{R-h}},
    \qquad \qquad
    \norm{\nabla^i\tau}_{\setL^\infty(\Omega)} \lesssim \frac1{h^i}, \qquad i=1,2.
  \]
  The support properties of $\tau$ guarantee that $\support \PiSZ(\tau z) \subseteq B^\Omega_R$. We define
  \[
    z_h \coloneqq \PiSZ(\tau z) + \PiP(z-\PiSZ(\tau z)),
  \]
  where $\PiP$ is defined in Lemma~\ref{lem:poly-proj-def}.
  If necessary, we increase $\Lambda$ here once more.
  
  Our next step is to obtain some regularity estimates on $z$ and $z_h$. Since $B_R^\Omega$ is convex, by standard regularity theory, we have that $z \in \setW^{2,2}(B^\Omega_R)$ and
  \begin{equation}
    \label{eq:dual-problem-stability}
    \frac1{R} \norm{\nabla z}_{\setL^2(B^\Omega_R)} + \norm{\nabla^2 z}_{\setL^2(B^\Omega_R)}
    \lesssim \norm{\tiluh-p_0-w}_{\setL^2(B^\Omega_R)}.
  \end{equation}
  
  Define the annulus $A_R \coloneqq B^\Omega_R \setminus B_{R-2h}$ and observe that, by definition, $\support \nabla \tau  \subset A_R$. Next, using a trace inequality,
  we estimate
  \[
    \norm{\nabla z}_{\setL^2(A_R)}^2
    = \int_0^{2h} \int_{\partial B_{R-t} \cap \Omega} \abs{\nabla z}^2 \diff{s} \diff{t}
    \lesssim h \bigg( \frac{1}{R} \norm{\nabla z}^2_{\setL^2(B^\Omega_R)} + R \norm{\nabla^2 z}^2_{\setL^2(B^\Omega_R)}\bigg).
  \]
  Notice that, while the first step is similar to using polar coordinates, there is no transformation determinant as the ball surfaces are still scaled to size.
  Taking square roots and applying \eqref{eq:dual-problem-stability} leads to
  \begin{equation}
    \label{eq:z-grad-annulus}
    \norm{\nabla z}_{\setL^2(A_R)}
    \lesssim (hR)^{1/2} \norm{\tiluh-p_0-w}_{\setL^2(B^\Omega_R)}.
  \end{equation}
  This last estimate, together with the $\setW^{1,2}(B_R^\Omega)$-stability of $\PiP$, approximation properties of the Scott--Zhang interpolant, and the definition of $\tau$, implies
  \begin{align*}
    \norm{\nabla(z-z_h)}_{\setL^2(B^\Omega_R)}^2
    &\lesssim \norm{\nabla(z-\PiSZ(\tau z))}_{\setL^2(B^\Omega_R)}^2 \\
    &\lesssim \norm{\nabla(\tau z-\PiSZ(\tau z))}_{\setL^2(B^\Omega_R)}^2 + \norm{\nabla((1-\tau)z)}_{\setL^2(B^\Omega_R)}^2 \\
    &\lesssim h^2 \norm{\nabla^2(\tau z)}_{\setL^2(B^\Omega_R)}^2 + \norm{\nabla z}_{\setL^2(A_R)}^2 + \frac1{h^2} \norm{z}_{\setL^2(A_R)}^2 \\
    &\lesssim h^2 \norm{\nabla^2 z}_{\setL^2(B^\Omega_R)}^2 + \norm{\nabla z}_{\setL^2(A_R)}^2 + \frac1{h^2} \norm{z}_{\setL^2(A_R)}^2.
  \end{align*}
  Since $z = 0$ on $\partial B^\Omega_R$, \Poincare's inequality, \eqref{eq:dual-problem-stability}, and \eqref{eq:z-grad-annulus} give us
  \begin{equation}
  \label{eq:z-approximation}
    \norm{\nabla(z-z_h)}_{\setL^2(B^\Omega_R)}^2
    \lesssim \Big( h^2 + hR \Big) \norm{\tiluh-p_0-w}_{\setL^2(B^\Omega_R)}^2.
  \end{equation}
  Finally, by the stability in $\setW^{1,2}(B_R^\Omega)$ of $\PiP$, we have
  \begin{align}
    \label{eq:zh-stability}
    \begin{aligned}
      \norm{\nabla z_h}_{\setL^2(B^\Omega_R)}
      &\lesssim \norm{\nabla \PiSZ(\tau z)}_{\setL^2(B^\Omega_R)} + \norm{\nabla z}_{\setL^2(B^\Omega_R)} \\
      &\lesssim \frac1{h} \norm{z}_{\setL^2(A_R)} + \norm{\nabla z}_{\setL^2(B^\Omega_R)} \\
      &\lesssim \norm{\nabla z}_{\setL^2(B^\Omega_R)}
      \lesssim R \norm{\tiluh-p_0-w}_{\setL^2(B^\Omega_R)}.
    \end{aligned}
  \end{align}
  
  Having established all these bounds we can estimate the lower order term by duality. Setting $\phi = \tiluh - p_0 - w \in \setW^{1,2}_0(B_R^\Omega)$ in \eqref{eq:DefOfDualProblem} gives
  \begin{multline*}
    \int_{B^\Omega_R} \abs{\tiluh-p_0-w}^2 \diff{x}
    = \int_{B^\Omega_R} \nabla(\tiluh-p_0-w) \cdot \nabla z \diff{x} \\
    = \int_{B^\Omega_R} \nabla(\tiluh-p_0-w) \cdot \nabla(z-z_h) \diff{x} + \int_{B^\Omega_R} \nabla(\tiluh-p_0-w)\cdot \nabla z_h \diff{x}
    \eqqcolon \mathrm{I} + \mathrm{II}.
  \end{multline*}
  For the term $\mathrm{I}$, notice $z-z_h = (z-\PiSZ(\tau z)) - \PiP(z-\PiSZ(\tau z))$. Therefore, we use that $p_H$ is harmonic, the defining property of $\PiP$, Young's inequality, and \eqref{eq:z-approximation}, to obtain
  \begin{align*}
    \mathrm{I}
    &= \int_{B^\Omega_R} \nabla(\tiluh-p_0-w)\cdot \nabla(z-z_h) \diff{x}
    = \int_{B^\Omega_R} \nabla(\tiluh-\barp-w)\cdot \nabla(z-z_h) \diff{x} \\
    &=\int_{B^\Omega_R} \nabla(\tiluh-w) \cdot \nabla(z-z_h) \diff{x} \\
    &\leq \delta R^2 \int_{B^\Omega_R} \abs{\nabla(\tiluh-w)}^2 \diff{x} + \frac1{4\delta R^2} \int_{B^\Omega_R} \abs{\nabla(z-z_h)}^2 \diff{x} \\
    &\lesssim \delta R^2 \int_{B^\Omega_R} \abs{\nabla(\tiluh-w)}^2 \diff{x} + \frac1{\delta}\Big( \Big(\frac{h}{R} \Big)^2 +\frac{h}{R} \Big) \int_{B^\Omega_R} \abs{\tiluh-p_0-w}^2 \diff{x}.
  \end{align*}
  We choose $\delta = \frac1{\epsilon} \Big( \Big(\frac{h}{R} \Big)^2 +\frac{h}{R} \Big)$, and using Lemma \ref{lem:properties-w} on the first term after adding $\pm \xi$ results in
  \begin{align*}
    \mathrm{I}
    &\lesssim \frac1{\epsilon} R^2\Big( \Big(\frac{h}{R} \Big)^2 +\frac{h}{R} \Big) \int_{B^\Omega_R} \abs{\nabla \tiluh -\xi}^2 \diff{x} + \epsilon \int_{B^\Omega_R} \abs{\tiluh-p_0-w}^2 \diff{x}.
  \end{align*}
  For $\epsilon$ small enough, the last term can be absorbed on the left hand side.
  
  For the term $\mathrm{II}$, we use that $w$ and $p_H$ are harmonic as well as the equation \eqref{eq:tiluh-def-equation}, Young's inequality, and \eqref{eq:zh-stability} to estimate
  \begin{align*}
    \mathrm{II}
    & = \int_{B^\Omega_R} \nabla(\tiluh-p_0-w)\cdot \nabla z_h \diff{x}
    = \int_{B^\Omega_R} \nabla(\tilu-p_0-p_H)\cdot \nabla z_h \diff{x} \\
    &\lesssim c_\epsilon R^2 \int_{B^\Omega_R} \abs{\nabla(\tilu-\barp)}^2 \diff{x} + \frac{\epsilon}{R^2} \int_{B^\Omega_R} \abs{\nabla z_h}^2 \diff{x} \\
    &\lesssim c_\epsilon R^2 \int_{B^\Omega_R} \abs{\nabla(\tilu-\barp)}^2 \diff{x} + \epsilon \int_{B^\Omega_R} \abs{\tiluh-p_0-w}^2 \diff{x}.
  \end{align*}
  Choosing $\epsilon$ small enough allows absorption of the last term on the left hand side. This finishes the proof.
\end{proof}

Finally, we are able to compare $\tiluh$ and $w$.

\begin{lemma}[Comparison]
\label{lem:comparison}
With the previous definitions, for every $\xi \in \RRn$, we have
  \[
    \int_{B^\Omega_R} \zeta_h^2 \abs{\nabla (\tiluh-p_0-w)}^2\diff{x}
    \lesssim \frac{h}{R} \int_{B^\Omega_R} \abs{\nabla \tiluh -\xi }^2\diff{x}
    + \int_{B^\Omega_R} \abs{\nabla (\tilu - \barp)}^2\diff{x}.
  \]
\end{lemma}
\begin{proof}
  Denote by $G$ the left hand side of the claim. Next, recall that, owing to \eqref{eq:def-psi},
  \[
    \nabla \psi = \zeta_h^2\nabla(\tiluh-p_0-w)+2\zeta_h\nabla\zeta_h(\tiluh-p_0-w).
  \]
  This allows us to write
  \begin{align*}
    G &=\int_{B^\Omega_R} \zeta_h^2\nabla(\tiluh-p_0-w)\cdot \nabla(\tiluh-p_0-w) \diff{x} \\
    &= \int_{B^\Omega_R} \nabla\psi \cdot \nabla(\tiluh-p_0-w) \diff{x} - 2\int_{B^\Omega_R} \nabla(\tiluh-p_0-w)\cdot \zeta_h \nabla\zeta_h(\tiluh-p_0-w) \diff{x} \\
    &= \int_{B^\Omega_R} \nabla(\psi-\psi_h)\cdot \nabla(\tiluh-p_0-w) \diff{x} + \int_{B^\Omega_R} \nabla \psi_h \cdot \nabla(\tiluh-p_0-w) \diff{x} \\
    &- 2\int_{B^\Omega_R} \nabla(\tiluh-p_0-w) \cdot \zeta_h\nabla\zeta_h(\tiluh-p_0-w) \diff{x} \\
    &\eqqcolon \mathrm{I}+\mathrm{II}+\mathrm{III}.
  \end{align*}
  We now estimate the terms $\mathrm{I}$, $\mathrm{II}$, and $\mathrm{III}$ separately. Fix $\xi \in \RRn$. A suitable bound of the term $\mathrm{I}$ was carried out in Lemma~\ref{lem:mixed-term}. Choosing, in that result, $\epsilon >0$ but sufficiently small, we may conclude that
  \[
    G \lesssim  \frac{h^2}{R^2} \int_{B^\Omega_R} \abs{\nabla \tiluh -\xi}^2 \diff{x}
    + \frac{1}{R^2} \int_{B^\Omega_R} \abs{\tiluh-p_0-w}^2 \diff{x}
    + \mathrm{II} + \mathrm{III}.
  \]

  We now proceed with the term $\mathrm{II}$. We use Lemma~\ref{lem:grad-psi-estimate} to obtain
  \begin{align*}
    \int_{B_R^\Omega}\abs{ \nabla\psi_h}^2 \diff x &\lesssim \int_{B_R^\Omega}\abs{ \nabla(\psi-\psi_h) }^2 \diff x + \int_{B_R^\Omega}\abs{ \nabla\psi}^2 \diff x \\
    &\lesssim G
    + \frac{h^2}{R^2} \int_{B^\Omega_R} \abs{\nabla \tiluh -\xi}^2 \diff{x}
    + \frac{1}{R^2} \int_{B^\Omega_R} \abs{\tiluh-p_0-w}^2 \diff{x}.
  \end{align*}
  Using the fact $w$ and $p_H$ are harmonic in $B_R^\Omega$ and using the equation for $\tiluh$ we get
  \begin{align*}
    \mathrm{II}
    &= \int_{B^\Omega_R} \nabla \psi_h \cdot \nabla(\tiluh-p_0-w) \diff{x}
    = \int_{B^\Omega_R} \nabla \psi_h \cdot \nabla (\tiluh-\barp) \diff{x}\\
    &= \int_{B^\Omega_R} \nabla(\tilu-\barp)\cdot \nabla\psi_h \diff{x}
    \leq \epsilon \int_{B_R^\Omega} \abs{\nabla \psi_h}^2 \diff x +  c_\epsilon\int_{B^\Omega_R} \abs{\nabla(\tilu-\barp)}^2 \diff{x}.
  \end{align*}
  The previously obtained estimate and an application of Young's inequality lead to
  \[
    G \lesssim  \frac{h^2}{R^2} \int_{B^\Omega_R} \abs{\nabla \tiluh -\xi}^2 \diff{x}
    + \frac{1}{R^2} \int_{B^\Omega_R} \abs{\tiluh-p_0-w}^2 \diff{x}
    + \int_{B^\Omega_R} \abs{\nabla (\tilu - \barp)}^2 \diff{x}
    + \mathrm{III}.
  \]

  For the third term, we recall that $\norm{\nabla\zeta_h}_{\setL^\infty(B_R^\Omega)} \lesssim R^{-1}$ as well as, again, Young's inequality to get
  \[
    \mathrm{III} \leq \epsilon G + \frac{c_\epsilon}{R^2} \int_{B^\Omega_R} \abs{\tiluh-p_0-w}^2 \diff{x}.
  \]
  A suitable choice of $\epsilon$ then implies
  \[
    G \lesssim \frac{h^2}{R^2} \int_{B^\Omega_R} \abs{\nabla \tiluh -\xi }^2 \diff{x}
    + \frac{1}{R^2} \int_{B^\Omega_R} \abs{\tiluh-p_0-w}^2 \diff{x}
    + \int_{B^\Omega_R} \abs{\nabla (\tilu - \barp)}^2 \diff{x}.
  \]
  Finally, we invoke Lemma~\ref{lem:lower-order-estimate} to conclude. This finishes the proof.
\end{proof}

\subsection{Proof of oscillation reduction estimate}

We now have all the ingredients to prove the main results of this section, namely, the oscillation estimates in Lemma~\ref{lem:campanato}.

\begin{proof}[Proof of Lemma~\ref{lem:campanato}]
  We prove each statement separately.

  \begin{enumerate}[label={\alph{*})},leftmargin=*]
    \item For the case of interior balls we begin by recalling that it suffices to prove \eqref{eq:campanato-tilde-estimate}. Let
    \[
      q_w \in \argmin_{q \in \mathcal{P}_m}  \dashint_{B^\Omega_\rho} \abs{\nabla(w-q)}^2 \diff{x}.
    \]
    The choice $q = q_w + p_0$ and the triangle inequality give us
    \begin{align*}
      \inf_{q \in \mathcal{P}_m} \dashint_{B^\Omega_\rho} \abs{\nabla (\tiluh-q)}^2 \diff{x}
      & \lesssim \dashint_{B^\Omega_\rho} \abs{\nabla (w - q_w)}^2 \diff{x}
      + \dashint_{B^\Omega_\rho} \abs{\nabla(\tiluh-p_0-w)}^2 \diff{x} \\
      & \eqqcolon \mathrm{I} + \mathrm{II}.
    \end{align*}
    We invoke Lemma~\ref{lem:properties-w} item \ref{itm:w-interior-regularity} to estimate
    \[
      \mathrm{I}
      = \dashint_{B^\Omega_\rho} \abs{\nabla (w - q_w)}^2 \diff{x}
      \lesssim \frac{\rho^{2m}}{R^{2m}}\dashint_{B^\Omega_R} \abs{\nabla w}^2 \diff{x}
      \lesssim \frac{\rho^{2m}}{R^{2m}}\dashint_{B^\Omega_R} \abs{\nabla \tiluh}^2 \diff{x}.
    \]
    For the term $\mathrm{II}$, let us recall the function $\zeta_h \in \mathcal{L}_1^1(\mathcal{T}_h)$ introduced in \eqref{eq:DefOfzetahCutOff}. Since $\rho < \frac12 R$ we may write
    \begin{align*}
      \mathrm{II} &= \dashint_{B^\Omega_\rho} \abs{\nabla(\tiluh-p_0-w)}^2 \diff{x} = \dashint_{B^\Omega_\rho} \zeta_h^2 \abs{\nabla(\tiluh-p_0-w)}^2 \diff{x}
      \\
      &\lesssim \frac{R^n}{\rho^n} \dashint_{B^\Omega_R} \zeta_h^2 \abs{\nabla(\tiluh-p_0-w)}^2 \diff{x}.
    \end{align*}
    We may then continue this estimate by invoking Lemma~\ref{lem:comparison} with $\xi = 0$ to obtain
    \begin{align*}
      \mathrm{II}
      &\lesssim \frac{h}{R} \frac{R^n}{\rho^n}\dashint_{B^\Omega_R} \abs{\nabla \tiluh}^2\diff{x}
      + \frac{R^n}{\rho^n}\dashint_{B^\Omega_R} \abs{\nabla(\tilu-\barp)}^2\diff{x}.
    \end{align*}
    Combining all the estimates therefore implies \eqref{eq:campanato-tilde-estimate}. This finishes the case of interior balls.

    \item For balls that intersect the boundary we must make a few adjustments. First recall that, in this case, $\tiluh \coloneqq R_h u$ and $\tilu \coloneqq u$. Next, we fix $m=1$ so that $\barp \in \mathcal{P}_1$ is itself harmonic. This implies $p_H= \barp$ and $p_0 = 0$. In a similar manner to the previous step we choose
    \[
      q_w \in \argmin_{q \in \mathcal{P}_1} \dashint_{B^\Omega_\rho} \abs{\nabla(w-q_w)}^2 \diff{x}.
    \]
    In view of Lemma~\ref{lem:inf_equiv} (which is trivial for $m=1$) we get
    \begin{align*}
      \inf_{\xi \in \RRn} \dashint_{B^\Omega_\rho} \abs{\nabla \tiluh-\xi}^2 \diff{x}
      &=\inf_{q \in \mathcal{P}_1} \dashint_{B^\Omega_\rho} \abs{\nabla (\tiluh-q)}^2 \diff{x} \\
      & \lesssim \dashint_{B^\Omega_\rho} \abs{\nabla (w - q_w)}^2 \diff{x}
      + \dashint_{B^\Omega_\rho} \abs{\nabla(\tiluh-p_0-w)}^2 \diff{x}.
    \end{align*} 
    For the first term we have that $w = \tiluh = R_h u = 0$ on $\partial\Omega \cap B_\rho^\Omega$, and therefore, with Lemma \ref{lem:properties-w} item \ref{itm:w-boundary-regularity}, we get
    \begin{align*}
      \dashint_{B^\Omega_\rho} \abs{\nabla (w - q_w)}^2 \diff{x}
      &\lesssim \frac{\rho^{2\sigma}}{R^{2\sigma}} \inf_{\xi \in \RRn}\dashint_{B^\Omega_R} \abs{\nabla w - \xi}^2 \diff{x}
      \lesssim \frac{\rho^{2\sigma}}{R^{2\sigma}} \inf_{\xi \in \RRn} \dashint_{B^\Omega_R} \abs{\nabla \tiluh -\xi}^2 \diff{x}.
    \end{align*}
    The second term like before is handled with Lemma~\ref{lem:comparison} to estimate
    \begin{align*}
      \dashint_{B^\Omega_\rho} \abs{\nabla(\tiluh-p_0-w)}^2 \diff{x}
      &\lesssim \rho^{-n} \int_{B^\Omega_R} \zeta_h^2 \abs{\nabla(\tiluh-p_0-w)}^2 \diff{x} \\
      &\lesssim \frac{hR^n}{R\rho^n} \inf_{\xi \in \RRn}\dashint_{B^\Omega_R} \abs{\nabla \tiluh -\xi}^2\diff{x}
      + \frac{R^n}{\rho^n}\dashint_{B^\Omega_R} \abs{\nabla(\tilu-\barp)}^2\diff{x}.
    \end{align*}
    Again with Lemma~\ref{lem:inf_equiv}, combining all estimates finishes the proof for balls that intersect the boundary.
  \end{enumerate}
\end{proof}

\section{Pointwise estimates and applications}
\label{sub:pointwise-estimate}

In this section the oscillation reduction estimate of Theorem~\ref{thm:osc-estimate} is used to establish a pointwise estimate similar to \eqref{eq:PtWiseEstimate}, which we previously obtained under the assumption that $n \leq 3$; see \cite{MR4748207}. This allows to obtain stability of the Ritz projection in a wide range of Banach function spaces, but now without the dimensional restriction.

\subsection{Maximal functions}

To provide the needed pointwise estimates we must begin by recalling some maximal functions, their relations, and properties. For $\tau \in [1,\infty)$ we define the Hardy--Littlewood maximal function of $v \in \setL^\tau_\loc(\Omega)$ via
\begin{equation}
\label{eq:DefOfMtau}
  M_\tau[v](x) \coloneqq \sup_{r>0} \bigg( \dashint_{B^\Omega_r(x)} |v|^\tau \diff y \bigg)^{1/\tau}, \qquad a.e.~x \in \Omega.
\end{equation}
As usual, we write $M \coloneqq M_1$.
Notice that, by Jensen's inequality, if $\tau_1 \leq \tau_2$, we have the pointwise bound
\[
  M_{\tau_1}[v](x) \leq M_{\tau_2}[v](x).
\]
In addition, it is well known that for $\rho \in (\tau,\infty]$ the maximal function $M_\tau[\cdot]$ is bounded in $\setL^\rho(\Omega)$, that is,
\begin{equation}
\label{eq:HLmaxBddLrho}
  \norm{ M_\tau[v] }_{\setL^\rho(\Omega)} \lesssim \norm{v}_{L^\rho(\Omega)}, \qquad \forall v \in \setL^\rho(\Omega).
\end{equation}
The proof of this result for the case $\tau=1$ can be found in many references; see, for instance, \cite[Theorem 2.1.6]{MR3243734}. For $\tau > 1$, notice $(M_\tau[v])^\tau = M[v^\tau]$. Then the proof follows from the $\tau=1$ case.
Next we define, for $\tau \in [1,\infty)$, the sharp maximal function of $v \in \setL^\tau_\loc(\Omega)$ via
\[
  M^{\sharp}_{\tau} [v](x) \coloneqq \sup_{r > 0} \osckp{1}{\tau} v(x,r), \qquad a.e.~x \in \Omega.
\]
Again, we write $M^\sharp \coloneqq M^\sharp_1$. The sharp maximal function also obeys the pointwise bound
\begin{equation}
\label{eq:MsharptauMonotoneInTau}
  M_{\tau_1}^\sharp[v](x) \leq M_{\tau_2}^\sharp[v](x), \qquad \tau_1 \leq \tau_2,
\end{equation}
and it is controlled by the Hardy--Littlewood maximal via
\begin{equation}
\label{eq:MsharptauBddByM}
  M^{\sharp}_{\tau} [v](x) \leq 2 M_\tau[v](x), \qquad a.e.~x \in \Omega.
\end{equation}
From this, the boundedness of $M_\tau^\sharp[\cdot]$ in $\setL^\rho(\Omega)$, for $\rho \in (\tau,\infty]$, immediately follows. The reverse pointwise estimate is clearly not valid (set $v\equiv 1$). However, a local version of the Fefferman--Stein inequality asserts that, for $\rho \in (\tau,\infty)$,
\begin{equation}
\label{eq:LrhoFeffermanStein}
  \biggnorm{ v - \dashint_\Omega v \diff y }_{\setL^\rho(\Omega)} \lesssim \norm{ M_\tau^\sharp[v] }_{\setL^\rho(\Omega)}.
\end{equation}
See, for instance, \cite[Theorem 5.23]{MR2643399} for a proof in the case $\tau = 1$, from which the general case here follows by the monotonicity \eqref{eq:MsharptauMonotoneInTau}.

\subsection{\texorpdfstring{Pointwise estimate and stability in $\setBMO$ and $\setL^\rho$}{Pointwise estimate and stability in BMO and Lp}}
\label{sec:BMO-and-lp-spaces}

Our first pointwise estimate now follows. We immediately note that, in Theorem~\ref{thm:pointwise-sharp-maximal-ii}, we will provide an alternative version of this result.

\begin{theorem}[Pointwise estimate]
  \label{thm:pointwise-sharp-maximal}
  Let $u \in \setW^{1,2}_0(\Omega)$ and $R_h u \in V_h$ be defined via~\eqref{eq:RitzProjection}. For almost every $x \in \Omega$ we have
  \begin{align}
    M^{\sharp}_{2} [\nabla R_h u](x) \lesssim M^{\sharp}_{2} [\nabla u](x).
  \end{align}
\end{theorem}
\begin{proof}
  In view of Theorem~\ref{thm:osc-estimate}, fix $\beta\in (0,\sigma)$ and obtain $\theta \in (0,1)$. There holds
  \begin{align*}
    \sup_{r> \theta\Lambda h} \osckp{1}{2} \nabla R_h u(x,r) \leq c_\beta \sup_{r > \Lambda h} \osckp{1}{2} \nabla u(x,r) + \theta^\beta \sup_{r > \Lambda h} \osckp{1}{2} \nabla R_h u(x,r).
  \end{align*}
  The last term on the right hand side can be absorbed on the left to get
  \[
    \sup_{r> \theta\Lambda h} \osckp{1}{2} \nabla R_h u(x,r) \lesssim M^{\sharp}_{2} [\nabla u](x).
  \]
  It remains to bound the terms $r \leq \theta\Lambda h$. Those balls only contain a uniformly bounded number of simplices, so by standard arguments for finite dimensional spaces one can conclude that, for $r \leq \theta\Lambda h$,
  \[
    \osckp{1}{2} \nabla R_h u(x,r) \lesssim \osckp{1}{2} \nabla R_h u(x,\theta\Lambda h).
   \qedhere
  \]
\end{proof}

The rest of this section is dedicated to obtaining applications of this result.

As a first application of Theorem~\ref{thm:pointwise-sharp-maximal} we obtain stability in $\setBMO(\Omega)$ and $\setL^\rho(\Omega)$, $\rho \in (1,\infty)$. We reiterate that the novelty in our results is not only in the approach and higher polynomial degree, but in that we are not placing any dimensional restriction.

For completeness we recall that $\setBMO(\Omega)$ can be equivalently characterized via, see \cite[Section 5.8]{MR3024912},
\begin{align*}
  \setBMO(\Omega) &\coloneqq \Big\{ v\in \setL^2(\Omega) \ : \ M_2^\sharp[v] \in \setL^\infty(\Omega) \Big\}, \\
  \norm{v}_{\setBMO(\Omega)} &\coloneqq \norm{v}_{\setL^2(\Omega)} + \norm{ M_2^\sharp[v] }_{\setL^\infty(\Omega)}.
\end{align*}

\begin{corollary}[$\setBMO$-stability]
  \label{cor:bmo-stab}
  Let $u \in \setW^{1,1}_0(\Omega)$ be such that $\nabla u \in \setBMO(\Omega)$. There holds
  \begin{align}
    \label{eq:BMO-stab}
    \norm{\nabla R_h u}_{\setBMO(\Omega)} \lesssim \norm{\nabla u}_{\setBMO(\Omega)}.
  \end{align}
\end{corollary}
\begin{proof}
  With Theorem~\ref{thm:pointwise-sharp-maximal} and the trivial $\setW^{1,2}$-boundedness of the Ritz projection 
  we estimate
  \begin{align*}
    \norm{\nabla R_h u}_{\setBMO(\Omega)}
    &= \norm{\nabla R_h u}_{\setL^2(\Omega)} +  \norm{M^\sharp_2[\nabla R_h u]}_{\setL^\infty(\Omega)} \\
    &\lesssim \norm{\nabla u}_{\setL^2(\Omega) }+ \norm{M^\sharp_2[\nabla u]}_{\setL^\infty(\Omega)}
    = \norm{\nabla u}_{\setBMO(\Omega)},
  \end{align*}
  as claimed.
\end{proof}

We now prove $\setL^\rho$-stability.

\begin{corollary}[$\setL^\rho$-stability]
  \label{cor:lp-stab}
  Let $\rho \in (1,\infty)$, $u \in \setW^{1,\rho}_0(\Omega)$. There holds
  \begin{align}
    \label{eq:lp-stab}
    \norm{\nabla R_h u}_{\setL^\rho(\Omega)} \lesssim \norm{\nabla u}_{\setL^\rho(\Omega)}.
  \end{align}
\end{corollary}
\begin{proof}
The case $\rho =2$ is trivial by definition of the Ritz projection.

In the case $\rho > 2$ we begin by invoking \eqref{eq:LrhoFeffermanStein} with $\tau = 2$ for $\nabla R_h u$, that is,
\[
    \norm{\nabla R_h u}_{\setL^\rho(\Omega)}
    \lesssim \norm{M^{\sharp}_{2} [\nabla R_h u]}_{\setL^\rho(\Omega)}.
\]
Notice that $R_h u = 0$ on $\partial\Omega$ implies that $\nabla R_h u$ has zero average. We continue with Theorem~\ref{thm:pointwise-sharp-maximal}; then we apply the pointwise bound \eqref{eq:MsharptauBddByM}; and, finally,
\eqref{eq:HLmaxBddLrho} with $\tau = 2$ and $\rho >2$. In summary, we have
\begin{align*}
    \norm{\nabla R_h u}_{\setL^\rho(\Omega)}
    &\!\lesssim\! \norm{M^{\sharp}_{2} [\nabla R_h u]}_{\setL^\rho(\Omega)}
    \!\lesssim \! \norm{M^{\sharp}_{2} [\nabla u]}_{\setL^\rho(\Omega)}
    \!\lesssim \! \norm{M_{2} [\nabla u]}_{\setL^\rho(\Omega)}
     \!\lesssim \! \norm{\nabla u}_{\setL^\rho(\Omega)}.
\end{align*}

For the remaining case, that is $\rho \in (1,2)$, we argue by duality.
We have, see \cite[Theorem 1]{MR1198129} and \cite{MR1156467}, owing to the fact that $\Omega$ is convex that
\[
  \norm{ \nabla v }_{\setL^\rho(\Omega)} \eqsim \sup_{\psi \in \setW^{1,\rho'}_0(\Omega)} \frac{ \int_\Omega \nabla v \cdot \nabla \psi \diff x}{ \norm{ \nabla\psi }_{\setL^{\rho'}(\Omega) } }, \qquad \forall v \in \setW^{1,\rho}_0(\Omega).
\]
For smooth domains see~\cite[Theorem 7.2]{MR473503}.
We apply this to $R_h u$ to get
  \begin{align*}
    \norm{\nabla R_h u}_{\setL^\rho(\Omega)}
    &\eqsim \sup_{\psi \in \setW_0^{1,\rho'}(\Omega)} \frac{\int_\Omega \nabla R_h u \cdot \nabla \psi \diff{x}}{\norm{\nabla \psi}_{\setL^{\rho'}(\Omega)}}
    = \sup_{\psi \in \setW_0^{1,\rho'}(\Omega)} \frac{\int_\Omega \nabla u \cdot \nabla R_h \psi \diff x}{\norm{\nabla \psi}_{\setL^{\rho'}(\Omega)}} \\
    &\lesssim \norm{\nabla u}_{\setL^\rho(\Omega)} \sup_{\psi \in \setW_0^{1,\rho'}(\Omega)} \frac{\norm{\nabla R_h \psi}_{\setL^{\rho'}(\Omega)}}{\norm{\nabla \psi}_{\setL^{\rho'}(\Omega)}}.
  \end{align*}
  Since $\rho' > 2$, the supremum is bounded by the previous step, giving the desired statement.
\end{proof}

\subsection{Reducing the oscillation exponent}
\label{sec:reducing-osc-exp}

Theorem~\ref{thm:osc-estimate} as well as the pointwise estimate of Theorem~\ref{thm:pointwise-sharp-maximal} are formulated for oscillations with integrability exponent equal to $2$. Using the stability results of Corollary~\ref{cor:lp-stab}, we are able to reduce the needed exponent to any $\tau \in (1,2)$, which will also strengthen the pointwise estimate of Theorem~\ref{thm:pointwise-sharp-maximal}. We begin with an auxiliary result.

\begin{lemma}[Oscillation exponent reduction]
  \label{lem:reduce-osc-exp-aux}
  Let $\beta>0$, $\eta \in (0,1)$, $\Lambda >1$, and $m \in \setN$. Fix $B_R \subset \Omega$ with $R = \diameter B_R \geq \Lambda h$. Assume $z_h \in V_h$ is such that, for every $x \in B_{2R}$ and all $\rho \in (\Lambda h, R)$ such that $B_\rho(x) \subset B_{2R}$,
  \begin{align}
    \label{eq:ass-reduction}
  \osckp{m}{2} \nabla z_h(x,\eta \rho) \leq \eta^\beta \osckp{m}{2} \nabla z_h(x,\rho).
  \end{align}
  Then, for every $\tau \in [1,2)$ and $\theta \in (0,1)$, there holds
  \begin{align*}
  \osckp{m}{\tau} \nabla z_h(x,\theta R) \lesssim \theta^\beta \osckp{m}{\tau} \nabla z_h(x,2R).
  \end{align*}
\end{lemma}
\begin{proof}
 Let $p \in (\mathcal{P}_{m-1})^n$ be arbitrary. We define $f = \nabla z_h - p$.
 By definition of $f$, from~\eqref{eq:ass-reduction} and basic calculations, for $\Lambda h \leq \rho \leq R$, we get
 \begin{align*}
   \osckp{m}{2} f(x,\eta \rho) &\leq \eta^\beta \osckp{m}{2} f(x,\rho), \\
   \Bigg(\dashint_{B_{\eta\rho}(x)} \abs{f}^2 \diff{z}\Bigg)^{1/2} &\lesssim c_\eta \osckp{m}{2} f(x,\rho) +
 \Bigg(\dashint_{B_{\rho}(x)} \abs{f}^2 \diff{z}\Bigg)^{1/2}.
 \end{align*}
 Iteration therefore gives that for $\Lambda h \leq \rho < r < R$ there holds
 \begin{align}
   \label{eq:f-iteration}
   \begin{aligned}
 \Bigg(\dashint_{B_{\rho}(x)} \abs{f}^2 \diff{z}\Bigg)^{1/2}
 &\lesssim \sum_{j} c_\eta \eta^{j\beta} \osckp{m}{2} f(x,r) + 
 \Bigg(\dashint_{B_{r}(x)} \abs{f}^2 \diff{z}\Bigg)^{1/2} \\
 &\lesssim \Bigg(\dashint_{B_{r}(x)} \abs{f}^2 \diff{z}\Bigg)^{1/2},
   \end{aligned}
 \end{align}
 where we use that this sum is a convergent geometric series.

  For the next step, for $0 < t < s < 2R$ with $s-t > \Lambda h$, we will prove
  \begin{align}
    \label{eq:Linfty-L2-est}
    \norm{f}_{\setL^\infty(B_t(x))}
    &\lesssim
    (s-t)^{-n/2}\norm{f}_{\setL^2(B_s(x))}.
  \end{align}
  To see this, fix some $y \in B_t(x)$. By inverse estimates, there holds
  \[
    \abs{f(y)} \lesssim \bigg(\dashint_{B_{\Lambda h}(y)} \abs{f}^2 \diff{z}\bigg)^{1/2}.
  \]
  Furthermore, from \eqref{eq:f-iteration} we get
  \begin{align*}
    \bigg(\dashint_{B_{\Lambda h}(y)} \abs{f}^2 \diff{z}\bigg)^{1/2}
    &\lesssim
  \bigg(\dashint_{B_{s-t}(y)} \abs{f}^2 \diff{z}\bigg)^{1/2}.
  \end{align*}
  Since $B_{s-t}(y) \subset B_s(x)$, putting this together gives
  \[
    |f(y)| \lesssim (s-t)^{-n/2} \norm{f}_{\setL^2(B_s(x))},
  \]
  proving~\eqref{eq:Linfty-L2-est}.

Now we define
\[
  \gamma \coloneqq n\bigg( \frac1\tau - \frac12 \bigg), \quad Z(t) \coloneqq \norm{f}_{\setL^2(B_t(x))}, \quad A \coloneqq \norm{f}_{\setL^\tau(B_{2R}(x))}.
\]
By interpolation and equation~\eqref{eq:Linfty-L2-est}, for $R < t < s < 3R/2$, with $s-t > \Lambda h$ there holds
\begin{align*}
  Z(t) \lesssim (s-t)^{-\gamma\tau/2}Z(s)^{1-\tau/2}A^{\tau/2}.
\end{align*}
With Young's inequality, this leads to
\begin{align}
  \label{eq:giusti-iter}
  Z(t) \lesssim \epsilon Z(s) + c_\epsilon (s-t)^{-\gamma}A.
\end{align}
For $s-t \leq \Lambda h$, on the other hand, by simple reverse inequalities there holds
\begin{align*}
  Z(t) \lesssim h^{-\gamma} A \lesssim (s-t)^{-\gamma}A.
\end{align*}
Therefore,~\eqref{eq:giusti-iter} is actually valid for all $R < t < s < 3R/2$. By the iteration argument of~\cite[Lemma~6.1]{MR1962933} we can therefore conclude
  \begin{align*}
    Z(R) \lesssim R^{-\gamma}A.
  \end{align*}
  With Jensen's inequality and iterative application of~\eqref{eq:ass-reduction}, this leads to
  \begin{align*}
  \osckp{m}{\tau} \nabla z_h(x,\theta R)
  &\lesssim \osckp{m}{2} \nabla z_h(x,\theta R)
  \lesssim
  \theta^\beta \osckp{m}{2} \nabla z_h(x,R) \\
  &\leq \theta^\beta R^{-n/2} \norm{f}_{\setL^2(B_R(x))}
  \lesssim \theta^\beta R^{-n/\tau} \norm{f}_{\setL^\tau(B_{2R}(x))} \\
  &= \theta^\beta R^{-n/\tau} \norm{\nabla z_h - p}_{\setL^\tau(B_{2R}(x))}.
  \end{align*}
  Since $p$ was arbitrary in the beginning, we may choose it as the best approximation of $\nabla z_h$ in $\setL^\tau(B_{2R}(x))$.
\end{proof}

Equipped with this auxiliary result we can formulate a direct improvement over Theorem~\ref{thm:osc-estimate}.

\begin{theorem}[Oscillation reduction]
  \label{thm:osc-estimate-ii}
  Let $\tau \in (1,2)$ and $\beta \in (0,k)$. There are $\theta \in (0,1)$ and  $\Lambda >1$ such that, for $u \in \setW^{1,\tau}_0(\Omega)$, the following estimates hold.
  \begin{enumerate}[label={\alph{*})},leftmargin=*]
    \item If $x \in \Omega$ and $R > \Lambda h$ are such that $B_R(x) \Subset \Omega$, then
    \[
      \oscp{\tau} \nabla R_h u(x,\theta R) \leq \theta^\beta \oscp{\tau} \nabla R_h u (x,R) + c_{\beta,\tau} \oscp{\tau} \nabla u (x,R).
    \]
    
    \item There exists $\sigma \in (0,1]$ depending on the geometry of $\Omega$ such that provided $\beta \in (0,\sigma)$, for every $x \in \bar\Omega$ and all $R > \Lambda h$, there holds
    \[
      \osckp{1}{\tau} \nabla R_h u(x,\theta R) \leq \theta^\beta \osckp{1}{\tau} \nabla R_h u (x,R) + c_{\beta,\tau} \osckp{1}{\tau} \nabla u (x,R).
    \]
  \end{enumerate}
  In these estimates the constants denoted $c_{\beta,\tau}$ blow up as $\tau \downarrow 1$.
\end{theorem}
\begin{proof}
  We only consider the interior case to simplify the presentation. The modifications needed to obtain the boundary case are straightforward. Furthermore, we assume $u \in \setW^{1,2}_0(\Omega)$; the general case follows by density. Let $u_h \coloneqq R_h u$.

  Fix a ball $B_R \coloneqq B_R(x)$ such that $B_{2R} \subset \Omega$.
  Let $p \in \mathcal{P}_{k}$. Define
  \begin{align*}
    F \coloneqq \indicator_{B_{2R}^\complement} \nabla u + \indicator_{B_{2R}}\nabla p.
  \end{align*}
  Let $z \in \setW^{1,2}_0(\Omega)$ verify $-\Delta z = -\divergence F$, and let $z_h =R_h z$ be its Ritz projection.
  Next we observe that, while Theorem~\ref{thm:osc-estimate} is formulated for $\nabla z$ and $\nabla z_h$, its conclusion remains valid when replacing $\nabla z$ with $\nabla p$. That is because everything that is needed is $\int_{B_{2R}} \nabla z_h \cdot \nabla \phi \diff{y} =\int_{B_{2R}} \nabla z \cdot \nabla \phi \diff{y} $ for all $\phi \in \setW^{1,2}_0(B_{2R})$. Therefore, we may replace $\nabla z$ by $F = \nabla p$ on $B_{2R}$.

  The previous considerations show that, for $\gamma \in (\beta, k)$, $\rho \in (\Lambda h,R)$, and every $y \in B_R(x)$ there exists $\eta \in (0,1)$ such that
  \begin{align*}
    \oscp{2} \nabla z_h(y,\eta \rho) &\leq \eta^\gamma \oscp{2} \nabla z_h (y,\rho) + c_\gamma \oscp{2} \nabla p (y,\rho) \\
    &= \eta^\gamma \oscp{2} \nabla z_h (y,\rho),
  \end{align*}
  where the oscillation of $\nabla p$ vanishes since it is a polynomial.
  This allows the application of Lemma~\ref{lem:reduce-osc-exp-aux}. By the triangle inequality, we get
  \begin{align*}
    \oscp{\tau} \nabla u_h(x,\theta R)
    &\leq
    \oscp{\tau} \nabla z_h(x,\theta R) + (\theta R)^{-n/\tau}\norm{\nabla(u_h-z_h)}_{\setL^\tau(\Omega)} \\
    &\leq C_{\tau,\gamma}
    \theta^\gamma\oscp{\tau} \nabla z_h(x,2 R) + (\theta R)^{-n/\tau}\norm{\nabla(u_h-z_h)}_{\setL^\tau(\Omega)} \\
    &\leq C_{\tau,\gamma}
    \theta^\gamma\oscp{\tau} \nabla u_h(x,2 R) + (c_\theta R)^{-n/\tau}\norm{\nabla(u_h-z_h)}_{\setL^\tau(\Omega)}.
  \end{align*}
  We may now choose $\theta$ sufficiently small so that $C_{\tau,\gamma}\theta^\gamma \leq \theta^\beta$ and, if necessary, further increase $\Lambda$. It remains to deal with the last term. By Corollary~\ref{cor:lp-stab} as well as the fact that $\Omega$ is convex (see again \cite[Theorem 1]{MR1198129} and \cite{MR1156467}) there holds
  \begin{align*}
    R^{-n/\tau}\norm{\nabla(u_h-z_h)}_{\setL^\tau(\Omega)}
    &\lesssim
    R^{-n/\tau}\norm{\nabla(u-z)}_{\setL^\tau(\Omega)} \\
    &\lesssim
    R^{-n/\tau}\norm{\nabla u- F}_{\setL^\tau(\Omega)} \\
    &\lesssim
    R^{-n/\tau}\norm{\indicator_{B_{2R}} (\nabla u - \nabla p)}_{\setL^\tau(\Omega)}.
  \end{align*}
  Together with Lemma~\ref{lem:inf_equiv}, all that is left is choosing $p$ as the best approximation in $\setL^\tau(B_{2R})$ and rescaling from $B_{2R}$ to $B_R$.
\end{proof}

As a direct consequence we are also able to improve Theorem~\ref{thm:pointwise-sharp-maximal}. 

\begin{theorem}[Pointwise estimate II]
  \label{thm:pointwise-sharp-maximal-ii}
  Let $\tau \in (1,2)$, $u \in \setW^{1,\tau}_0(\Omega)$, and $R_h u \in V_h$ be defined via~\eqref{eq:RitzProjection}. For almost every $x \in \Omega$ there holds
  \begin{align}
    M^{\sharp}_{\tau} [\nabla R_h u](x) \lesssim M^{\sharp}_{\tau} [\nabla u](x).
  \end{align}
\end{theorem}
\begin{proof}
  The proof repeats that of Theorem~\ref{thm:pointwise-sharp-maximal} by using the improved oscillation reduction of Theorem~\ref{thm:osc-estimate-ii}.
\end{proof}

\subsection{Banach function spaces}
\label{sec:banach}

We now generalize the stability of the Ritz projection to a very large class of Banach function spaces in the sense of \cite[Definition 2.1]{Lerner2017}. The essence of the argument in the proof of Corollary~\ref{cor:lp-stab} is an application of the Fefferman--Stein inequality \eqref{eq:LrhoFeffermanStein}, Theorem~\ref{thm:pointwise-sharp-maximal}, the ideal property, and boundedness of the maximal function \eqref{eq:HLmaxBddLrho}. Thus, the stability of the Ritz projection can be generalized to any Banach function space that verifies these properties. We axiomatize these ideas in the following definition.

\begin{definition}[Maximal Fefferman--Stein]
\label{def:PropFrakRFuncSpace}
  Let $\setX(\Omega)$ be a Banach function space on $\Omega$. We say that this space satisfies the maximal Fefferman--Stein property if:
  \begin{enumerate}[label={\alph{*})},leftmargin=*]
    \item The Hardy--Littlewood maximal function is bounded, that is,
    \begin{equation}
    \label{eq:AbstractHLbound}
      \norm{ M[v] }_{\setX(\Omega)} \lesssim \norm{ v}_{\setX(\Omega)}, \qquad \forall v \in \setX(\Omega).
    \end{equation}

    \item The following version of the Fefferman--Stein inequality holds
    \begin{equation}
    \label{eq:AbstractFeffermanStein}
    \biggnorm{ v - \dashint_\Omega v \diff{x} }_{\setX(\Omega)} \lesssim \norm{ M^\sharp[v] }_{\setX(\Omega)}, \qquad \forall v \in \setX(\Omega).
    \end{equation}
  \end{enumerate}
\end{definition}

The following Theorem is formulated for bounded John domains, see~\cite{MR2643399}. Note that every bounded polytopal domain is a John domain. If $\setX(\RRn)$ is a Banach function space, then $\setX(\Omega)$ is just its restriction to $\Omega$.

\begin{theorem}[Dual maximal boundedness]
  \label{thm:whole-space-restriction}
  Assume $\setX(\RRn)$ is a Banach function space such that $M$ is bounded on $\setX(\RRn)$ and $\setX'(\RRn)$, where $\setX'(\RRn)$ is the K\"othe dual of $\setX(\RRn)$. Furthermore let $\Omega \subset \RRn$ be a bounded John domain. Then, $\setX(\Omega)$ satisfies the maximal Fefferman--Stein property of Definition~\ref{def:PropFrakRFuncSpace}.
\end{theorem}
\begin{proof}
  Property~\eqref{eq:AbstractHLbound} follows directly by restriction. It remains to prove the Fefferman--Stein type inequality \eqref{eq:AbstractFeffermanStein}.
  A proof of~\eqref{eq:AbstractFeffermanStein} for weighted Lebesgue spaces on John domains can be found in~\cite{MR2643399}. Rubio de Francia extrapolation (see~\cite[Theorem 10.1]{CruzUribe2017} and~\cite{CruzUribeMartellPerez2011}) then finishes the proof.
\end{proof}

\begin{remark}[Connections]
  By the results of \cite{MR2763004}, the dual maximal boundedness of Theorem~\ref{thm:whole-space-restriction} is equivalent to $M$ being bounded on $\setX(\RRn)$ and the Fefferman--Stein inequality holding on $\setX(\RRn)$. Furthermore, $M$ being bounded on $\setX(\RRn)$ and $\setX(\RRn)$ being $s$-concave is a stronger condition and implies the dual maximal boundedness; see \cite{lerner2026dual}. Essentially, the $s$-concavity ensures $\setX(\RRn)$ is not too close to $\setL^\infty(\RRn)$. The weak Lebesgue spaces $\setL^{\rho,\infty}(\RRn)$, for $1 < \rho < \infty$, are an example of a non $s$-concave space that still satisfies the dual maximal boundedness.
\end{remark}

For Banach function spaces with the maximal Fefferman--Stein property the stability of the Ritz projection immediately follows.

\begin{theorem}[$\setX$-stability]
\label{cor:StabRitzBanachFunctionSpaces}
  Let $\setX(\Omega)$ be a Banach function space that satisfies the maximal Fefferman--Stein property, and let $u \in \setW^{1,1}_0(\Omega)$ be such that $\nabla u \in \setX(\Omega)$. Then we have
  \[
    \norm{ \nabla R_h u }_{\setX(\Omega)} \lesssim \norm{ \nabla u }_{\setX(\Omega)},
  \]
  with a constant that is independent of $h$. Moreover,
  \[
    \norm{ \nabla (u - R_h u) }_{\setX(\Omega)} \lesssim \inf_{v_h \in V_h} \norm{ \nabla (u -v_h) }_{\setX(\Omega)}.
  \]
\end{theorem}
\begin{proof}
  As announced, the proof abstracts that of Corollary~\ref{cor:lp-stab}.
  According to \cite[Corollary 1.3]{MR2375698}, from \eqref{eq:AbstractHLbound} it follows the existence of a $\tau > 1$ for which $M_\tau$ is bounded on $\setX(\Omega)$.
  By Jensen's inequality, $M^\sharp[\nabla R_h u] \leq M^\sharp_\tau[\nabla R_h u]$. Also there holds $\int_\Omega \nabla R_h u \diff{x} = 0$.
  We invoke \eqref{eq:AbstractFeffermanStein}, the ideal property, Theorem~\ref{thm:pointwise-sharp-maximal} or \ref{thm:pointwise-sharp-maximal-ii}, the ideal property again, and finally the mentioned boundedness of $M_\tau$ to obtain
  \begin{equation}
    \label{eq:fefferman-stein-hardy-little}
    \begin{aligned}
      \norm{\nabla R_h u}_{\setX(\Omega)}
      &\lesssim \norm{M^\sharp[\nabla R_h u]}_{\setX(\Omega)} 
      \lesssim \norm{M^\sharp_{\tau}[\nabla R_h u]}_{\setX(\Omega)} 
      \lesssim \norm{M^\sharp_{\tau}[\nabla u]}_{\setX(\Omega)} \\
      &\lesssim \norm{M_{\tau}[\nabla u]}_{\setX(\Omega)}
      \lesssim \norm{\nabla u}_{\setX(\Omega)}.
    \end{aligned}
  \end{equation}
  This proves the first part. For the second, we invoke the triangle inequality, the projection property of $R_h$ and the first inequality to get
  \begin{align*}
    \norm{\nabla (u-R_h u)}_{\setX(\Omega)}
    &\leq
    \inf_{v_h \in V_h}
    \big(\norm{\nabla (u-v_h)}_{\setX(\Omega)} + \norm{\nabla R_h (u-v_h)}_{\setX(\Omega)}\big) \\
    &\lesssim \inf_{v_h \in V_h} \norm{ \nabla (u -v_h) }_{\setX(\Omega)}.
    \qedhere
  \end{align*}
\end{proof}

\subsubsection{Muckenhoupt weighted spaces}%
\label{sec:muckenhoupt}

As an example of Banach function spaces that satisfy the maximal Fefferman--Stein property we obtain stability of the Ritz projection for Muckenhoupt weighted spaces. This extends our result from~\cite{MR4748207} to all dimensions.

We recall that a function $0 \leq \omega \in \setL^1_\loc(\RRn)$ that is almost everywhere positive is called a weight. For $\rho \in (1,\infty)$ we say that a weight $\omega$ belongs to the Muckenhoupt class $\mathcal{A}_\rho$ if
\[
  [\omega]_{\mathcal{A}_\rho} \coloneqq
    \sup_{B} \bigg( \dashint_B \omega(x) \diff{x} \bigg)\bigg( \dashint_B \omega^{-\tfrac1{\rho-1}}(x) \diff{x} \bigg)^{\rho-1},
\]
is finite. Here the supremum is over all balls $B$. Weighted Lebesgue spaces are defined, for $\rho \in (1,\infty)$ and $\omega \in \mathcal{A}_\rho$, as
\begin{align*}
  \| f \|_{\setL^\rho_\omega(\Omega)} &\coloneqq \bigg[ \int_\Omega \abs{f(x)}^\rho \omega(x) \diff{x} \bigg]^{1/\rho}, \\
  \setL^\rho_\omega(\Omega) &\coloneqq \big\{ f \in \setL^0(\Omega) : \| f \|_{\setL_\omega^\rho(\Omega)}< \infty \big\}.
\end{align*}

Let us now show that the maximal Fefferman--Stein property is satisfied by Muckenhoupt weighted spaces and that, as a consequence, the Ritz projection is stable. It is well known that $M$ is bounded on $\setL^\rho_\omega(\Omega)$ and its dual; see, for instance, \cite[Theorem 7.1.9]{MR3243734}. Therefore, we have the dual maximal-boundedness of Theorem~\ref{thm:whole-space-restriction}, which implies the maximal Fefferman--Stein property of Definition~\ref{def:PropFrakRFuncSpace}. Theorem~\ref{cor:StabRitzBanachFunctionSpaces} gives stability of the Ritz projection.

\begin{corollary}[Weighted $\setL^\rho$-spaces]
  Let $\rho \in (1,\infty)$ and $\omega \in \mathcal{A}_\rho$. Then $\setL^\rho_\omega(\Omega)$ satisfies the maximal Fefferman--Stein property of Definition~\ref{def:PropFrakRFuncSpace}. As a consequence, for $u \in \setW^{1,1}_0(\Omega)$ such that $\nabla u\in \setL^\rho_\omega(\Omega)$, we have
  \begin{align*}
    \norm{\nabla R_h u}_{\setL^\rho_\omega(\Omega)} \lesssim \norm{\nabla u}_{\setL^\rho_\omega(\Omega)}.
  \end{align*}
\end{corollary}

\begin{remark}[Left-openness]
  In the case of Muckenhoupt weighted spaces, the existence of $\tau$ from the proof of Theorem~\ref{cor:StabRitzBanachFunctionSpaces} can be made qualitative. In particular, $\tau$ only depends on $n, \rho$, and the Muckenhoupt constant $[\omega]_{\mathcal{A}_\rho}$. This is a consequence of the left-openness of the Muckenhoupt class; see \cite[Corollary 7.2.6]{MR3243734}.
\end{remark}

\subsubsection{Other function spaces}

As the maximal function is a classical operator, it is well studied for a wide range of function spaces $\setX(\RRn)$. Together with the $\setX(\RRn)$-stability of Theorem~\ref{cor:StabRitzBanachFunctionSpaces}, this immediately implies the stability of the Ritz projection. We give a short, non-exhaustive list of examples.

\begin{enumerate}[label={$\bullet$},leftmargin=*]
  \item \textbf{Lorentz and Lorentz--Zygmund spaces:} If $1 < \rho < \infty$, $1 \leq q \leq \infty$ and $\alpha \in \RR$, then the maximal operator is bounded on $\setL^{\rho,q}(\RRn)$ and $\setL^{\rho,q}(\log \setL)^\alpha(\RRn)$ and their K\"othe dual, respectively; see~\cite[Theorems A, section 5.3]{MR1156767} as well as~\cite{BennettRudnick1980}.

  \item \textbf{Orlicz spaces:} Under the $\Delta_2$- and $\nabla_2$-conditions on $\phi$, the maximal operator is bounded on $\setL^\phi(\RRn)$ and its K\"othe dual; see~\cite[Theorem 1.2.1(v)]{MR1156767}. 

  \item \textbf{Variable exponent spaces:} If $\rho(\cdot)$ is $\log$-H\"older continuous and $1 < \rho^- \leq \rho^+ < \infty$, the maximal operator is bounded on $\setL^{\rho(\cdot)}(\RRn)$ and its K\"othe dual; we refer to~\cite{MR2790542,MR3026953}. The $\log$-H\"older continuity can be replaced by the Muckenhoupt condition 
  $\norm{\indicator_B}_{\rho(\cdot)} \norm{\indicator_B}_{\rho'(\cdot)} \eqsim \abs{B}$
  and the Nekvinda decay condition; see~\cite[Theorem 4.52]{MR3026953} and~\cite{MR4951202}.
  
  \item \textbf{Weighted variable exponent spaces:} If $\rho(\cdot)$ is $\log$-Hölder continuous and $1 < \rho^- \leq \rho^+ < \infty$, and additionally $\omega \in \mathcal{A}_{\rho(\cdot)}$, the maximal operator is bounded on $\setL_\omega^{\rho(\cdot)}(\RRn)$ and its K\"othe dual; we refer to~\cite{MR2837636} and \cite{MR2790542}.

  \item \textbf{Weighted Orlicz spaces:} If $\phi$ is a Young function with Simonenko indices $1<p \leq q < \infty$ and $\omega \in \mathcal{A}_p$, then the maximal operator is bounded on $\setL^\phi_\omega(\RRn)$ and its K\"othe dual, see~\cite{KermanTorchinsky1982}.

  \item \textbf{Mixed-norm Lebesgue spaces:} Let $\boldsymbol{\rho}=(\rho_1,\ldots,\rho_n)$ with $1<\rho_j<\infty$ for every $j$. The maximal operator is bounded on $\setL^{\boldsymbol{\rho}}(\RRn)$ and its K\"othe dual $\setL^{\boldsymbol{\rho}'}(\RRn)$; see~\cite[Theorem 4.1]{MR909030}, and \cite[Theorem 1]{MR2316449} for the weighted case.

  \item \textbf{Rearrangement-invariant spaces:} If $\setX(\RRn)$ is a rearrangement-invariant Banach function space \cite[Definition 2.4.1]{BennettSharpley1988} with nontrivial Boyd indices \cite[Definition 3.5.12]{BennettSharpley1988}; see also \cite[Chapter II]{MR506343}, then the maximal operator is bounded on $\setX(\RRn)$ and its K\"othe dual; see~\cite[Theorem 3.5.17]{BennettSharpley1988}.
\end{enumerate}

\section{Discrete smoothness spaces}
\label{sec:finite-frequency}

In this section, we develop discrete analogues of Besov and Triebel--Lizorkin spaces. Our main motivation is to measure the smoothness of finite element functions and their gradients at scales down to the mesh size $h$. Indeed, when $s\rho>1$, we have $\mathcal{L}^0_m(\mathcal{T}_h) \not\subset \setW^{s,\rho}(\Omega)$, due to jumps across simplex faces. The optimal regularity of those functions in terms of Besov spaces is $\setB^{1/2}_{2,\infty}(\Omega)$; see~\cite{MR4238777}.

Nevertheless, we will show that finite element solutions possess higher-order smoothness in a suitable discrete sense. Moreover, local mollifications and other smooth approximations of these solutions belong to standard higher-order smoothness spaces, with bounds uniform in $h$.

Related constructions of discrete smoothness spaces appear in different settings. Oswald \cite{Oswald1990,Oswald1994} introduced Besov-type approximation spaces associated with nested finite element spaces in the context of multilevel methods.  Martin and Perkowski \cite{MartinPerkowski2019} developed discrete Besov spaces on lattices using a Littlewood--Paley decomposition truncated at the lattice scale. Our approach is based on local oscillations and does not resolve scales below the mesh size.

Mesh-dependent norms that measure differentiability beyond the global Sobolev regularity of the discrete functions have a long tradition in discontinuous Galerkin and interior penalty methods, see~\cite{Arnold1982,ArnoldBrezziCockburnMarini2002}. For continuous piecewise polynomial functions measured in $\setW^{2,p}(\Omega)$ via jumps of the derivatives, see~\cite{BrennerSung2005,FengHenningsNeilan2017}.  It is also possible to measure discrete $\setW^{2,2}$-regularity by the norm of the Riesz representative of the weak Laplacian in the finite element space, using either the standard $\setL^2$-inner product~\cite[Section~2.2]{MR2461254} or a mass-lumped counterpart~\cite[Section~3]{BartelsDolzmannNochetto2010}. Estimates controlling elementwise Laplacians and jumps of normal derivatives by weighted norms of the discrete Laplacian and the gradient are established on convex polyhedral domains in~\cite[Lemma~2.6]{LeykekhmanVexler2016}.  Discrete Bessel potential spaces on uniform grids are defined in~\cite[Section~2.5.3]{MR3136501} by means of the Fourier transform.

The results of this section hold on arbitrary regular domains $\Omega$. In particular, convexity is not required.

\subsection{Measuring smoothness via oscillations}
\label{sub:Brudnyi}

We begin by recalling the definition of the oscillation introduced in \eqref{eq:LocOscillation}.
Notice that, owing to \eqref{eq:BiteOffEstimate}, for every $x \in \overline{D}$ and $t>0$ we have
\[
  \mathrm{osc}_{m,p}^Dv(x,t) \lesssim \mathrm{osc}_{m,p}v(x,t).
\]
We will focus on oscillations with $p=1$, and so we set
\begin{align*}
  \mathrm{osc}_{m} v(\cdot,\cdot) \coloneqq \mathrm{osc}_{m,1} v(\cdot,\cdot), \qquad \mathrm{osc}_{m}^D v(\cdot,\cdot) \coloneqq \mathrm{osc}_{m,1}^D v(\cdot,\cdot).
\end{align*}
With these at hand we can give an equivalent characterization of Besov and Triebel--Lizorkin spaces. For smooth domains these definitions are rather classical, see \cite{MR290001,MR1023118}. For Lipschitz domains, we refer the reader to \cite{MR2017700,MR1984176,MR4801010} for a discussion of this equivalence. For regular domains, see~\cite{MR2247585}.

\begin{definition}[Besov and Triebel--Lizorkin spaces]
\label{def:besovTriebel}
  Let $D \subseteq \Omega$ be a regular domain with $R = \diameter D$.
  Assume that $q,\rho \in [1,\infty]$ and $s \in (0,k)$.
  Set
  \[
    |g|_{\setB^s_{\rho,q}(D)} \coloneqq
    \begin{dcases}
    \bigg( \int_0^R t^{-sq} \bignorm{ \osc^D g(\cdot, t) }_{\setL^\rho(D)}^q \frac{\diff{t}}t \bigg)^{1/q}, &  q<\infty,
      \\
    \sup_{t \in (0,R)} t^{-s} \bignorm{ \osc^D g(\cdot, t) }_{\setL^\rho(D)},  & q = \infty.
  \end{dcases}
  \]
  The Besov spaces are defined as
  \[
    \setB_{\rho,q}^s(D) \coloneqq \Big\{ g \in \setL^\rho(D) \ : \ \abs{g}_{\setB_{\rho,q}^s(D)} < \infty \Big\},
    \quad
    \norm{g}_{\setB^s_{\rho,q}(D)} \coloneqq \norm{g}_{\setL^\rho(D)} + |g|_{\setB^s_{\rho,q}(D)}.
  \]

  Assume now that the indices are such that $\rho \in [1,\infty)$, $q \in [1,\infty]$, and $s \in (0,k)$.
  Set
  \[
    \abs{g}_{\setF_{\rho,q}^s(D)} \coloneqq
    \begin{dcases}
      \biggnorm{ \bigg( \int_0^R \Big( \frac{\osc^D g(\cdot,t)}{t^s} \Big)^q \frac{\diff{t} }t \bigg)^{1/q}  }_{\setL^\rho(D)}, & q<\infty,
      \\
      \biggnorm{ \sup_{t \in (0,R)} t^{-s} \osc^D g(\cdot,t) }_{\setL^\rho(D)}, & q = \infty.
    \end{dcases}
  \]
  The Triebel--Lizorkin spaces are
  \[
    \setF_{\rho,q}^s(D) \coloneqq \Big\{ g \in \setL^\rho(D) \ : \ \abs{g}_{\setF^s_{\rho,q}(D)} < \infty \Big\},
    \quad
    \norm{g}_{\setF^s_{\rho,q}(D)} \coloneqq \norm{g}_{\setL^\rho(D)} + |g|_{\setF^s_{\rho,q}(D)} .
  \]
\end{definition}

The interested reader is referred to the aforementioned literature for details.
In particular it is shown that using, in these definitions, oscillations of different orders gives rise to the same space and equivalent norms, as long as the order of the oscillation is larger than $s$. In addition, provided $r \in [1,\infty)$ satisfies
\begin{equation}
\label{eq:BesovIntoL2}
  n\bigg( \frac1{\rho}-\frac1r \bigg)_+ < s,
\end{equation}
for Besov spaces, and
\begin{equation}
\label{eq:TriebelIntoL2}
  n\bigg( \frac1{\rho}-\frac1r \bigg)_+ < s, \qquad n\bigg( \frac1{q}-\frac1r\bigg)_+ < s,
\end{equation}
for Triebel--Lizorkin spaces, we may equivalently use $\oscp{r}^D$ in the definition of the seminorm.

\begin{remark}[Conciseness]
  As it has become tradition when dealing with these classes of function spaces, we will use $\setA_{\rho,q}^s(D)$, with $\setA \in \{\setB,\setF\}$, to state most of our results, with the understanding that in every statement the meaning of $\setA$ is fixed. In addition, to avoid repetition in the proofs, we will only provide details for the case of Triebel--Lizorkin spaces and $\rho,q < \infty$. Unless stated otherwise, the Besov space case will follow by a straightforward modification.
\end{remark}

From the definition of the seminorms it is clear that, if $\varepsilon>0$,
\[
 \setA^{s+\varepsilon}_{\rho,q}(\Omega) \hookrightarrow \setA^s_{\rho,q}(\Omega),
 \qquad
 \setA^{s}_{\rho+\varepsilon,q}(\Omega) \hookrightarrow \setA^s_{\rho,q}(\Omega),
 \qquad
 \setA^{s}_{\rho,q}(\Omega) \hookrightarrow \setA^s_{\rho,q+\varepsilon}(\Omega),
\]
provided all the indices are in the allowed range. Further properties of these scales of spaces are collected below.

\begin{remark}[Embeddings]
\label{rem:sobolev-embeddings}
  Let $s\in (0,k)$, $\rho_1,\rho_2 \in [1,\infty]$, $q_1,q_2 \in [1,\infty]$. Then
  \begin{equation}
  \label{eq:BintoFintoB}
    \setA^s_{\rho_2,q_2}(\Omega) \hookrightarrow \setA^s_{\rho_1,q_1}(\Omega),
  \end{equation}
  if and only if $q_2 \leq q_1$ and $\rho_2 \geq \rho_1$. Together with $\setB^s_{\rho,\rho}(\Omega) = \setF^s_{\rho,\rho}(\Omega)$ this covers the case of same differentiability.

  For $i = 1,2$ let $\rho_i,q_i \in [1,\infty]$ and $s_i \in (0,k)$. Assume, in addition, that $s_1 < s_2$. The continuous embedding
  \begin{equation}
  \label{eq:SobEmbeddingWithAs}
    \setA^{s_2}_{\rho_2,q_2}(\Omega) \hookrightarrow \setA^{s_1}_{\rho_1,q_1}(\Omega),
  \end{equation}
  holds under the following conditions:
  \begin{enumerate}[label={\alph{*})},leftmargin=*]
    \item (Subcritical case.) We have
    \[
      s_1 - \frac{n}{\rho_1} < s_2 - \frac{n}{\rho_2}.
    \]
    \item (Critical case.) We have
    \[
      s_1 - \frac{n}{\rho_1} = s_2 - \frac{n}{\rho_2},
    \]
    for $\setA = \setB$ we additionally require $q_1 \geq q_2$.
  \end{enumerate}
  In all those cases, if $\setA=\setF$, we must also have $\rho_i < \infty$.

  The reader is referred to \cite[Theorem 3.3.1]{MR781540} for a proof. In addition, one may combine \eqref{eq:BintoFintoB} with \cite[Theorem 1.97]{MR2250142} to assert that, in the subcritical case, the embedding in \eqref{eq:SobEmbeddingWithAs} is compact.
\end{remark}

\begin{remark}[Some standard identifications]
\label{rem:StandardReductions}
  Note, first of all, that by Fubini's theorem we have
  \[
    \setF^s_{\rho,\rho}(\Omega) = \setB^s_{\rho,\rho}(\Omega).
  \]
  In addition, owing to \cite[(4.20)]{MR2250142}, if $\rho \in (1,\infty)$ and $m \in \setN$, we have $\setF^m_{\rho,2}(\Omega) = \setW^{m,\rho}(\Omega)$, where on the right we have the classical Sobolev spaces of weakly differentiable functions. Finally, \cite[Theorem 6.7]{MR1152321} shows that if $s \notin \setN$
  \[
    \setB^s_{\rho,\rho}(\Omega) = \setW^{s,\rho}(\Omega),
  \]
  where $\setW^{s,\rho}(\Omega)$ denotes the fractional Sobolev spaces: let $s = m + \nu$, $m \in \setN_0$ and $\nu \in (0,1)$; then
  \begin{align*}
    \setW^{s,\rho}(\Omega) &\coloneqq \Big\{ g \in \setW^{m,\rho}(\Omega) \ : \ | g |_{\setW^{s,\rho}(\Omega)} < \infty \Big\}, \\
    |g|_{\setW^{s,\rho}(\Omega)} &\coloneqq \Bigg( \sum_{|j|=m} \iint_{\Omega \times \Omega} \frac{ | \partial^j g(x) - \partial^j g(y) |^\rho }{|x-y|^{n + \nu \rho} } \diff x \diff y
    \Bigg)^{1/\rho}.
  \end{align*}
\end{remark}

\begin{remark}[The case $s=0$]
  \label{rem:triebel-lizorkin-zero}
  Our treatment of Besov and Triebel--Lizorkin spaces via oscillations does not include the case $s=0$. However, the space $\setF^0_{\rho,q}(\Omega)$ for a bounded Lipschitz domain $\Omega$ can be defined by other means, like the Fourier transform combined with restriction. In particular, if $\rho \in (1,\infty)$, then $\setF^0_{\rho,2}(\Omega) = \setL^\rho(\Omega)$. If $\rho = \infty$ then
  \[
    \setF_{\infty,2}^0(\Omega) = \setBMO(\Omega);
  \]
  see~\cite[Sections 1.7.1 and 1.7.2]{TriebelFunctionSpacesII}.
\end{remark}

\begin{assumption}[Indices]
  For the rest of our work we will operate under the assumptions on $\rho$, $q$, $s$ that are stated in Definition~\ref{def:besovTriebel}.
\end{assumption}

\subsection{Discrete Besov and Triebel--Lizorkin spaces}

The discrete version of our smoothness spaces is defined next. From here on, $\Omega$ is a regular polytopal domain. For the most part, $h > 0$ is a free parameter, with the intention of it becoming the mesh size of the triangulation $\mathcal{T}_h$.

\begin{definition}[Discrete Besov and Triebel--Lizorkin spaces]
\label{def:disc-besovTriebel}
Let $D \subseteq \Omega$ be a regular subdomain with $R = \diameter D$, and let $h \in (0,\tfrac{R}{4})$.
  Under the conditions of Definition~\ref{def:besovTriebel} define the seminorms
  \[
    |g|_{\setB^s_{\rho,q,h}(D)} \coloneqq
    \begin{dcases}
    \bigg( \int_h^R t^{-sq} \bignorm{ \osc^D g(\cdot, t) }_{\setL^\rho(D)}^q \frac{\diff{t}}t \bigg)^{1/q}, & q<\infty,
      \\
    \sup_{t \in (h,R)} t^{-s} \bignorm{ \osc^D g(\cdot, t) }_{\setL^\rho(D)}, & q = \infty,
    \end{dcases}
  \]
  and
  \[
    \abs{g}_{\setF_{\rho,q,h}^s(D)} \coloneqq
    \begin{dcases}
      \biggnorm{ \bigg( \int_h^R \Big( \frac{\osc^D g(\cdot,t)}{t^s} \Big)^q \frac{\diff{t} }t \bigg)^{1/q}  }_{\setL^\rho(D)}, & q<\infty,
      \\
      \biggnorm{ \sup_{t \in (h,R)} t^{-s} \osc^D g(\cdot,t) }_{\setL^\rho(D)}, & q = \infty.
    \end{dcases}
  \]
  The discrete Besov and Triebel--Lizorkin spaces are, respectively,
  \begin{align*}
    \setB_{\rho,q,h}^s(D) &\coloneqq \Big\{ g \in \setL^\rho(D) \ : \ \abs{g}_{\setB_{\rho,q,h}^s(D)} < \infty \Big\},
    \\
    \norm{g}_{\setB^s_{\rho,q,h}(D)} &\coloneqq \norm{g}_{\setL^\rho(D)} + |g|_{\setB^s_{\rho,q,h}(D)},
    \\
    \setF_{\rho,q,h}^s(D) &\coloneqq \Big\{ g \in \setL^\rho(D) \ : \ \abs{g}_{\setF^s_{\rho,q,h}(D)} < \infty \Big\},
    \\
    \norm{g}_{\setF^s_{\rho,q,h}(D)} &\coloneqq \norm{g}_{\setL^\rho(D)} + |g|_{\setF^s_{\rho,q,h}(D)} .
  \end{align*}
\end{definition}

Note that, as in Remark~\ref{rem:StandardReductions}, we have
\[
  \setF^s_{\rho,\rho,h}(D) = \setB^s_{\rho,\rho,h}(D).
\]
In addition, for any $C>0$ such that $Ch < \tfrac{R}{4}$,
\[
  \setA^s_{\rho,q,h}(D) = \setA^s_{\rho,q,Ch}(D),
\]
with an equivalence constant that depends on $C$ but not on $h$. In what follows we will repeatedly use this equivalence without explicit mention.

\subsection{\texorpdfstring{Uniform embeddings and compactness for $s\rho <1$}{Uniform embeddings and compactness for sp<1}}

We next wish to study some of the properties of the scales of spaces $\{\setA^s_{\rho,q,h}(\Omega)\}_{h > 0}$. First of all, it is clear that
\begin{equation}
\label{eq:ObviousEmbedding}
  \norm{g}_{\setA^s_{\rho,q,h}(\Omega)} \leq \norm{g}_{\setA^s_{\rho,q}(\Omega)}.
\end{equation}
Next, since in the definition of the $\setA^s_{\rho,q,h}(\Omega)$ seminorm the oscillations are measured down to scale $h>0$, we must elucidate what happens to these norms as $h \searrow 0$. The following results are in the spirit of those in the seminal work \cite{MR3586796}; see also \cite{MR4654593} for the vectorial case.

\begin{lemma}[Monotone convergence]
\label{lem:norm-convergence-h0}
  Let $g \in \setL^\rho(\Omega)$ be such that
  \[
    \sup_{h > 0} \norm{g}_{\setA^s_{\rho,q,h}(\Omega)} < \infty.
  \]
  Then $g \in \setA^s_{\rho,q}(\Omega)$ with
  \[
    \norm{g}_{\setA^s_{\rho,q}(\Omega)}
     = \lim_{h \searrow 0} \norm{g}_{\setA^s_{\rho,q,h}(\Omega)}.
  \]
\end{lemma}
\begin{proof}
  Notice that $\norm{g}_{\setA^s_{\rho,q}(\Omega)} = \norm{g}_{\setA^s_{\rho,q,0}(\Omega)}$. Observe also that, for $h_1 \leq h_2$, we have
  \begin{align*}
    \norm{g}_{\setA^s_{\rho,q,h_1}(\Omega)} \geq \norm{g}_{\setA^s_{\rho,q,h_2}(\Omega)}.
  \end{align*}
  Therefore, the result follows from monotone convergence.
\end{proof}

The next result states a sort of lower semicontinuity of the discrete smoothness norms. In turn, this provides a compactness criterion for smoothness spaces. Once again, this result follows the spirit of those in \cite{MR3586796,MR4654593}.
\begin{theorem}[Lower semicontinuity]
\label{thm:norm-convergence-h0-discr}
  Assume the family $\{g_h \in \setA^s_{\rho,q,h}(\Omega)\}_{h > 0}$ is such that there is $g \in \setL^1(\Omega)$ for which, as $h\searrow 0$, we have $g_h \to g $ in $\setL^1(\Omega)$. Then
  \begin{align*}
    \abs{g}_{\setA^s_{\rho,q}(\Omega)} \leq \liminf_{h \searrow 0} \abs{g_h}_{\setA^s_{\rho,q,h}(\Omega)}.
  \end{align*}
\end{theorem}
\begin{proof}
  Observe that, for every $x\in \Omega$ and $t>0$,
  \[
    \osc g(x,t) \leq \osc g_h(x,t) + C_1(x,t) \norm{ g - g_h }_{\setL^1(\Omega)},
  \]
  where by $C_1(x,t)$ we denote a constant that may depend on $x$ and $t$.
  Consequently, the convergence $g_h \to g $ in $L^1(\Omega)$ implies that, for almost every $x \in \Omega$ and $t>0$,
  \[
    \osc g(x,t) \leq \liminf_{h \searrow 0} \osc g_h(x,t).
  \]
  Therefore, an application of Fatou's lemma implies
  \[
    |g|^\rho_{\setF^s_{\rho,q}(\Omega)} \leq \liminf_{h \searrow 0} \int_\Omega \bigg( \int_h^R \Big( \frac{\osc g_h(x,t)}{t^s}\Big)^q \frac{ \diff{t} }{t}\bigg)^{\frac{\rho}{q}} \diff{x},
  \]
  finishing the proof.
\end{proof}

Our next step is to study when finite element functions belong to the classes $\setA^s_{\rho, q}(\Omega)$.

\begin{lemma}[Semi-norm equivalence]
\label{lem:norm-equivalence-no-trace}
Let $\mathcal{T}_h$ be a triangulation on $\Omega$. If $s\rho < 1$, then $\big(\mathcal{L}^0_{k-1}(\mathcal{T}_h),\norm{\cdot}_{\setA^s_{\rho,q,h}}\big) \hookrightarrow \setA^s_{\rho,q}(\Omega)$, and we have
  \[
    \abs{g_h}_{\setA^s_{\rho,q}(\Omega)} \lesssim \abs{g_h}_{\setA^s_{\rho,q,h}(\Omega)},
    \qquad
    \forall g_h \in \mathcal{L}^0_{k-1}(\mathcal{T}_h),
  \]
  where the hidden constant depends on the regularity constant of $\Omega$.
\end{lemma}
\begin{proof}
  Let $g_h \in \mathcal{L}^0_{k-1}(\mathcal{T}_h)$. In the case $\setA = \setF$ and $q<\infty$, it is sufficient to show that
  \[
      \int_\Omega \bigg( \int_0^h \Big( \frac{\osc g_h(x,t)}{t^s}\Big)^q \frac{ \diff{t} }{t}\bigg)^{\frac{\rho}{q}} \diff{x} \lesssim
      \int_\Omega \bigg( \int_h^R \Big( \frac{\osc g_h(x,t)}{t^s}\Big)^q \frac{ \diff{t} }{t}\bigg)^{\frac{\rho}{q}} \diff{x}.
  \]
  
  We first observe that, for almost every $x \in \Omega$, there is a unique $T_x \in  \mathcal{T}_h$ such that $x \in \mathring{T}_x$. In the finite element language we can say that these points are not in the skeleton of the mesh. For such points, define $d_x \coloneqq \distance(x, \partial T_x)$ and observe then that, for $t \in (0,d_x]$, we have
  $
    B_t(x) \subset T_x
  $
  and so, since $g_{h|T_x} \in \mathcal{P}_{k-1}$,
  \[
    \osc g_h(x,t) = 0, \qquad \forall t \in (0,d_x].
  \]

    Next, let $x \in \Omega$ be a point not in the skeleton of the mesh and $t \in (0,h]$.  Using that $g_h$ is piecewise polynomial we may invoke local inverse estimates to assert for every $p \in \mathcal{P}_{k-1}$ that
  \begin{align*}
    \abs{(g_h-p)(z)} \lesssim \dashint_{B^\Omega_h(z)} \abs{g_h-p} \diff{y}
  \end{align*}
  for every $z \in B_t(x)$, and therefore with $B_h(z) \subset B_{3h}(y)$ for $y \in T_x$ that
  \[
    \osc g_h(x,t) \lesssim  \indicator_{t \geq d_x} \dashint_{T_x} \osc g_h(y,3h) \diff{y}.
  \]
With these observations, we estimate
  \begin{align*}
    \MoveEqLeft \int_\Omega \bigg( \int_0^h \Big( \frac{\osc g_h(x,t)}{t^s}\Big)^q \frac{ \diff{t} }{t}\bigg)^{\frac{\rho}{q}} \diff{x} \\
    &\lesssim
    \int_\Omega \dashint_{T_x} \osc g_h(y,3h)^\rho \bigg( \int_{0}^h \Big( \frac{\indicator_{t\ge d_x}}{t^s}\Big)^q \frac{ \diff{t} }{t}\bigg)^{\frac{\rho}{q}} \diff{y}\diff{x} \\
    &\lesssim
    \int_\Omega \dashint_{T_x} \osc g_h(y,3h)^\rho d_x^{-s\rho} \diff{y}\diff{x}
    = \int_\Omega \osc g_h(y,3h)^\rho \dashint_{T_y} d_x^{-s\rho} \diff{x}\diff{y} \\
    &\lesssim h^{-s\rho} \int_\Omega \osc g_h(y,3h)^\rho \diff{y}.
  \end{align*}
  Note that $s\rho < 1$ was necessary for the integral over $d_x^{-s\rho}$ to be finite. Finally, we use that, for $3h \leq t \leq 4h$, we have
  \[
    \osc g_h(x,3h) \lesssim \osc g_h(x,t)
  \]
  to get
  \begin{align*}
    h^{-s\rho} \int_\Omega \osc g_h(x,3h)^\rho \diff{x}
    &\lesssim h^{-s\rho} \int_\Omega \bigg( \int_{3h}^{4h} \osc g_h(x,t)^q \frac{ \diff{t} }{t} \bigg)^\frac{\rho}{q} \diff{x} \\
    &\lesssim
    \int_\Omega \bigg( \int_{3h}^{4h} \Big( \frac{\osc g_h(x,t)}{t^s}\Big)^q \frac{ \diff{t} }{t}\bigg)^{\frac{\rho}{q}} \diff{x}.
  \end{align*}
Per Definition~\ref{def:discrete-domain}, $4h \leq R$. This finishes the case of Triebel--Lizorkin spaces.
  
To conclude we quickly comment on the Besov space case, that is $\setA = \setB$. In this case, for $t \in (0,h]$, similar to above we calculate
\begin{align*}
  \norm{ \osc g_h(\cdot, t) }_{\setL^\rho(\Omega)}^\rho
  &\lesssim
  \int_\Omega \indicator_{d_x \leq t} \bigg( \dashint_{T_x} \osc g_h(y,3h) \diff{y} \bigg)^\rho \diff{x} \\
  &\lesssim
  \int_\Omega \osc g_h(y,3h)^\rho \dashint_{T_y} \indicator_{d_x \leq t} \diff{x} \diff{y} \\
  &\lesssim
  \frac{t}{h} \norm{\osc g_h(\cdot,3h)}_{\setL^\rho(\Omega)}^\rho,
\end{align*}
  so that
  \begin{align*}
    \int_0^h t^{-sq} \norm{ \osc g_h(\cdot, t) }_{\setL^\rho(\Omega)}^q \frac{\diff{t}} t
    &\lesssim \norm{\osc g_h(\cdot,3h)}_{\setL^\rho(\Omega)}^q \int_0^h \Big(\frac{t}{h}\Big)^{q/\rho} \frac{\diff{t}}{t^{1+sq}} \\
    &\lesssim h^{-sq} \norm{ \osc g_h(\cdot, 3h) }_{\setL^\rho(\Omega)}^q ,
  \end{align*}
  where, once again, the condition $s\rho <1$ guarantees convergence of the integral. At this point we may continue with similar arguments as given before.

  This concludes the proof.
\end{proof}

In view of Lemma~\ref{lem:norm-equivalence-no-trace}, for $s\rho < 1$, it is possible to use the standard uniform embedding and compactness results of  $\setA^{s}_{\rho,q}(\Omega)$, see~Remark~\ref{rem:sobolev-embeddings}, to conclude embedding and compactness properties of the scale $\setA^{s}_{\rho,q,h}(\Omega)$. For brevity we do not state these results since, below, we will remove the restriction $s\rho <1$.

\subsection{Smooth companion}
\label{sub:SmoothCompation}

We observe that, if $s \rho \geq1$, some of the results of the previous section are not applicable. For this reason here we introduce the notion of smooth companion, that is, a regularized version of $g_h \in \mathcal{L}_{k-1}^0(\mathcal{T}_h)$ that behaves similarly but which uniformly belongs to $\setA^s_{\rho,q}(\Omega)$.

\begin{definition}[Smooth companion]
\label{def:SmoothCompanion}
We call a family of linear maps $\{\mathcal{S}_h\}_{h>0}$,
  \[
    \mathcal{S}_h : \setL^1(\Omega) \to \setC^{k}(\overline{\Omega}),
  \]
  a smooth companion operator, if it fulfills the following properties:
  \begin{enumerate}[label=\alph{*}), leftmargin=*]
    \item \label{itm:SmoothCompanionPolynomial}
      \textbf{Polynomial preserving:} For all $p \in \mathcal{P}_{k-1}$ there holds $\mathcal{S}_h p = p$.
      
    \item \label{itm:SmoothCompanionL1stab}
      \textbf{Local $L^1$-stability:} There is $\varrho \in \setN$ such that, for every $T \in \mathcal{T}_h$,
      \[
        \norm{\mathcal{S}_h g }_{\setL^1(T)} \lesssim \norm{ g }_{\setL^1(\omega^{\varrho}(T))}, \qquad \forall g \in \setL^1(\Omega).
      \]

    \item \label{itm:SmoothCompanionOscSmallt} \textbf{Inverse estimate:} For every $T \in \mathcal{T}_h$,
    \[
      \| \nabla^k \mathcal{S}_h g \|_{\setL^\infty(T)} \lesssim h^{-k} \dashint_{\omega^{\varrho}(T)} \osc g (x, h) \diff x, \qquad \forall g \in \setL^1(\Omega).
    \]
  \end{enumerate}
\end{definition}

  The notion of a smooth companion is closely related to conforming reconstruction and enrichment operators. Reconstructions of discontinuous Galerkin functions into continuous Lagrange spaces are studied in~\cite{MR2034620} and in the appendix of~\cite{DieningKroenerRuzickaToulopoulos2014}. For fourth-order problems, enrichment and recovery operators into $\setW^{2,2}$-conforming spaces are developed in~\cite[Section~4]{MR2777246} and~\cite[Section~3]{MR2755946}, while conforming companions for the Morley element are constructed in~\cite[Section~2.3]{Gallistl2015} and~\cite[Section~5.1]{CarstensenPuttkammer2020}. Beyond restoring conformity, companions can also provide estimates in stronger standard Sobolev norms, as shown for the Crouzeix--Raviart and Morley methods in~\cite[Section~4.2]{CarstensenNataraj2021}.

A related approach is smoothing on a spatially varying length scale. Such operators, with local approximation estimates and bounds in higher-order Sobolev norms, are constructed in~\cite{KarkulikMelenk2015}. Spline quasi-interpolants provide another related construction; see~\cite[Sections~5.7, 10.3, and 17.4]{MR2355272}. Our definition specifies the locality, polynomial preservation, and regularity properties needed below to estimate smooth reconstructions in Besov and Triebel--Lizorkin spaces.

The parameter $\varrho$
accommodates the enlarged patches involved in such local constructions. In particular, $(\mathcal{S}_h g)_{|T}$ may depend on the values of $g$ on $\omega^{\varrho}(T)$ for some $\varrho \geq 1$. To simplify the presentation, we assume below that $\varrho = 1$; the arguments extend to the general case with the corresponding patch enlargements.

This definition directly implies the following properties.

\begin{proposition}[Properties of $\mathcal{S}_h$]
\label{prop:PropertiesSh}
  The smooth companion, introduced in Definition~\ref{def:SmoothCompanion}, satisfies the following properties.
  \begin{enumerate}[label=\alph{*}), leftmargin=*]
    \item \textbf{Local $\setL^1$-stability:} For a discrete subdomain $D \subseteq \Omega$, $t \geq h$, $x \in D$, and $g \in \setL^1(\Omega)$ we have
    \begin{align*}
      \norm{\mathcal{S}_h g}_{\setL^1(B_t(x) \cap D)} \lesssim \norm{g}_{\setL^1(B_{t+h}(x) \cap D^+)},
    \end{align*}
    where
    \[
      D^+ \coloneqq \big\{ x \in \Omega \ : \ \distance( x, D ) \leq 2h \big\}.
    \]
    Here, the implicit constant does not depend on $h$.
    
    \item \textbf{Large-scale stability:} For $t \geq h$, and $g \in \setL^1(\Omega)$ we have, for every $x \in D$, the oscillation estimate
    \begin{align*}
       \osc^D \mathcal{S}_h g(x,t) \lesssim \osc^{D^+} g(x,t+ h).
    \end{align*}
    Here, the implicit constant does not depend on $h$.
    
    \item \textbf{Small-scale smoothing:} If $D \subseteq \Omega$ is a discrete subdomain, $x \in D$, $t \leq h$, and $g \in \setL^1(\Omega)$ there holds the restricted oscillation estimate
      \begin{align*}
        \osc^D \mathcal{S}_h g(x,t) \lesssim \Big(\frac{t}{h}\Big)^k \osc^{{D^+}} g(x,3h),
      \end{align*}
      where, the implicit constant does not depend on $h$ either, and $D^+$ is defined as before.
  \end{enumerate}
\end{proposition}
\begin{proof}
  We prove each result separately.
  \begin{enumerate}[label=\alph{*}), leftmargin=*]
    \item One merely needs to add the local $\setL^1$-stability over $T \in \mathcal{T}_h$ that intersect $B_t(x)$.
    
    \item Let $p \in \mathcal{P}_{k-1}$. By the previous property and the $\setL^1$-stability, there holds
    \begin{align*}
      \norm{\mathcal{S}_h g - p}_{\setL^1(B_t(x) \cap D)}
      &=
      \norm{\mathcal{S}_h (g - p)}_{\setL^1(B_t(x) \cap D)}
      \lesssim
      \norm{g - p}_{\setL^1(B_{t+ h}(x) \cap D^+) },
    \end{align*}
    where we used that $t \eqsim t + h$ and \eqref{eq:BiteOffEstimate}.
    Taking the infimum on both sides gives the desired oscillation estimate.
    
    \item Since $\mathcal{S}_h g \in \setC^k(\overline\Omega)$, we have that
    \[
      \osc^D \mathcal{S}_h g(x,t) \lesssim t^k \norm{\nabla^k \mathcal{S}_h g }_{\setL^\infty(B_t \cap D)}.
    \]
    The inverse estimate then implies the result. \qedhere
  \end{enumerate}
\end{proof}

For completeness, we provide an explicit example of a smooth companion.

\begin{example}
  Let $\{\phi_T\}_{T \in \mathcal{T}_h}$ be a partition of unity on~$\overline{\Omega}$ aligned with the triangulation $\mathcal{T}_h$, such that $\phi_T \in \setC^\infty(\overline{\Omega})$, $\support_{\overline{\Omega}} \phi_T \subset \omega(T)$, and, for $i \in \mathbb{N}$,
  \[
    \norm{\nabla^i \phi_T}_{\setL^\infty(\Omega)} \lesssim h^{-i}.
  \]
  Furthermore, for every element $T \in \mathcal{T}_h$, let $\Pi_T \colon \setL^1(T) \to \mathcal{P}_{k-1}$ be the orthogonal projection with respect to the $\setL^2(T)$ inner product. Then, for $g \in \setL^1(\Omega)$, define
  \begin{align*}
    \mathcal{S}_h g(x) \coloneqq \sum_{T \in \mathcal{T}_h} \phi_T(x) (\Pi_T g)(x).
  \end{align*}
  The projection property $\mathcal{S}_h p = p$ for $p \in \mathcal{P}_{k-1}$ follows immediately. Similarly, one has local $\setL^1$-stability of the form
  \[
    \norm{\mathcal{S}_h g}_{\setL^1(T)} \lesssim \norm{g}_{\setL^1(\omega(T))}.
  \]
  It remains to verify property~\ref{itm:SmoothCompanionOscSmallt}. Fix $T_0 \in \mathcal{T}_h$ and $x \in T_0$. For an arbitrary $p \in \mathcal{P}_{k-1}$ we may compute
  \begin{align*}
    |\nabla^k \mathcal{S}_h g(x)| &= \abs{ \nabla^k ( \mathcal{S}_h g(x) - p(x) )} 
    = \abs{ \nabla^k  \mathcal{S}_h (g - p )(x)} \\
    &\lesssim \sum_{T \in \mathcal{T}_h} \sum_{l=0}^k \abs{\nabla^l \varphi_T(x)} \abs{\nabla^{k-l}\Pi_T (g-p)(x)} \\
    &\lesssim \sum_{l=0}^k \sum_{T \in \omega(T_0)} \frac1{h_T^l} \dashint_{T} \abs{\nabla^{k-l}\Pi_T(g-p)} \diff y \\
    &\lesssim h^{-k}\sum_{T \in \omega(T_0) } \dashint_T\abs{\Pi_T(g-p)} \diff y
    \lesssim h^{-k}\dashint_{\omega(T_0) } \abs{g-p} \diff y.
  \end{align*}
  Since $p$ is arbitrary, a chaining ball argument for polynomial approximations gives
    \begin{align*}
      \inf_{p \in \mathcal{P}_{k-1}} \dashint_{\omega(T_0) } \abs{g-p} \diff y
      &\lesssim \dashint_{\omega(T_0) } \osc g(y,h) \diff y,
    \end{align*}
  and the property follows; compare to~\cite[Section 7]{MR559195} and \cite{MR2643399} for similar covering and chaining arguments.
\end{example}

For functions in $\mathcal{L}_{k-1}^0(\mathcal{T}_h)$ the smooth companion gives a way to remove the restriction $s \rho <1$ when comparing our discrete scales of spaces with the standard ones.

\begin{proposition}[Stability]
\label{proplem:stab-smoothing}
  Let $D \subseteq \Omega$ be a discrete subdomain and let
  \[
    D^+ \coloneqq \big\{x \in \Omega \ : \ \distance(x,D) < 2 h \big\}.
  \]
  Then we have
  \[
    \abs{\mathcal{S}_h g}_{\setA^s_{\rho,q}(D)} \lesssim \abs{g}_{\setA^s_{\rho,q,h}(D^+)}, \qquad \forall g \in \setL^1(\Omega).
  \]
\end{proposition}
\begin{proof}
  Let $R = \diameter D$. It suffices to estimate the seminorm
  \begin{align*}
    \MoveEqLeft \abs{\mathcal{S}_h g}^\rho_{\setF^s_{\rho,q}(D)}
    \lesssim \int_D \bigg( \int_0^R \Big( \frac{\osc^D \mathcal{S}_h g(x,t)}{t^s}\Big)^q \frac{ \diff{t} }{t}\bigg)^{\frac{\rho}{q}} \diff{x} \\
    &\lesssim \int_D \biggabs{ \int_0^{h} \Big( \frac{\osc^D \mathcal{S}_h g(x,t)}{t^s}\Big)^q \frac{ \diff{t} }{t}}^{\frac{\rho}{q}} \diff{x} 
    + \int_D \biggabs{ \int_{h}^R \Big( \frac{\osc^D \mathcal{S}_h g(x,t)}{t^s}\Big)^q \frac{ \diff{t} }{t}}^{\frac{\rho}{q}} \diff{x} \\
    &\eqqcolon \mathrm{I} + \mathrm{II},
  \end{align*}
  and we consider each term separately.
  
  By Proposition~\ref{prop:PropertiesSh} we have, if $t \leq h$, that
  \begin{align*}
    \osc^D \mathcal{S}_h g(x,t)
    &\lesssim \frac{t^k}{h^k} \osc^{D^+} g (x,3h).
  \end{align*}
  Therefore, we have
  \begin{align*}
    \mathrm{I}
    &\lesssim h^{-k\rho}\int_D \bigg( \int_0^h t^{(k-s)q-1} \diff{t} \bigg)^{\frac{\rho}{q}} \osc^{D^+} g(x,3h)^\rho \diff{x} \lesssim h^{-s\rho} \int_D  \osc^{D^+} g(x,3h)^\rho \diff{x}.
   \end{align*}
  Next we use that, for $t \in (3h,4h)$, we have
  \[
    \osc^{D^+} g(x,3h) \lesssim \osc^{D^+} g (x,t)
  \]
  Consequently, we can finish estimating the first term as follows:
  \begin{align*}
    \mathrm{I} &\lesssim h^{-s\rho} \int_D \osc^{D^+} g(x,3h)^\rho \diff{x}
    \lesssim h^{-s\rho} \int_D \bigg( \int_{3h}^{4h} \osc^{D^+} g(x,t)^q \frac{ \diff{t} }{t} \bigg)^\frac{\rho}{q} \diff{x}
    \\
   &\lesssim \int_D \bigg( \int_{3h}^{4h} \Big( \frac{\osc^{D^+} g(x,t)}{t^s}\Big)^q \frac{ \diff{t} }{t}\bigg)^{\frac{\rho}{q}} \diff{x}
    \leq \abs{g}^\rho_{\setF^s_{\rho,q,h}(D^+)}.
  \end{align*}

  It remains to bound the term $\mathrm{II}$. For every $t > h$ by Proposition~\ref{prop:PropertiesSh}, there holds
  \begin{align*}
    \osc^D \mathcal{S}_h g(x,t)
    &\lesssim
    \osc^{D^+} g(x,t+h).
    \end{align*}
  Using that
  \[
    \frac{t+h}t = 1 + \frac{h}t \leq 2,
  \]
  it is trivial to compute
  \begin{align*}
    \mathrm{II} &= \int_D \bigg( \int_{h}^R \Big( \frac{\osc^D \mathcal{S}_h g(x,t)}{t^s}\Big)^q \frac{ \diff{t} }{t}\bigg)^{\frac{\rho}{q}} \diff{x} \\
    &\lesssim \int_D \bigg( \int_h^R \frac{ \osc^{D^+} g(x,t+h)^q }{ (t+h)^{1+sq} } \Big( \frac{t+h}t \Big)^{1+sq} \diff{t}\bigg)^\frac\rho{q}  \diff{x} \\
    &\lesssim \int_D \bigg( \int_{2h}^{R+h} \Big( \frac{\osc^{D^+}  g(x,t)}{t^s}\Big)^q \frac{ \diff{t} }{t}\bigg)^{\frac{\rho}{q}} \diff{x} \\
    &\leq \abs{g}^\rho_{\setF^s_{\rho,q,h}(D^+)}.
  \end{align*}

  Collecting the obtained bounds for $\mathrm{I}$ and $\mathrm{II}$ finishes the proof.
\end{proof}

We next present an auxiliary result.

\begin{lemma}[Integral average]
  Let $D \subset \RRn$ be regular and $g \in \setL^1(D)$. For $t>0$ define
  \[
    g_t(x) \coloneqq \dashint_{B_t(x)\cap D} g(y) \diff y, \qquad x \in D.
  \]
  Then, for every $\rho \in [1,\infty]$, we have
  \begin{equation}
  \label{eq:AveragesEqNorm}
    \| g_t \|_{\setL^\rho(D)} \lesssim \| g \|_{\setL^\rho(D)},
  \end{equation}
  where the implicit constant is independent of $t$. 
\end{lemma}
\begin{proof}
  Clearly,
  \[
    \| g_t \|_{\setL^\infty(D)} \lesssim \| g \|_{\setL^\infty(D)}.
  \]
  Next we consider the case $\rho = 1$:
  \begin{align*}
    \| g_t \|_{\setL^1(D)} &\leq \int_D \dashint_{B_t(x) \cap D} |g(y)| \diff y \diff x = \int_D \int_D \frac1{|B_t(x) \cap D|} \indicator_{B_t(x) \cap D}(y) |g(y)| \diff y \diff{x}
    \\
    &= \int_D |g(y)| \int_D \frac1{|B_t(x) \cap D|} \indicator_{B_t(x) \cap D}(y) \diff{x} \diff{y}.
  \end{align*}
  Now, since $D$ is assumed to be regular, see \eqref{eq:BiteOffEstimate},
  \[
    \sup_{y \in D}\int_D \frac1{|B_t(x) \cap D|} \indicator_{B_t(x) \cap D}(y) \diff{x} \lesssim 1,
  \]
  so that
  \[
    \| g_t \|_{\setL^1(D)} \lesssim \| g \|_{\setL^1(D)}.
  \]
  To conclude \eqref{eq:AveragesEqNorm} one only needs to invoke the Riesz--Thorin interpolation theorem.
  
\end{proof}

Having proved this simple result we proceed to discuss approximation properties of the smooth companion. This will be useful to prove embedding results.

\begin{lemma}[Local error estimate]
  \label{lem:smooth_comp_approx}
  Let $D \subseteq \Omega$ be a discrete subdomain. We have
  \begin{align}
    \label{eq:smooth_comp_approx_1}
    \abs{g-\mathcal{S}_h g}_{\setA^s_{\rho,q,h}(D)} \lesssim h^{-s} \norm{\osc g(\cdot,2h)}_{\setL^\rho(D)}, \qquad \forall g \in \setL^\rho(\Omega).
  \end{align}
  In addition, if $g_h \in \mathcal{L}^0_{k-1}(\mathcal{T}_h)$ we have
  \begin{align}
    \label{eq:smooth_comp_approx_2}
    \norm{g_h -\mathcal{S}_h g_h}_{\setL^\rho(D)} \lesssim \norm{\osc g_h(\cdot,2h)}_{\setL^\rho(D)}.
  \end{align}
\end{lemma}
\begin{proof}
  Let $T \in \mathcal{T}_h$. Polynomial invariance and local $\setL^1$-stability of $\mathcal{S}_h$ give
  \[
    \dashint_T |g - \mathcal{S}_hg| \diff y \lesssim \dashint_{\omega(T)} \osc g(x,2h) \diff x.
  \]
  This estimate implies, for $t \geq h$,
  \[
    \dashint_{B_t(x) \cap D} |g - \mathcal{S}_hg| \diff y  \lesssim \dashint_{B_{t+h}(x) \cap D} \osc g(x,2h) \diff x.
  \]
  Next we observe that, for $x \in D$,
  \[
    \osc^D (g-\mathcal{S}_h g)(x,t) \leq \dashint_{B_t(x) \cap D} \abs{ g - \mathcal{S}_h g} \diff y .
  \]
  Therefore, estimate \eqref{eq:AveragesEqNorm} gives
  \[
    \norm{ \osc^D (g-\mathcal{S}_h g)(\cdot,t) }_{\setL^\rho(D)} \lesssim \norm{ \osc g(\cdot,2h) }_{\setL^\rho(D)}.
  \]
  We now use this estimate in the definition of the discrete Besov seminorm and get
  \begin{align*}
    \abs{g-\mathcal{S}_h g}_{\setB^s_{\rho,q,h}(D)}
    &=
    \bigg( \int_h^R t^{-sq} \norm{\osc^D (g-\mathcal{S}_h g)(\cdot, t)}_{\setL^\rho(D)}^q \frac{\diff{t}}t \bigg)^{1/q} \\
    &\lesssim
    \bigg( \int_h^R t^{-sq} \norm{\osc g(\cdot, 2h)}_{\setL^\rho(D)}^q \frac{\diff{t}}t \bigg)^{1/q} \\
    &\lesssim
    h^{-s}\norm{\osc g(\cdot, 2h)}_{\setL^\rho(D)}.
  \end{align*}
  We now consider the Triebel--Lizorkin seminorm. Observe that, for every $x \in D$ and all $t \leq r$ with $r \eqsim t$, we have
  \[
    \osc v(x,t) \lesssim \osc v(x,r).
  \]
  Let $I \in \setN$ be the smallest integer such that $2^Ih \eqsim R$. With the help of a dyadic partition and the embedding $\ell^1 \hookrightarrow \ell^q$ we are able to calculate
\begin{align*}
  \abs{g-\mathcal{S}_h g}_{\setF^s_{\rho,q,h}(D)}
  &=
  \Biggnorm{\bigg(\int_h^R \Big( \frac{\osc^D (g-\mathcal{S}_h g)(\cdot,t)}{t^s} \Big)^q \frac{\diff{t} }t\bigg)^{\frac1{q}}}_{\setL^\rho(D)} \\
  &\lesssim h^{-s} \Biggnorm{\Bigg(\sum_{i=0}^I\Big( \frac{\osc^D (g-\mathcal{S}_h g)(\cdot,2^i h)}{2^{is}} \Big)^q\Bigg)^{\frac1{q}}}_{\setL^\rho(D)} \\
  &\leq h^{-s} \Bignorm{\sum_{i =0}^I \frac{\osc^D (g-\mathcal{S}_h g)(\cdot,2^i h)}{2^{is}}}_{\setL^\rho(D)} \\
  &\leq h^{-s} \norm{\osc g(\cdot,2h)}_{\setL^\rho(D)} \sum_{i= 0}^\infty \frac1{2^{is}}.
\end{align*}
  Since the series appearing in the last term is convergent, this shows \eqref{eq:smooth_comp_approx_1}.
  
  We now observe that, for $T \in \mathcal{T}_h$ and $v \in \setC^k(T)$, we may choose $p\in \mathcal{P}_{k-1}$ to show that
  \begin{align}
  \label{eq:rev-estimate-Ck}
    \begin{aligned}
    \norm{v}_{\setL^\infty(T)} &\leq \norm{p}_{\setL^\infty(T)} + \norm{v-p}_{\setL^\infty(T)} \lesssim h^{-n}\norm{p}_{\setL^1(T)} + h^k \norm{\nabla^kv}_{\setL^\infty(T)} \\
    &\lesssim h^{-n} \norm{v}_{\setL^1(T)}+ h^k \norm{\nabla^kv}_{\setL^\infty(T)}.
    \end{aligned}
  \end{align}
  Set $v = \mathcal{S}_h (g - p)$. The local $\setL^1$-stability \ref{itm:SmoothCompanionL1stab} and inverse estimate \ref{itm:SmoothCompanionOscSmallt} give, for every $x \in T$,
  \begin{equation}
  \label{eq:UsefulLater}
    \abs{\mathcal{S}_h (g - p)(x)} \lesssim \dashint_{ \omega(T) } \osc g(y, 2h) \diff y.
  \end{equation}
  Let now $g_h \in \mathcal{L}_{k-1}^0(\mathcal{T}_h)$. By inverse estimates for finite element functions and \eqref{eq:UsefulLater} we derive
  \[
    |g_h(x) - \mathcal{S}_h g_h(x)| \leq \abs{g_h(x) - p(x)} + \abs{\mathcal{S}_h(g-p)(x)}\lesssim \dashint_{ \omega(T) } \osc g_h(y, 2h) \diff y.
  \]
  Taking $\setL^\rho$-norm and using \eqref{eq:AveragesEqNorm} proves \eqref{eq:smooth_comp_approx_2}.
\end{proof}

\begin{corollary}[$\setL^\rho$-stability]
\label{col:SmoothIsLrhostable}
  For every $\rho \in [1,\infty]$ the smooth companion is $\setL^\rho$-stable, that is
  \[
    \norm{\mathcal{S}_h g }_{\setL^\rho(\Omega)} \lesssim \norm{ g }_{\setL^\rho(\Omega)} .
  \]
\end{corollary}
\begin{proof}
  We may use~\eqref{eq:rev-estimate-Ck}, the local $\setL^1$-stability \ref{itm:SmoothCompanionL1stab} and inverse estimate \ref{itm:SmoothCompanionOscSmallt} to obtain
  \begin{align*}
    \abs{\mathcal{S}_h g(x)}
    &\lesssim \dashint_{\omega(T)} \abs{g(y)} \diff{y} + \dashint_{\omega(T)} \osc g(y,2h) \diff{y}
    \lesssim \dashint_{B^\Omega_{3h}(x)} \abs{g(y)} \diff{y}
  \end{align*}
  for $x \in T$. We integrate to get
  \[
    \norm{\mathcal{S}_h g }_{\setL^\rho(\Omega)}
    \lesssim \biggnorm{ \dashint_{B^\Omega_{3h}(\cdot)} \abs{g(y)} \diff y }_{\setL^\rho(\Omega)},
  \]
  and invoke \eqref{eq:AveragesEqNorm} to conclude. The case $\rho = \infty$ is simpler and thus omitted.
\end{proof}

We can now extend the embedding results in the discrete spaces and drop the assumption $s_i\rho_i < 1$.

\begin{lemma}[Uniform embedding]
  \label{lem:UniformEmbeddingImprov}
  Let, for $i =1,2$, $\rho_i$, $q_i$, and $s_i$ be as in Definition~\ref{def:besovTriebel}.
  \begin{enumerate}[label={\alph{*})},leftmargin=*]
    \item Let the indices be such that the embedding
    \[
      \setA^{s_2}_{\rho_2,q_2}(\Omega) \hookrightarrow \setA^{s_1}_{\rho_1,q_1}(\Omega)
    \]
    holds. Then we have
    \begin{align*}
      \norm{g_h}_{\setA^{s_1}_{\rho_1,q_1,h}(\Omega)} \lesssim \norm{g_h}_{\setA^{s_2}_{\rho_2,q_2,h}(\Omega)}, \qquad \forall g_h \in \mathcal{L}_{k-1}^0(\mathcal{T}_h),
    \end{align*}
    with a constant that is uniform in $h$.

    \item Assume that the indices allow for the compact embedding
    \[
      \setA^{s_2}_{\rho_2,q_2}(\Omega) \hookrightarrow\hookrightarrow \setA^{s_1}_{\rho_1,q_1}(\Omega).
    \]
    Let the family of functions $\{g_h \in \mathcal{L}^0_{k-1}(\mathcal{T}_h)\}_{h > 0}$ be such that
    \begin{align*}
      \sup_{h > 0} \norm{g_h}_{\setA^{s_2}_{\rho_2,q_2,h}(\Omega)} \leq C < \infty.
    \end{align*}
    Then, $\{g_h\}_{h > 0}$ is precompact in $\setL^{\rho_1}(\Omega)$. In addition, any limit point $g \in \setL^{\rho_1}(\Omega)$ belongs to $\setA^{s_2}_{\rho_2,q_2}(\Omega)$ and there holds
    \begin{align*}
      \liminf_{h \searrow 0} \norm{g-g_h}_{\setA^{s_1}_{\rho_1,q_1,h}(\Omega)} = 0.
    \end{align*}
  \end{enumerate}
\end{lemma}
\begin{proof}
  We prove each statement separately.
  \begin{enumerate}[label={\alph{*})},leftmargin=*]
    \item   We introduce the smooth companion and apply the triangle inequality to estimate
    \begin{align*}
      \norm{g_h}_{\setA^{s_1}_{\rho_1,q_1,h}(\Omega)}
      &\leq
      \norm{\mathcal{S}_h g_h}_{\setA^{s_1}_{\rho_1,q_1,h}(\Omega)}
      + \norm{g_h - \mathcal{S}_h g_h}_{\setA^{s_1}_{\rho_1,q_1,h}(\Omega)}.
    \end{align*}
    For the first term we employ the assumed embedding and the stability result of Proposition \ref{proplem:stab-smoothing} and Corollary~\ref{col:SmoothIsLrhostable} to arrive at
    \begin{align*}
      \norm{\mathcal{S}_h g_h}_{\setA^{s_1}_{\rho_1,q_1,h}(\Omega)}
      &\leq
      \norm{\mathcal{S}_h g_h}_{\setA^{s_1}_{\rho_1,q_1}(\Omega)}
      \lesssim
      \norm{\mathcal{S}_h g_h}_{\setA^{s_2}_{\rho_2,q_2}(\Omega)}
      \lesssim
      \norm{g_h}_{\setA^{s_2}_{\rho_2,q_2,h}(\Omega)}.
    \end{align*}
    Now, to estimate the second term we observe that
    \[
      \osc g_h(x,t) \lesssim \osc g_h(y,\lambda t)
    \]
    for some $\lambda > 1$ and $\abs{x-y} \leq (\lambda-1)t$. Averaging over $y$ leads to
    \begin{align*}
      \osc g_h(x,t) \lesssim (t^{-n}\indicator_{B_{(\lambda-1)t}(0)}) \star \osc g_h(x,\lambda t).
    \end{align*}
    Next we observe that, without loss of generality $\rho_2 \leq \rho_1$. Therefore, Young's convolution estimate implies
    \begin{align*}
      \norm{\osc g_h(\cdot,2h)}_{\setL^{\rho_1}(\Omega)} \lesssim h^{\frac{n}{\rho_1}-\frac{n}{\rho_2}} \norm{\osc g_h(\cdot,3h)}_{\setL^{\rho_2}(\Omega)}.
    \end{align*}
    This last estimate, combined with Lemma \ref{lem:smooth_comp_approx}, allows us to conclude
    \begin{align*}
        \norm{g_h - \mathcal{S}_h g_h}_{\setA^{s_1}_{\rho_1,q_1,h}(\Omega)}
        &\lesssim h^{-s_1} \norm{\osc g_h (\cdot,2h)}_{\setL^{\rho_1}(\Omega)}
        \\
        &\lesssim h^{-s_1+\frac{n}{\rho_1}-\frac{n}{\rho_2}} \norm{\osc g_h (\cdot,3h)}_{\setL^{\rho_2}(\Omega)}.
    \end{align*}
    Since, for all $3h \leq t \leq 4h$, we have $\osc g_h(\cdot,3h) \lesssim \osc g_h(\cdot,t)$ we may then finally estimate
    \begin{equation}
    \label{eq:smooth_comp_approx_sobolev}
    \begin{aligned}
      \norm{g_h - \mathcal{S}_h g_h}_{\setF^{s_1}_{\rho_1,q_1,h}(\Omega)}
      &\lesssim h^{-s_1+\frac{n}{\rho_1}-\frac{n}{\rho_2}} \norm{\osc g_h (\cdot,3h)}_{\setL^{\rho_2}(\Omega)} \\
      &\lesssim h^{s_2-s_1+\frac{n}{\rho_1}-\frac{n}{\rho_2}} \Biggnorm{\bigg[ \int_{3h}^{4h} \biggabs{ \frac{\osc g_h(\cdot,t)}{t^{s_2}}}^{q_2} \frac{\diff t}{t}\bigg]^\frac1{q_2}}_{\setL^{\rho_2}(\Omega)} \\
      &\leq h^{s_2-s_1+\frac{n}{\rho_1}-\frac{n}{\rho_2}} \abs{g_h}_{\setF^{s_2}_{\rho_2,q_2,h}(\Omega)}.
    \end{aligned}
  \end{equation}
  By assumption, the power of $h$ is nonnegative. For the Besov case, the last step follows similarly.

\item In this case we note that the compact embedding assumption and the stability result of Proposition~\ref{proplem:stab-smoothing} imply that the family $\{\mathcal{S}_h g_h\}_{h > 0}$ is precompact in $\setA^{s_1}_{\rho_1,q_1}(\Omega)$. Let $g \in \setA^{s_1}_{\rho_1,q_1}(\Omega)$ be any limit point of this family. In other words, there exists a sequence $\{h_j\}_{j \in \setN}$ with $h_j \to 0$ and $\mathcal{S}_{h_j} g_{h_j} \to g$ in $\setA^{s_1}_{\rho_1,q_1}(\Omega)$ as $j \to \infty$. We have
  \begin{align*}
    \lim_{j \to \infty} \norm{g-g_{h_j}}_{\setA^{s_1}_{\rho_1,q_1,{h_j}}(\Omega)}
    &\leq
    \lim_{j\to\infty}\norm{g-\mathcal{S}_{h_j} g_{h_j}}_{\setA^{s_1}_{\rho_1,q_1,{h_j}}(\Omega)}
    \\
    &+ \lim_{j\to\infty}\norm{g_{h_j}-\mathcal{S}_{h_j} g_{h_j}}_{\setA^{s_1}_{\rho_1,q_1,{h_j}}(\Omega)}.
  \end{align*}
  Since $g$ is the limit point, the first term is zero. For the second term, we reuse estimate \eqref{eq:smooth_comp_approx_sobolev} to get
 \begin{align*}
   \lim_{j\to\infty}\norm{g_{h_j}-\mathcal{S}_{h_j} g_{h_j}}_{\setA^{s_1}_{\rho_1,q_1,{h_j}}(\Omega)}
  & \lesssim
  \lim_{j\to\infty}h_j^{s_2-s_1+\frac{n}{\rho_1}-\frac{n}{\rho_2}} \abs{g_{h_j}}_{\setA^{s_2}_{\rho_2,q_2,{h_j}}(\Omega)} \\
  &\leq C \lim_{j\to\infty}h_j^{s_2-s_1+\frac{n}{\rho_1}-\frac{n}{\rho_2}}  = 0.
 \end{align*}
 Here, we used that for a compact embedding, there must hold $s_2-\frac{n}{\rho_2} > s_1-\frac{n}{\rho_1}$. Finally, $g \in \setA^{s_2}_{\rho_2,q_2}(\Omega)$ follows from Theorem \ref{thm:norm-convergence-h0-discr}.
  \end{enumerate}

  The proof is complete.
\end{proof}

\begin{remark}[The case $s_1 =0$]
  While Lemma \ref{lem:UniformEmbeddingImprov} is only stated for embeddings into spaces with $s_1 > 0$, the exact same techniques also give embeddings into $\setL^\rho(\Omega)$.
\end{remark}

\section{Regularity transfer}
\label{sec:StabRitz}

In this section we use the oscillation reduction estimate of Theorem~\ref{thm:osc-estimate} to prove the stability of the Ritz projection, defined in \eqref{eq:RitzProjection}, in the discrete smoothness norms we defined in Section~\ref{sec:finite-frequency}. Since the bounds on the Ritz projection are in terms of smoothness norms of the original function, this can be understood as a regularity transfer.

To state it we recall that, in Definition \ref{def:disc-besovTriebel}, the oscillations are measured down to scale $h$. To apply Theorem~\ref{thm:osc-estimate}, however, we need oscillations down to scale $\theta\Lambda h$. Since $\theta$ and $\Lambda$ are fixed, those norms are equivalent on balls of certain minimum size. In the following, we will therefore still write the discrete norms down to scale $h$ while starting the oscillations at $\theta\Lambda h$, as justified by their equivalence.
In the following result, the restriction on the indices guarantees that $\setA^s_{\rho, q} \hookrightarrow \setL^2$; see Remark~\ref{rem:sobolev-embeddings}. As we have discussed above, in this scenario, the $\setA^s_{\rho,q}$-scales can be equivalently defined using $\osckp{m}{2}$ for any $m>s$.

\begin{theorem}[Discrete stability]
\label{thm:stab-estimate-discr}
Let $\setA \in \{\setB, \setF\}$. Assume that the indices $\rho,q \in [1,\infty]$, and $s \in (0,k)$ satisfy \eqref{eq:BesovIntoL2} (with $r=2$ if $\setA=\setB$), or \eqref{eq:TriebelIntoL2} (again with $r = 2$) and $\rho < \infty$ if $\setA = \setF$. Let $\Lambda >1$ be as in Theorem~\ref{thm:osc-estimate}. For $u \in \setW^{1,2}_0(\Omega)$, every $x_0 \in \Omega$, $R>\Lambda h$, and every ball $B \coloneqq B_R(x_0)$ with $2B \coloneqq B_{2R}(x_0) \Subset \Omega$ we have
  \begin{equation}
  \label{eq:LocEstimate}
    \abs{\nabla R_h u}_{\setA^s_{\rho,q,h}(B)} \lesssim \abs{\nabla u}_{\setA^s_{\rho,q}(2B)} +
    R^{\frac{n}{\rho}-s}
    \inf_{p \in (\mathcal{P}_{k-1})^n} \Big( \dashint_{2B} \abs{\nabla R_h u -p}^2 \diff{x} \Big)^{\frac{1}2},
  \end{equation}
  where $R_h u\in V_h$ is defined in \eqref{eq:RitzProjection}, and the hidden constant is independent of $h$.
  Correspondingly, if we further restrict $s \in (0,\sigma)$, for $\sigma$ as in Theorem~\ref{thm:osc-estimate}, then, with $R = \diameter \Omega$,
  \begin{equation}
  \label{eq:GlobEstimate}
    \abs{\nabla R_h u}_{\setA^s_{\rho,q,h}(\Omega)} \lesssim \abs{\nabla u}_{\setA^s_{\rho,q}(\Omega)} +
    R^{\frac{n}{\rho}-s}
    \inf_{ \xi \in \RRn } \Big( \dashint_{\Omega} \abs{\nabla R_h u -\xi}^2 \diff{x} \Big)^{\frac{1}2}.
  \end{equation}
\end{theorem}

\begin{remark}[Improvements]
  \label{rem:improve-r}
  For simplicity, Theorem~\ref{thm:stab-estimate-discr} is stated with $r=2$, i.e., using $\oscp{2}$. By replacing Theorem~\ref{thm:osc-estimate} in the proof with Theorem~\ref{thm:osc-estimate-ii} those results can be extended to $r=\tau$, using $\oscp{\tau}$, for $\tau \in (1,2)$.
\end{remark}

\begin{proof}[Proof of Theorem~\ref{thm:stab-estimate-discr}]
  To avoid unnecessary repetition, as we have done before, we only carry out the proof for $\setA = \setF$. The Besov space case follows similarly by simply reversing the order of integration. In addition, we only consider the case $q<\infty$. The remaining cases can be proved similarly.
  
  Let $u_h \coloneqq R_h u$, and let $x \in B$ and $t \in (0,R)$. By Jensen's inequality and since $|B \cap B_t(x)| \eqsim |B_t(x)|$, see \eqref{eq:BiteOffEstimate}, we have that
  \[
    \osc^{B} \nabla u_h(x,t) \leq \oscp{2}^{B} \nabla u_h(x,t) \lesssim \oscp{2} \nabla u_h(x,t).
  \]
  Therefore, we may estimate
  \begin{align*}
      |\nabla u_h|_{\setF^s_{\rho,q,h}(B)}
      &\lesssim \bigg(\int_{B} \bigg( \int_{\theta\Lambda h}^R \Big( \frac{\oscp{2} \nabla u_h(x,t)}{t^s}\Big)^q \frac{ \diff{t} }{t}\bigg)^{\frac{\rho}{q}} \diff{x}\bigg)^{1/\rho} \\
      &\leq \bigg(\int_{B} \bigg( \int_{\theta\Lambda h}^{\theta R} \Big( \frac{\oscp{2} \nabla u_h(x,t)}{t^s}\Big)^q \frac{ \diff{t} }{t}\bigg)^{\frac{\rho}{q}} \diff{x}\bigg)^{1/\rho}
      \\
      &+ \bigg(\int_{B} \bigg( \int_{\theta R}^R \Big( \frac{\oscp{2} \nabla u_h(x,t)}{t^s}\Big)^q \frac{ \diff{t} }{t}\bigg)^{\frac{\rho}{q}} \diff{x}\bigg)^{1/\rho}
      \eqqcolon \mathrm{I} + \mathrm{II},
  \end{align*}
  and bound each term separately.

  For the first term, $\mathrm{I}$, we rescale the inner integral and apply Theorem~\ref{thm:osc-estimate} to arrive at
  \begin{align*}
    \mathrm{I}
    &= \bigg(\int_{B} \bigg( \int_{\theta\Lambda h}^{\theta R} \Big( \frac{\oscp{2} \nabla u_h(x,t)}{t^s}\Big)^q \frac{ \diff{t} }{t}\bigg)^{\frac{\rho}{q}} \diff{x}\bigg)^{1/\rho} \\
    &= \bigg(\int_{B} \bigg( \int_{\Lambda h}^R \Big( \frac{\oscp{2} \nabla u_h(x,\theta t)}{(\theta t)^s}\Big)^q \frac{\diff{t} }{t}\bigg)^{\frac{\rho}{q}} \diff{x}\bigg)^{1/\rho} \\
    &\leq \theta^{\beta-s}\bigg(\int_{B} \bigg( \int_{\Lambda h}^R \Big( \frac{\oscp{2} \nabla u_h(x,t)}{t^s}\Big)^q \frac{\diff{t}}{t}\bigg)^{\frac{\rho}{q}} \diff{x}\bigg)^{1/\rho} \\
    &+ c_\beta \theta^{-s} \bigg(\int_{B} \bigg( \int_{\Lambda h}^R \Big( \frac{\oscp{2} \nabla u(x,t)}{t^s}\Big)^q \frac{\diff{t}}{t}\bigg)^{\frac{\rho}{q}} \diff{x}\bigg)^{1/\rho} \eqqcolon \mathrm{I}_a + \mathrm{I}_b.
  \end{align*}
  Next, we write
  \begin{align*}
    \mathrm{I}_a &\leq \theta^{\beta-s}\bigg(\int_{B} \bigg( \int_{\Lambda h}^{\theta R} \Big( \frac{\oscp{2} \nabla u_h(x,t)}{t^s}\Big)^q \frac{\diff{t}}{t}\bigg)^{\frac{\rho}{q}} \diff{x}\bigg)^{1/\rho} \\
                 &+ \theta^{\beta-s}\bigg(\int_{B} \bigg( \int_{\theta R}^R \Big( \frac{\oscp{2} \nabla u_h(x,t)}{t^s}\Big)^q \frac{\diff{t}}{t}\bigg)^{\frac{\rho}{q}} \diff{x} \bigg)^{1/\rho} \\
    &\leq C_\alpha \mathrm{I} + C_\beta \theta^{\beta-s} \mathrm{II},
  \end{align*}
  where, since $\theta^{\beta-s} <1$, we have that $C_\alpha <1$ and that term can be reabsorbed on the left hand side of this inequality. In addition, we use that, for $x \in B$ and $t <R$, $\oscp{2} \nabla u(x,t) = \oscp{2}^{2B} \nabla u(x,t)$. Combining these leads to
  \begin{align*}
    \mathrm{I}
    &\lesssim \bigg(\int_{B} \bigg( \int_{\theta R}^R \Big( \frac{\oscp{2} \nabla u_h(x,t)}{t^s}\Big)^q \frac{ \diff{t} }{t}\bigg)^{\frac{\rho}{q}} \diff{x}\bigg)^{1/\rho} \\
    &+ \bigg(\int_{2B} \bigg( \int_{\Lambda h}^{R} \Big( \frac{\oscp{2}^{2B} \nabla u(x,t)}{t^s}\Big)^q \frac{\diff{t} }{t}\bigg)^{\frac{\rho}{q}} \diff{x}\bigg)^{1/\rho} \\
    &\leq \mathrm{II} + \abs{\nabla u}_{\setF^s_{\rho,q}(2B)},
  \end{align*}
  where in the last step we used that, since the indices satisfy \eqref{eq:TriebelIntoL2}, the $\setF^s_{\rho,q}$-seminorm can be equivalently characterized by $\oscp{2}$.
  For $x \in B$ and $t \in (\theta R,R)$ there holds $\abs{B_t(x)} \eqsim \abs{2B}$, see again \eqref{eq:BiteOffEstimate}, and therefore
  \begin{align*}
  \oscp{2} \nabla u_h(x,t) \lesssim \inf_{p \in (\mathcal{P}_{k-1})^n} \bigg( \dashint_{2B} \abs{\nabla u_h -p}^2 \diff{x} \bigg)^{\frac12}.
  \end{align*}
  Hence,
  \begin{align*}
    \mathrm{II}^\rho
    &\lesssim \int_{B} \bigg( \int_{\theta R}^R \Big( \frac{1}{t^s}\Big)^q \frac{\diff{t} }{t}\bigg)^{\frac{\rho}{q}} \diff{x}
    \inf_{p \in \mathcal{P}_{k-1}} \bigg( \dashint_{2B} \abs{\nabla u_h -p}^2 \diff{x} \bigg)^{\frac{\rho}2} \\
    &\lesssim 
    R^{n-s\rho} 
    \inf_{p \in \mathcal{P}_{k-1}} \bigg( \dashint_{2B} \abs{\nabla u_h -p}^2 \diff{x} \bigg)^{\frac{\rho}2}.
  \end{align*}
  This finishes the proof for $B \Subset 2B \Subset\Omega$.
  
  The case of up to the boundary estimates is proven by straightforward modifications.
\end{proof}

As a simple corollary we can combine this with Lemma \ref{lem:norm-equivalence-no-trace}.

\begin{corollary}[Uniform stability]
  Under the conditions of Theorem \ref{thm:stab-estimate-discr} assume, in addition, that we have $s\rho < 1$ and that $\nabla u \in \setA^s_{\rho,q}(\Omega)$. Furthermore, let $D \subset B$ be a discrete subdomain. Then we have
  \[
    \abs{\nabla R_h u}_{\setA^s_{\rho,q}(D)} \lesssim \abs{\nabla u}_{\setA^s_{\rho,q}(2B)} +
    R^{\frac{n}{\rho}-s}
    \inf_{p \in \mathcal{P}_k} \Big( \dashint_{2B} \abs{\nabla (R_h u -p)}^2 \diff{x} \Big)^{\frac12}.
  \]
  Moreover, if we  restrict $s \in (0,\sigma)$, then
  \[
    \abs{\nabla R_h u}_{\setA^s_{\rho,q}(\Omega)} \lesssim \abs{\nabla u}_{\setA^s_{\rho,q}(\Omega)} +
    R^{\frac{n}{\rho}-s}
    \inf_{\xi \in \RRn} \Big( \dashint_{\Omega} \abs{\nabla u -\xi}^2 \diff{x} \Big)^{\frac12}. 
  \]
  Here $\sigma \in (0,1]$ is, again, the domain dependent constant introduced in Theorem~\ref{thm:osc-estimate}.
\end{corollary}
\begin{proof}
  It suffices to recall that, for the given choice of indices, we have $\mathcal{L}^0_{k-1}(\mathcal{T}_h) \subseteq \setA^s_{\rho,q}(\Omega)$. We then combine the estimates of Theorem~\ref{thm:stab-estimate-discr} with Lemma~\ref{lem:norm-equivalence-no-trace} and \ref{lem:inf_equiv}.
\end{proof}

\begin{remark}[Particular cases]
\label{rem:SpecCasesWeRecover}
  Notice that this result generalizes \cite[(7.5)]{MR1052086}, see also \cite[(2.7)]{MR2461254}, where this result is proved for $\setA = \setB$, and $\rho = q = 2$.
\end{remark}

For the smooth companion, introduced in Section~\ref{sub:SmoothCompation}, we have a similar stability result.

\begin{corollary}[Smooth companion]
\label{cor:SmoothCompanionStab}
In the setting of Theorem~\ref{thm:stab-estimate-discr} assume, in addition, that $\nabla u \in\setA^s_{\rho,q}(\Omega)$. Let $D \subseteq \Omega$ be a discrete subdomain such that $D^+ \subset B$, where $D^+$ is as in Lemma~\ref{proplem:stab-smoothing}. Then we have
  \[
    \abs{\mathcal{S}_h \nabla R_h u}_{\setA^s_{\rho,q}(D)} \lesssim \abs{\nabla u}_{\setA^s_{\rho,q}(2B)} +
     R^{\frac{n}{\rho}-s}
    \inf_{p \in \mathcal{P}_k} \Big( \dashint_{2B} \abs{\nabla (R_h u -p)}^2 \diff{x} \Big)^{\frac12}.
  \]
  As usual, by restricting $s \in (0,\sigma)$ and the polynomial degree to $1$, this estimate is valid up to the boundary.
\end{corollary}
\begin{proof}
  The proof is a simple consequence of Theorem \ref{thm:stab-estimate-discr} in combination with the stability of the smooth companion, established in Proposition~\ref{proplem:stab-smoothing}.
\end{proof}

To conclude our discussion on higher order regularity transfer we remark that, from the stated stability results, optimal orders of convergence in the spaces $\setA^s_{\rho,q}(\Omega)$, provided $s \rho <1$ and $s < \sigma$, can be obtained. To obtain convergence in higher order smoothness spaces we make use of the smooth companion $\mathcal{S}_h \nabla R_h u$.
The following results bear resemblance to so-called local estimates for the Ritz projection, see \cite{MR1115238}.

\begin{corollary}[Convergence rates]
  \label{cor:convergence-rates}
  Under the conditions of Theorem~\ref{thm:stab-estimate-discr} assume, in addition, that $\nabla u \in\setA^{s}_{\rho,q}(\Omega)$. Let $D \subseteq \Omega$ be a discrete subdomain such that $D^+ \subset B$, where $D^+$ is as in Lemma~\ref{proplem:stab-smoothing}. For $s_1 \in (0 , s)$ there holds
  \[
    \abs{\nabla u-\mathcal{S}_h \nabla R_h u}_{\setA^{s_1}_{\rho,q}(D)} \lesssim h^{s-s_1}\abs{\nabla u}_{\setA^{s}_{\rho,q}(2B)} + R^{\frac{n}{\rho}-s_1-\frac{n}2}\norm{\nabla u - \nabla R_h u}_{\setL^2(2B)}.
  \]
  As usual, by restricting $0 < s_1 < \sigma$ this estimate can be made global. Here $\sigma \in (0,1]$ is again the domain depending constant as in Theorem~\ref{thm:osc-estimate}.
\end{corollary}
\begin{proof}
  To shorten our discussion and to avoid repetition, we shall only prove the global estimate. The remaining cases follow from small modifications. We write $u_h \coloneqq R_h u$ and $\abs{\cdot}_{\setA^\gamma_{\rho,q;I}}$ for the seminorm with the scale variable restricted to $I$.
  Using linearity of $\mathcal{S}_h$ and the triangle inequality we get
  \begin{align*}
    \abs{\nabla u-\mathcal{S}_h \nabla u_h}_{\setA^{s_1}_{\rho,q}(\Omega)} &\leq
    \abs{\nabla u-\mathcal{S}_h \nabla u}_{\setA^{s_1}_{\rho,q}(\Omega)} + \abs{ \mathcal{S}_h \nabla  u-\mathcal{S}_h \nabla u_h}_{\setA^{s_1}_{\rho,q}(\Omega)}
    \\
    &\eqqcolon \mathrm{I} + \mathrm{II}.
  \end{align*}
  We further split the first term as
  \begin{align*}
    \mathrm{I}
    &\lesssim \abs{\nabla u -\mathcal{S}_h \nabla u}_{\setA^{s_1}_{\rho,q;(0, h)}(\Omega)}
    + \abs{\nabla u -\mathcal{S}_h \nabla u}_{\setA^{s_1}_{\rho,q,h}(\Omega)}
    \eqqcolon \mathrm{I}_a + \mathrm{I}_b.
  \end{align*}
  Using the triangle inequality and carefully revisiting the proof of Proposition \ref{proplem:stab-smoothing}, there holds
  \begin{align*}
    \mathrm{I}_a
    &\lesssim \abs{\nabla u}_{\setA^{s_1}_{\rho,q;(0,3h)}(\Omega)}
    \lesssim h^{s-s_1}\abs{\nabla u}_{\setA^{s}_{\rho,q;(0,3h)}(\Omega)}
    \leq h^{s-s_1}\abs{\nabla u}_{\setA^{s}_{\rho,q}(\Omega)}.
  \end{align*}
  For the term $\mathrm{I}_b$, we apply Lemma \ref{lem:smooth_comp_approx} to estimate
  \begin{align*}
    \mathrm{I}_b
    &= \abs{\nabla u-\mathcal{S}_h \nabla u}_{\setA^{s_1}_{\rho,q,h}(\Omega)}
    \lesssim h^{-s_1} \norm{\osc \nabla u(\cdot,2 h)}_{\setL^\rho(\Omega)}
    \lesssim h^{s-s_1} \abs{\nabla u}_{\setA^{s}_{\rho,q;(2h,3h)}(\Omega)} \\
    &\leq h^{s-s_1} \abs{\nabla u}_{\setA^{s}_{\rho,q}(\Omega)}.
  \end{align*}
Here, the second step follows from $\osc\nabla v(x,2 h)\lesssim\osc\nabla v(x,t)$ for $2 h<t<3 h$, taking the respective spatial and scale norms.

  Now for the term $\mathrm{II}$. Going back to the proof of Theorem~\ref{thm:stab-estimate-discr}, only scales down to size $h$ appear and the result can be improved to
  \begin{align}
    \label{eq:improv-disc-stab}
    \begin{aligned}
    \abs{\nabla R_h u}_{\setF^{s_1}_{\rho,q,h}(\Omega)}
    &\lesssim
    \Biggabs{\int_\Omega \Bigabs{ \int_h^R \bigg(\frac{\oscp{2} \nabla u(x,t)}{t^{s_1}}\bigg)^q \frac{\diff{t}}{t} }^{\frac{\rho}{q}} \diff{x}}^{1/\rho} \\
    &\quad + R^{\frac{n}{\rho}-s_1} \inf_{ \xi \in \RRn } \Big( \dashint_{\Omega} \abs{\nabla R_h u -\xi}^2 \diff{x} \Big)^{\frac{1}2},
    \end{aligned}
  \end{align}
  same for the Besov case.
  So for $\rho < \infty$, after fixing some $s_0 \in (0,\min\{s_1,\frac1{2\rho}\})$, we may trivially estimate the first term on the right by $h^{s_0-s_1} \abs{\cdot}_{\setA^{s_0}_{\rho,q}(\Omega)}$. We apply linearity of the smooth companion, Proposition \ref{proplem:stab-smoothing}, and linearity of the Ritz projection to get
  \begin{align*}
    \mathrm{II} &= \abs{\mathcal{S}_h \nabla (u-u_h)}_{\setA^{s_1}_{\rho,q}(\Omega)}
  \lesssim \abs{\nabla (u-u_h)}_{\setA^{s_1}_{\rho,q,h}(\Omega)} \\
  &= \abs{\nabla (u-\PiSZ u-(u_h-\PiSZ u))}_{\setA^{s_1}_{\rho,q,h}(\Omega)} \\
  &\leq \abs{\nabla (u-\PiSZ u)}_{\setA^{s_1}_{\rho,q,h}(\Omega)} + \abs{\nabla R_h (u-\PiSZ u)}_{\setA^{s_1}_{\rho,q,h}(\Omega)}.
  \end{align*}
  Another reverse type inequality gives \begin{align*}
    \abs{\nabla (u-\PiSZ u)}_{\setA^{s_1}_{\rho,q,h}(\Omega)} &\lesssim h^{s_0-s_1} \abs{\nabla (u-\PiSZ u)}_{\setA^{s_0}_{\rho,q,h}(\Omega)} \\
    &\leq h^{s_0-s_1} \abs{\nabla (u-\PiSZ u)}_{\setA^{s_0}_{\rho,q}(\Omega)}.
  \end{align*}  Therefore, applying~\eqref{eq:improv-disc-stab}, and the mentioned reverse type inequalities, leads to
  \begin{align*}
    \mathrm{II}
  &\lesssim h^{s_0-s_1} \abs{\nabla (u-\PiSZ u)}_{\setA^{s_0}_{\rho,q}(\Omega)} + R^{\frac{n}{\rho}-s_1-\frac{n}{2}} \norm{\nabla (u_h-\PiSZ u)}_{\setL^2(\Omega)}.
\end{align*}
Stability of the Scott--Zhang projection yields
\[
  \norm{\nabla (u_h-\PiSZ u)}_{\setL^2(\Omega)}
  = \norm{\nabla \PiSZ(u_h- u)}_{\setL^2(\Omega)}
  \lesssim \norm{\nabla (u_h-u)}_{\setL^2(\Omega)}.
\]
Finally, since $s_0 \rho < \frac12$, by Theorem~\ref{thm:sz-bf-approximation}, we get
\begin{align*}
  h^{s_0-s_1} \abs{\nabla (u-\PiSZ u)}_{\setA^{s_0}_{\rho,q}(\Omega)}
  &\lesssim
  h^{s-s_1} \abs{\nabla u}_{\setA^{s}_{\rho,q}(\Omega)}.
\end{align*}
Collecting all estimates finishes the proof. In the Besov case for $\rho=\infty$, some modifications are needed as no $s_0$ exists. Instead, we use
\begin{align*}
  \norm{\oscp{2} \nabla (u-\PiSZ u)(\cdot,t)}_{\setL^\infty(\Omega)}
  &\lesssim
  \norm{\oscp{2} \nabla u(\cdot,2h)}_{\setL^\infty(\Omega)}, \qquad t \geq h.
\end{align*}
Starting from the use of~\eqref{eq:improv-disc-stab}, we then continue similarly to the proof of~\ref{thm:sz-bf-approximation}, using $\oscp{2}$-oscillations instead. This is justified, as~\eqref{eq:BesovIntoL2} is satisfied.
\end{proof}

\appendix
\section{Polynomial projection}

In the context of Section~\ref{sec:Oscillation} here we construct, given a sufficiently large ball, a projection from $\mathcal{P}_k$ onto the subspace of $V_h$ consisting of functions supported on the given ball. 
\begin{lemma}[Approximate polynomial]
  \label{lem:poly-approx}
  There exists $\Lambda > 1$ such that, if $R \geq \Lambda h$, then for every ball $B_R$ centered in $\Omega$ and polynomial $p \in \mathcal{P}_k$, there exists $p_h \in \mathcal{L}^1_k(\mathcal{T}_h) \cap \setW^{1,2}_0(B^\Omega_R)$ such that
  \begin{equation}
    \begin{aligned}
      \int_{B^\Omega_R} p p_h \diff{x}
      &\eqsim \int_{B^\Omega_R} \abs{p}^2 \diff{x}, \\
      R^2 \int_{B^\Omega_R} \abs{\nabla p_h}^2 \diff{x}
      &\lesssim \int_{B^\Omega_R} \abs{p}^2 \diff{x}.
    \end{aligned}
  \end{equation}
\end{lemma}
\begin{proof}
  We state the proof for a ball $B_R \Subset \Omega$; the adaptation to balls centered on the boundary follows by small modifications. As in the proof of Lemma~\ref{lem:lower-order-estimate} we use a cut-off function $\tau \in \setC^1(B_R)$ such that $\indicator_{B_{R-\xi}} \leq \tau \leq \indicator_{B_R}$ and $\norm{\nabla \tau}_{\setL^\infty(B_R)} \lesssim \xi^{-1}$ for some $\xi \in (h, \frac9{10}R) $ to be later specified. We define $A_R \coloneqq B_R \setminus B_{R-\xi}$. Now let
  \begin{align*}
    p_h \coloneqq \PiSZ(\tau p).
  \end{align*}
  By construction of $\tau$, there holds $p_h \in \setW^{1,2}_0(B_R)$. The stability of the Scott--Zhang projection and basic estimates for polynomials imply
  \begin{align*}
    \int_{B_R} \abs{\nabla p_h}^2 \diff{x}
    &\lesssim \int_{B_R} \abs{ \nabla (\tau p)}^2 \diff{x}
    \lesssim \int_{B_R} \abs{\nabla p}^2 \diff{x} + \xi^{-2} \int_{B_R} \abs{p}^2 \diff{x} \\
    &\lesssim \frac1{\xi^2} \int_{B_R} \abs{p}^2 \diff{x}.
  \end{align*}
  Finally, there holds
  \begin{align}
    \label{eq:p-ph-equivalence}
    \int_{B_R} \abs{p}^2 \diff{x}
    &= \int_{B_R} p_h p \diff{x} + \int_{B_R} (p-p_h)p \diff{x}.
  \end{align}
  By approximation properties of the Scott--Zhang operator, we have
  \begin{align*}
    \norm{p-p_h}_{\setL^2(B_R)}
    &\lesssim \norm{\tau p -\PiSZ(\tau p)}_{\setL^2(B_R)} + \norm{p-\tau p}_{\setL^2(B_R)} \\
    &\lesssim h \norm{\nabla (\tau p)}_{\setL^2(B_R)} + \norm{p}_{\setL^2(A_R)} \\
    &\lesssim \frac{h}{\xi} \norm{p}_{\setL^2(B_R)} + \frac{\xi^{1/2}}{R^{1/2}} \norm{p}_{\setL^2(B_R)}.
  \end{align*}
  We fix $\xi = cR$ for some $c>0$ small enough and fix $\Lambda > 1$ big enough such that
  \[
    \norm{p-p_h}_{\setL^2(B_R)}
    \leq \frac12 \norm{p}_{\setL^2(B_R)}.
  \]
  In view of \eqref{eq:p-ph-equivalence}, this finishes the proof.
\end{proof}

\begin{lemma}[A local interpolant]
  \label{lem:poly-proj-def}
  There exists $\Lambda > 1$ such that if $R \geq \Lambda h$, then for every ball $B_R$ centered in $\Omega$, there exists an projection $\PiP \colon \setW^{1,2}_0(B_R^\Omega) \to \mathcal{L}^1_k(\mathcal{T}_h) \cap \setW^{1,2}_0(B^\Omega_{R})$ such that
  \begin{equation}
    \begin{aligned}
    \int_{B^\Omega_R} \nabla \PiP(z) \cdot \nabla p \diff{x}
    &= \int_{B^\Omega_R} \nabla z \cdot \nabla p \diff{x}, \qquad\quad \forall p \in \mathcal{P}_k, \\
    \int_{B^\Omega_R} \abs{\nabla\PiP(z)}^2 \diff{x}
    &\lesssim \int_{B^\Omega_R} \abs{\nabla z}^2 \diff{x}.
    \end{aligned}
  \end{equation}
  Here, the hidden constant depends only on $k$ and $\Lambda$.
\end{lemma}
\begin{proof}
  Integrating by parts, we see that the first condition will follow from
  \begin{align*}
    \begin{aligned}
    \int_{B^\Omega_R} \PiP(z) p \diff{x}
    &= \int_{B^\Omega_R} z p \diff{x}, \qquad \forall p \in \mathcal{P}_{k-2}.
    \end{aligned}
  \end{align*}
  To prove that this problem has a solution let us prove an inf-sup condition. Namely,
  \[
    \inf_{p \in \mathcal{P}_{k-2}} \sup_{v_h \in \mathcal{L}^1_k(\mathcal{T}_h) \cap \setW^{1,2}_0(B^\Omega_R)}
    \frac{\int_{B^\Omega_R} v_h p \diff{x}}{\norm{\nabla v_h}_{\setL^2(B^\Omega_R)}\norm{p}_{\setW^{-1,2}(B^\Omega_R)}} \gtrsim 1.
  \]
  Fix $p \in \mathcal{P}_{k-2}$ and let $p_h \in \mathcal{L}^1_k(\mathcal{T}_h) \cap \setW^{1,2}_0(B^\Omega_R)$ be the function constructed in Lemma \ref{lem:poly-approx}. With this choice
  \begin{align*}
     \sup_{v_h \in \mathcal{L}^1_k(\mathcal{T}_h) \cap \setW^{1,2}_0(B^\Omega_R)}
    \frac{\int_{B^\Omega_R} v_h p \diff{x}}{\norm{\nabla v_h}_{\setL^2(B^\Omega_R)}}
    &\geq
    \frac{\int_{B^\Omega_R} p_h p \diff{x}}{\norm{\nabla p_h}_{\setL^2(B^\Omega_R)}}
    \gtrsim
    \frac{\int_{B^\Omega_R} \abs{p}^2 \diff{x}}{\frac1{R} \norm{p}_{\setL^2(B^\Omega_R)} }
    \gtrsim R \norm{p}_{\setL^2(B^\Omega_R)} \\
    &\gtrsim \norm{p}_{\setW^{-1,2}(B^\Omega_R)},
  \end{align*}
  where all the constants are independent of $h$. The existence of $\PiP$ with the requisite properties then follows.
\end{proof}

\section{Properties of harmonic functions}
\label{sec:harmonic-functions}

Let $\Omega \subset \RRn$ be a convex polytopal domain. Fix $x_0 \in \overline\Omega$ and, for every $r>0$, set $B_r \coloneqq B_r(x_0)$.
Here we present some properties of the function $w \in \setW^{1,2}(B_R^\Omega)$, defined by
  \begin{equation}
    \label{eq:comparison-vw-appendix}
    \begin{dcases}
      -\Delta w = 0 & \text{in } B^\Omega_R,
      \\
      w = v & \text{on } \partial B^\Omega_R,
    \end{dcases}
\end{equation}
for some $v \in \setW^{1,2}_0(\Omega)$.
Observe that most of these are standard for harmonic functions. 

\begin{lemma}
\label{lem:properties-w}
  Let the function $w \in \setW^{1,2}(B^\Omega_R)$ be defined in \eqref{eq:comparison-vw-appendix}.
  \begin{enumerate}[label={\alph{*})},leftmargin=*]
    \item \label{itm:w-min-energy}
      For every $\xi \in \RRn$ we have
    \[
      \int_{B^\Omega_R} \abs{\nabla w - \xi}^2 \diff{x} \lesssim \int_{B^\Omega_R} \abs{\nabla v -\xi}^2 \diff{x}.
    \]
    \item \label{itm:w-caccioppoli}
      For $\rho \in (0,R)$ there holds the \Caccioppoli{} estimate
    \[
      \int_{B^\Omega_\rho} \abs{\nabla w}^2 \diff{x} \lesssim \frac1{(R-\rho)^{2}} \int_{B^\Omega_R} \abs{w}^2 \diff{x}.
    \]
    \item \label{itm:w-higher-caccioppoli}
      For $\rho \in (0,R)$ there holds the \Caccioppoli{}-like estimate
    \[
      \int_{B^\Omega_\rho} \abs{\nabla^2 w}^2 \diff{x}
      \lesssim \frac1{(R-\rho)^{2}}
      \inf_{\xi \in \RRn}
      \int_{B^\Omega_R} \abs{\nabla w - \xi}^2 \diff{x}.
    \]
    \item \label{itm:w-interior-regularity}
    If $B_R \Subset \Omega$ and $\rho \in (0,R)$ we have, for any $m \geq 1$, the estimate
    \[
      \inf_{q \in \mathcal{P}_m} \dashint_{B^\Omega_\rho} \abs{\nabla (w - q)}^2 \diff{x}
      \lesssim \frac{\rho^{2m}}{R^{2m}}
      \inf_{\xi \in \RRn}
      \dashint_{B^\Omega_R}
      \abs{\nabla w - \xi}^2 \diff{x}.
    \]
    \item \label{itm:w-boundary-regularity}
    There exists $\sigma \in (0,1]$ only depending on the geometry of the domain such that, for $\rho \in (0, \tfrac12R)$,
    \[
      \inf_{\xi \in \RRn} \dashint_{B^\Omega_\rho} \abs{\nabla w - \xi}^2 \diff{x}
      \lesssim \frac{\rho^{2\sigma}}{R^{2\sigma}}
      \inf_{\xi \in \RRn} \dashint_{B^\Omega_R}
      \abs{\nabla w -\xi}^2 \diff{x}.
    \]
  \end{enumerate}
\end{lemma}
\begin{proof}
  Estimates \ref{itm:w-min-energy}, \ref{itm:w-caccioppoli}, and \ref{itm:w-interior-regularity} are standard properties of harmonic functions. We refer to, for instance, \cite[Lemma 4.2]{MR1677478} for similar results. We provide short proofs for~\ref{itm:w-higher-caccioppoli} and~\ref{itm:w-boundary-regularity}. Starting with~\ref{itm:w-boundary-regularity}, we consider several cases.
  
  We first assume that $x_0 \in \partial\Omega$ lies on a face of dimension at most $n-2$. By the geometry of $\Omega$, there exist two hyperplanes, whose intersection contains this face, such that $\Omega$ is contained in the resulting dihedral wedge $W$.
  Let $z$ be a  harmonic function on $W$ such that $z_{|\partial W} = 0$. Since $W$ is a cone, in polar coordinates there holds the harmonic expansion
  \begin{align*}
    z(r,\theta) = \sum_{i = 1}^\infty a_i r^{\lambda_i} \phi_i(\theta).
  \end{align*}
  Here, the $\phi_i$ are the eigenfunctions to the Laplace--Beltrami operator on the unit sphere, with eigenvalues $\mu_i > 0$. They form an orthonormal basis of $\setL^2$ on the sphere. Furthermore, there holds $\lambda_i(\lambda_i+n-2) = \mu_i$. An analysis of this spectrum gives $\lambda_i \geq \frac{\pi}{\omega}$ for all $i \in \setN$, where $\omega$ is the angle of $W$, see \cite[Chapter 6, Proposition 18.6']{MR961439}. Because of the convexity of $\Omega$, we can conclude $\lambda_i > 1$, and therefore
  \begin{align}
    \label{eq:dauge}
    \dashint_{W \cap B_\rho} \abs{z}^2 \diff{x}
    &\lesssim
    \frac{\rho^{2(\sigma+1)}}{R^{2(\sigma+1)}} \dashint_{W \cap B_R} \abs{z}^2 \diff{x},
  \end{align}
  for some $\sigma>0$ that depends on the angle of $W$.
  Since $\Omega$ gives a finite number of such angles, we choose the worst value of $\sigma$ to obtain a universal bound. If $B_R^W \coloneqq W \cap B_R \subset \Omega$, we can apply this to $w$. Otherwise, we have to make a comparison estimate. Let $z$ solve
\begin{equation}
  \label{eq:comparison-wz}
    \begin{dcases}
      -\Delta z = 0, & \text{in } B^W_R,
      \\
      z = \abs{w} & \text{on } \partial B^W_R \cap \overline{\Omega},
      \\
      z = 0 & \text{on } \partial B^W_R \setminus \overline{\Omega}.
    \end{dcases}
\end{equation}
As $w_{|\partial\Omega} = 0$, those boundary values are a valid $\setW^{1,2}$-extension. By multiple applications of the maximum principle, there holds $\abs{w} \leq z$. Furthermore, Dirichlet's principle bounds the energy of $z$ on $B^W_R$ by that of $w$, see property~\ref{itm:w-min-energy}. \Caccioppoli{} estimate~\ref{itm:w-caccioppoli}, estimate~\eqref{eq:dauge}, and \Poincare{} with zero boundary values
therefore leads to
\begin{align*}
    \dashint_{B^\Omega_\rho} \abs{\nabla w}^2 \diff{x}
    &\lesssim
    \frac1{\rho^2} \dashint_{B^\Omega_{2\rho}} \abs{w}^2 \diff{x}
    \leq
    \frac1{\rho^2}\dashint_{B^\Omega_{2\rho}} \abs{z}^2 \diff{x}
    \lesssim
    \frac{\rho^{2\sigma}}{R^{2(\sigma+1)}}\dashint_{B^W_{R}} \abs{z}^2 \diff{x} \\
    &\lesssim
    \frac{\rho^{2\sigma}}{R^{2\sigma}}\dashint_{B^W_{R}} \abs{\nabla z}^2 \diff{x}
    \lesssim
    \frac{\rho^{2\sigma}}{R^{2\sigma}}\dashint_{B^W_{R}} \abs{\nabla w}^2 \diff{x}
    =
    \frac{\rho^{2\sigma}}{R^{2\sigma}}\dashint_{B^\Omega_{R}} \abs{\nabla w}^2 \diff{x}
\end{align*}
for $\rho < \tfrac12 R$.
Now finally, since there is no linear function that is zero on both hyperplanes making up $W$, a combination of trace inequality and \Poincare{} proves
\begin{align}
  \label{eq:boundary-sub-polynom}
  \xi^2
  \lesssim \dashint_{B^\Omega_R} \abs{\nabla w- \xi}^2 \diff{x},
\end{align}
where the constant again depends on the angle, of which there are finitely many.
Therefore,
\begin{align}
  \label{eq:boundary-best-approx}
  \dashint_{B^\Omega_{R}} \abs{\nabla w}^2 \diff{x}
  \lesssim
  \inf_{\xi \in \RRn}\dashint_{B^\Omega_{R}} \abs{\nabla w-\xi}^2 \diff{x}.
\end{align}
This finishes the case of balls centered on an $(n-2)$-dimensional face.

For balls $B_R \Subset \Omega$, \ref{itm:w-boundary-regularity} is just \ref{itm:w-interior-regularity}. The last special case we want to look at is the case of balls centered on an $(n-1)$-dimensional face that does not intersect another face. We argue by odd extension. The normal part of $\xi$ preserves the zero boundary condition, while the tangential part is handled with a similar trace argument as in~\eqref{eq:boundary-sub-polynom}.

Finally, we consider an arbitrary center of the balls. For sets of similar size, the estimate follows trivially by inclusion. Depending on the position of the center and size of the balls, we estimate by a slightly larger ball on the boundary, use the established estimate, and go back to the original center. One crucial ingredient here is the existence of a constant $c_\Omega$ such that, if $x_0$ has distance $d$ from all $(n-2)$-dimensional faces, the ball $B_{c_\Omega d}(x_0)$ intersects at most one $(n-1)$-dimensional face. A finite chaining of the cases of interior ball, ball on $(n-1)$-dimensional face, and wedge case covers all cases.

Now for~\ref{itm:w-higher-caccioppoli}. We take a suitable cut-off function $\eta$ with $\indicator_{B^\Omega_{\frac12 R}} \leq \eta \leq \indicator_{B^\Omega_R}$. By product rule, we can compute
\begin{align*}
  \int_{B^\Omega_R} \eta^2 \abs{\nabla^2 w}^2 \diff{x}
  &\lesssim \int_{B^\Omega_R} \abs{\nabla^2 (\eta w)}^2 \diff{x}
  + \int_{B^\Omega_R} \abs{\nabla w}^2 \abs{\nabla \eta}^2 \diff{x}
  + \int_{B^\Omega_R} \abs{ \nabla^2 \eta}^2 \abs{w}^2 \diff{x}, \\
  \int_{B^\Omega_R} \abs{\Delta (\eta w)}^2 \diff{x}
  &\lesssim \int_{B^\Omega_R} \eta^2 \abs{\Delta w}^2 \diff{x}
  + \int_{B^\Omega_R} \abs{\nabla w}^2 \abs{\nabla \eta}^2 \diff{x}
  + \int_{B^\Omega_R} \abs{\Delta \eta}^2 \abs{w}^2 \diff{x}.
\end{align*}
We choose $\eta$ such that $\abs{\nabla^i \eta} \lesssim R^{-i}$ for $i = 1,2$. Combining properties of $w$ and $\eta$, the product $\eta w$ has zero boundary values on $B^\Omega_R$. Because of convexity, there also holds the regularity result
\begin{align*}
  \int_{B^\Omega_R} \abs{\nabla^2 (\eta w)}^2 \diff{x}
  &\lesssim \int_{B^\Omega_R} \abs{\Delta (\eta w)}^2 \diff{x}.
\end{align*}
see~\cite[Theorem 3.2.1.2]{MR775683}. Collecting all of those, for $\rho \in (0,\frac12 R)$, there holds
\begin{align*}
  \int_{B^\Omega_\rho} \abs{\nabla^2 w}^2 \diff{x}
  &\lesssim
  \frac1{R^2} \int_{B^\Omega_R} \abs{\nabla w}^2 \diff{x}
  + \frac1{R^4} \int_{B^\Omega_R} \abs{w}^2 \diff{x}.
\end{align*}
Now we again consider the cases of an interior ball, ball on a single face, and ball in a wedge. For the interior case, we may replace $w$ by $w-\xi^\top x$ and apply \Poincare{} for mean values. In a wedge, the zero boundary values of $w$ directly allow for \Poincare{}, while~\eqref{eq:boundary-best-approx} allows us to introduce $\xi$. On a flat face, we combine those approaches as above to deal with the normal and the tangential part of $\xi$. Finally, the same geometric arguments as above allow us to conclude the result for $x_0$ and all $\rho \in (0,cR)$. Here, $c$ depends on the geometry of $\Omega$. A covering argument for larger $\rho$ finishes the proof.
\end{proof}

\section{A continuous sharp-maximal estimate}
\label{sec:cont-sharp-maxim}

The pointwise estimate of Theorem~\ref{thm:pointwise-sharp-maximal-ii} can also be stated for $M^\sharp_\tau[\nabla u]$ and $M^\sharp_\tau[F]$. This is in the same spirit as the classical whole space estimate
\begin{align*}
  M^\sharp[\nabla u]
  \lesssim M_\tau[F],
\end{align*}
see~\cite[Theorem 3.4.9]{Grafakos2014}.
We briefly record the corresponding estimate on convex polytopal
domains.

\begin{lemma}[Continuous sharp-maximal estimate]
\label{lem:continuous-sharp-maximal}
Let $\Omega$ be a convex polytopal domain and let
$1<\tau\le 2$.
Suppose that $u\in \setW^{1,\tau}_0(\Omega)$ solves
$-\Delta u=\divergence F$ with $F\in \setL^\tau(\Omega)$.
Then, for almost every $x\in\Omega$,
\begin{align*}
  M^\sharp_\tau[\nabla u](x)
  \lesssim M^\sharp_\tau[F](x).
\end{align*}
\end{lemma}
\begin{proof}
By \cite[Section 3]{MR1156467}, the estimate
\begin{align*}
  \norm{\nabla u}_{\setL^\tau(D)}
  \lesssim \norm{F}_{\setL^\tau(D)}
\end{align*}
holds for $-\Delta u=\divergence F$ with zero Dirichlet data
on any bounded convex domain $D$, with a constant
independent of $D$.

Fix $x \in \Omega$ and $B_R^\Omega = B_R(x)\cap\Omega$. Let $w$ be the harmonic
replacement of $u$ in $B_R^\Omega$.
Since $B_R^\Omega$ is convex and
$-\Delta(u-w)=\divergence(F-\xi)$ for every constant vector $\xi$,
we obtain
\begin{align*}
  \bigg(\dashint_{B_R^\Omega}|\nabla(u-w)|^\tau \diff{y} \bigg)^{1/\tau}
  \lesssim \osckp{1}{\tau}F(x,R).
\end{align*}
Combining this with the harmonic oscillation decay
of Lemma~\ref{lem:properties-w}(e), also valid for
$\setL^\tau$-oscillations, gives, for some $\sigma>0$,
\begin{align*}
  \osckp{1}{\tau}\nabla u(x,\theta R)
  \lesssim
  \theta^\sigma\osckp{1}{\tau}\nabla u(x,R)
  +\theta^{-n/\tau}\osckp{1}{\tau}F(x,R).
\end{align*}
Choosing $\theta$ sufficiently small, taking the
supremum over $R>0$, and absorbing proves the claim.
\end{proof}

\begin{remark}[A counterexample]
  \label{rem:counterexample-linfty}
  Unlike the counter part for the sharp maximal function, the estimate $M[\nabla u] \lesssim M_\tau[F]$ is false in general. As an example, consider $\Omega = B_1(0) \subset \mathbb{R}^2$. Then, $u = x_1 \log \abs{x}$ is the solution to $F = -2 \frac{x_1x}{\abs{x}^2}$. There holds $F \in L^\infty(\Omega)$, and especially $M_\tau[F] \leq 2$ for every $\tau \in (1,2]$. But $\abs{\nabla u}$ is unbounded and only in $\setBMO$, with $M[\nabla u](x) \geq \abs{\log \abs{x}} -1$.
\end{remark}

\section{Scott--Zhang approximation}
\label{sec:sz-bf-approximation}

As before, let $\PiSZ$ denote the usual Scott--Zhang projection onto the continuous Lagrange space of degree~$k$, which preserves discrete boundary data; see~\cite{MR1011446}. In this section we prove stability and approximation estimates for $\nabla\PiSZ v$. Gradient stability in the fractional Sobolev spaces $\setB^s_{\rho,\rho}(\Omega)=\setF^s_{\rho,\rho}(\Omega)$ for $1<\rho<\infty$ and $0<s<1/\rho$ is proved in~\cite[Appendix~A]{DieningStorn2026}. Related Besov stability estimates are established in~\cite[Theorem~2.2 and Corollary~3.7]{MR4238777} for the volume-averaging version of the Scott--Zhang projection, which is $\setL^1$-stable but does not in general preserve discrete boundary data.

\begin{theorem}[Stability and approximation]
  \label{thm:sz-bf-approximation}
  Let
  $\rho,q$ be as in Definition~\ref{def:besovTriebel}, and let $0<s_2\leq s<k$. For every $v\in\setW^{1,\rho}(\Omega)$, we have
  \begin{equation}
    \label{eq:disc-stab}
    \abs{\nabla(v-\PiSZ v)}_{\setA^{s_2}_{\rho,q,h}(\Omega)}
    \lesssim h^{s-s_2}\abs{\nabla v}_{\setA^s_{\rho,q,h}(\Omega)}.
  \end{equation}
  If $s_2\rho<1$, then also
  \begin{equation}
    \label{eq:cont-stab}
    \abs{\nabla(v-\PiSZ v)}_{\setA^{s_2}_{\rho,q}(\Omega)}
    \lesssim h^{s-s_2}\abs{\nabla v}_{\setA^s_{\rho,q}(\Omega)}.
  \end{equation}
  The case $s_2=s$ yields stability of $\PiSZ$ in the corresponding gradient seminorm. The constants are independent of $h$ and $v$.
\end{theorem}
\begin{proof}
  We write $\abs{\cdot}_{\setA^\gamma_{\rho,q};I}$ for the seminorm with the scale variable restricted to $I$. All scale intervals below lie in $(0,\diameter\Omega)$ by the standing convention that $4h<\diameter\Omega$.

  For every $y\in T$, we have $\omega(T)\subset B_{2h}(y)$. Subtracting a polynomial in $\mathcal P_k$, using polynomial reproduction and local gradient stability, and applying Lemma~\ref{lem:inf_equiv} with integrability exponent equal one, we obtain
  \begin{align}
    \label{eq:sz-local-short}
    \dashint_T\abs{\nabla(v-\PiSZ v)}\diff y
    &\lesssim\dashint_T\osc\nabla v(y,2h)\diff y,\\
    \label{eq:sz-osc-short}
    \osc\nabla\PiSZ v(x,t)
    &\lesssim\osc\nabla v(x,t+2h),
    \qquad t\geq h.
  \end{align}
  The proof of Lemma~\ref{lem:smooth_comp_approx}, starting from~\eqref{eq:sz-local-short}, therefore gives
  \begin{align}
    \label{eq:sz-large-short}
    \begin{aligned}
      \abs{\nabla(v-\PiSZ v)}_{\setA^{s_2}_{\rho,q,h}(\Omega)}
      &\lesssim h^{-s_2}\norm{\osc\nabla v(\cdot,2h)}_{\setL^\rho(\Omega)}\\
      &\lesssim h^{s-s_2}\abs{\nabla v}_{\setA^s_{\rho,q};(2h,3h)}
      \leq h^{s-s_2}\abs{\nabla v}_{\setA^s_{\rho,q,h}(\Omega)}.
    \end{aligned}
  \end{align}
  Here, the second step follows from $\osc\nabla v(x,2h)\lesssim\osc\nabla v(x,t)$ for $2h<t<3h$, taking the respective spatial and scale norms. This proves~\eqref{eq:disc-stab}.

  Suppose now $s_2\rho<1$. Adjusting the fixed scale ratios in the small-scale argument in the proof of Lemma~\ref{lem:norm-equivalence-no-trace} and using~\eqref{eq:sz-osc-short} with $t=h$, we obtain
  \begin{align}
    \label{eq:sz-small-short}
    \begin{aligned}
      \abs{\nabla\PiSZ v}_{\setA^{s_2}_{\rho,q};(0,h)}
      &\lesssim h^{-s_2}\norm{\osc\nabla\PiSZ v(\cdot,h)}_{\setL^\rho(\Omega)}\\
      &\lesssim h^{-s_2}\norm{\osc\nabla v(\cdot,3h)}_{\setL^\rho(\Omega)}\\
      &\lesssim h^{s-s_2}\abs{\nabla v}_{\setA^s_{\rho,q};(3h,4h)}.
    \end{aligned}
  \end{align}
  The triangle inequality and $t^{s-s_2}\leq h^{s-s_2}$ for $t<h$ now yield
  \begin{align*}
    \abs{\nabla(v-\PiSZ v)}_{\setA^{s_2}_{\rho,q};(0,h)}
    &\lesssim
      \abs{\nabla v}_{\setA^{s_2}_{\rho,q};(0,h)}
      +\abs{\nabla\PiSZ v}_{\setA^{s_2}_{\rho,q};(0,h)}\\
    &\lesssim h^{s-s_2}\abs{\nabla v}_{\setA^s_{\rho,q};(0,4h)}.
  \end{align*}
  Together with~\eqref{eq:sz-large-short}, this proves~\eqref{eq:cont-stab}. Taking $s_2=s$ and applying the triangle inequality gives stability. The argument includes $q=\infty$, with suprema in place of scale integrals.
\end{proof}

%% ------------------------------------------------------------

\section*{Acknowledgements}
The work of the first two authors was supported by the Deutsche Forschungsgemeinschaft
(German Research Foundation) grant SFB 1283/2 2021 - 317210226.
The work of AJS is partially supported by NSF grant DMS-2409918.

\section*{Statement about AI use}
The authors take full responsibility for the correctness of the mathematical arguments and results and for the accuracy, and completeness of the references. During the preparation of this work, OpenAI's ChatGPT was used to identify potentially relevant literature and to suggest approaches for adapting the auxiliary results on harmonic functions in Appendix~\ref{sec:harmonic-functions} to convex polytopal domains. The authors checked and adapted the resulting mathematical arguments and references. After a complete draft had been written, ChatGPT was also used for proofreading and language editing.

\printbibliography

@article {MR1011446,
    AUTHOR = {Scott, L. Ridgway and Zhang, Shangyou},
     TITLE = {Finite element interpolation of nonsmooth functions satisfying
              boundary conditions},
   JOURNAL = {Math. Comp.},
  FJOURNAL = {Mathematics of Computation},
    VOLUME = {54},
      YEAR = {1990},
    NUMBER = {190},
     PAGES = {483--493},
      ISSN = {0025-5718,1088-6842},
   MRCLASS = {65D05 (65N30)},
  MRNUMBER = {1011446},
MRREVIEWER = {Qian\ Li},
       DOI = {10.2307/2008497},
       URL = {https://doi.org/10.2307/2008497},
}

@article {MR1023118,
    AUTHOR = {Triebel, H.},
     TITLE = {Local approximation spaces},
   JOURNAL = {Z. Anal. Anwendungen},
  FJOURNAL = {Zeitschrift f\"ur Analysis und ihre Anwendungen},
    VOLUME = {8},
      YEAR = {1989},
    NUMBER = {3},
     PAGES = {261--288},
      ISSN = {0232-2064},
   MRCLASS = {46E35 (42B99)},
  MRNUMBER = {1023118},
MRREVIEWER = {Vladimir\ D.\ Stepanov},
       DOI = {10.4171/ZAA/353},
       URL = {https://doi.org/10.4171/ZAA/353},
}

@article {MR2017700,
    AUTHOR = {Dispa, Sophie},
     TITLE = {Intrinsic characterizations of {B}esov spaces on {L}ipschitz
              domains},
   JOURNAL = {Math. Nachr.},
  FJOURNAL = {Mathematische Nachrichten},
    VOLUME = {260},
      YEAR = {2003},
     PAGES = {21--33},
      ISSN = {0025-584X,1522-2616},
   MRCLASS = {46E35 (42B35 46F05)},
  MRNUMBER = {2017700},
       DOI = {10.1002/mana.200310101},
       URL = {https://doi.org/10.1002/mana.200310101},
}

@article {MR2247585,
    AUTHOR = {Shvartsman, Pavel},
     TITLE = {Local approximations and intrinsic characterization of spaces
              of smooth functions on regular subsets of {${\mathbb R}^n$}},
   JOURNAL = {Math. Nachr.},
  FJOURNAL = {Mathematische Nachrichten},
    VOLUME = {279},
      YEAR = {2006},
    NUMBER = {11},
     PAGES = {1212--1241},
      ISSN = {0025-584X,1522-2616},
   MRCLASS = {46E35 (26B35 41A30 46E30)},
  MRNUMBER = {2247585},
MRREVIEWER = {Dorothee\ D.\ Haroske},
       DOI = {10.1002/mana.200510418},
       URL = {https://doi.org/10.1002/mana.200510418},
}

@article {MR4801010,
    AUTHOR = {Hovemann, Marc and Weimar, Markus},
     TITLE = {Oscillations and differences in {T}riebel-{L}izorkin-{M}orrey
              spaces},
   JOURNAL = {Rev. Mat. Complut.},
  FJOURNAL = {Revista Matem\'atica Complutense},
    VOLUME = {37},
      YEAR = {2024},
    NUMBER = {3},
     PAGES = {735--782},
      ISSN = {1139-1138,1988-2807},
   MRCLASS = {46E35 (26B35 41A50 46E30)},
  MRNUMBER = {4801010},
MRREVIEWER = {Xue\ Xiu\ Zhong},
       DOI = {10.1007/s13163-024-00487-4},
       URL = {https://doi.org/10.1007/s13163-024-00487-4},
}

@incollection {MR1984176,
    AUTHOR = {Dispa, Sophie},
     TITLE = {Intrinsic descriptions using means of differences for {B}esov
              spaces on {L}ipschitz domains},
 BOOKTITLE = {Function spaces, differential operators and nonlinear analysis
              ({T}eistungen, 2001)},
     PAGES = {279--287},
 PUBLISHER = {Birkh\"auser, Basel},
      YEAR = {2003},
      ISBN = {3-7643-6935-3},
   MRCLASS = {46E35},
  MRNUMBER = {1984176},
}

@book {MR2373954,
    AUTHOR = {Brenner, Susanne C. and Scott, L. Ridgway},
     TITLE = {The mathematical theory of finite element methods},
    SERIES = {Texts in Applied Mathematics},
    VOLUME = {15},
   EDITION = {Third},
 PUBLISHER = {Springer, New York},
      YEAR = {2008},
     PAGES = {xviii+397},
      ISBN = {978-0-387-75933-3},
   MRCLASS = {65-01 (65-02)},
  MRNUMBER = {2373954},
       DOI = {10.1007/978-0-387-75934-0},
       URL = {https://doi.org/10.1007/978-0-387-75934-0},
}

@article {MR559195,
    AUTHOR = {Dupont, Todd and Scott, Ridgway},
     TITLE = {Polynomial approximation of functions in {S}obolev spaces},
   JOURNAL = {Math. Comp.},
  FJOURNAL = {Mathematics of Computation},
    VOLUME = {34},
      YEAR = {1980},
    NUMBER = {150},
     PAGES = {441--463},
      ISSN = {0025-5718,1088-6842},
   MRCLASS = {65D15 (41A10)},
  MRNUMBER = {559195},
       DOI = {10.2307/2006095},
       URL = {https://doi.org/10.2307/2006095},
}

@article {MR645661,
    AUTHOR = {Rannacher, Rolf and Scott, Ridgway},
     TITLE = {Some optimal error estimates for piecewise linear finite
              element approximations},
   JOURNAL = {Math. Comp.},
  FJOURNAL = {Mathematics of Computation},
    VOLUME = {38},
      YEAR = {1982},
    NUMBER = {158},
     PAGES = {437--445},
      ISSN = {0025-5718,1088-6842},
   MRCLASS = {65N30},
  MRNUMBER = {645661},
       DOI = {10.2307/2007280},
       URL = {https://doi.org/10.2307/2007280},
}

@article {MR2495783,
    AUTHOR = {Guzm\'an, J. and Leykekhman, D. and Rossmann, J. and Schatz,
              A. H.},
     TITLE = {H\"older estimates for {G}reen's functions on convex
              polyhedral domains and their applications to finite element
              methods},
   JOURNAL = {Numer. Math.},
  FJOURNAL = {Numerische Mathematik},
    VOLUME = {112},
      YEAR = {2009},
    NUMBER = {2},
     PAGES = {221--243},
      ISSN = {0029-599X,0945-3245},
   MRCLASS = {65N30 (65N12 65N15)},
  MRNUMBER = {2495783},
MRREVIEWER = {Ivan\ Yotov},
       DOI = {10.1007/s00211-009-0213-y},
       URL = {https://doi.org/10.1007/s00211-009-0213-y},
}

@article {MR2869035,
    AUTHOR = {Demlow, A. and Leykekhman, D. and Schatz, A. H. and Wahlbin,
              L. B.},
     TITLE = {Best approximation property in the {$W^{1}_{\infty}$} norm for
              finite element methods on graded meshes},
   JOURNAL = {Math. Comp.},
  FJOURNAL = {Mathematics of Computation},
    VOLUME = {81},
      YEAR = {2012},
    NUMBER = {278},
     PAGES = {743--764},
      ISSN = {0025-5718,1088-6842},
   MRCLASS = {65N30 (65N15 65N50)},
  MRNUMBER = {2869035},
MRREVIEWER = {Nicolae\ Pop},
       DOI = {10.1090/S0025-5718-2011-02546-9},
       URL = {https://doi.org/10.1090/S0025-5718-2011-02546-9},
}

@article {MR4060457,
    AUTHOR = {Drelichman, Irene and Dur\'an, Ricardo G. and Ojea, Ignacio},
     TITLE = {A weighted setting for the numerical approximation of the
              {P}oisson problem with singular sources},
   JOURNAL = {SIAM J. Numer. Anal.},
  FJOURNAL = {SIAM Journal on Numerical Analysis},
    VOLUME = {58},
      YEAR = {2020},
    NUMBER = {1},
     PAGES = {590--606},
      ISSN = {0036-1429,1095-7170},
   MRCLASS = {65N30 (35B45 35J25 65N15)},
  MRNUMBER = {4060457},
MRREVIEWER = {Jiguang\ Sun},
       DOI = {10.1137/18M1213105},
       URL = {https://doi.org/10.1137/18M1213105},
}

@article {MR4748207,
    AUTHOR = {Diening, Lars and Rolfes, Julian and Salgado, Abner J.},
     TITLE = {Pointwise gradient estimate of the {R}itz projection},
   JOURNAL = {SIAM J. Numer. Anal.},
  FJOURNAL = {SIAM Journal on Numerical Analysis},
    VOLUME = {62},
      YEAR = {2024},
    NUMBER = {3},
     PAGES = {1212--1225},
      ISSN = {0036-1429,1095-7170},
   MRCLASS = {65N30 (65N12 65N80)},
  MRNUMBER = {4748207},
       DOI = {10.1137/23M1571800},
       URL = {https://doi.org/10.1137/23M1571800},
}

@book {MR3024912,
    AUTHOR = {Pick, Lubo\v s{} and Kufner, Alois and John, Old\v rich and
              Fu\v c\'ik, Svatopluk},
     TITLE = {Function spaces. {V}ol. 1},
    SERIES = {De Gruyter Series in Nonlinear Analysis and Applications},
    VOLUME = {14},
   EDITION = {extended},
 PUBLISHER = {Walter de Gruyter \& Co., Berlin},
      YEAR = {2013},
     PAGES = {xvi+479},
      ISBN = {978-3-11-025041-1},
   MRCLASS = {46Exx (46-02)},
  MRNUMBER = {3024912},
MRREVIEWER = {Dumitru\ Popa},
}

@article {MR5022107,
    AUTHOR = {Drelichman, Irene and Dur\'an, Ricardo G.},
     TITLE = {Stability of the {R}itz projection in weighted {$W^{1,1}$}},
   JOURNAL = {Rev. Un. Mat. Argentina},
  FJOURNAL = {Revista de la Uni\'on Matem\'atica Argentina},
    VOLUME = {69},
      YEAR = {2026},
    NUMBER = {1},
     PAGES = {21--24},
      ISSN = {0041-6932,1669-9637},
   MRCLASS = {46E35 (42B25)},
  MRNUMBER = {5022107},
       DOI = {10.33044/revuma.5001},
       URL = {https://doi.org/10.33044/revuma.5001},
}

@article {MR1677478,
    AUTHOR = {Dolzmann, Georg},
     TITLE = {Optimal convergence for the finite element method in
              {C}ampanato spaces},
   JOURNAL = {Math. Comp.},
  FJOURNAL = {Mathematics of Computation},
    VOLUME = {68},
      YEAR = {1999},
    NUMBER = {228},
     PAGES = {1397--1427},
      ISSN = {0025-5718,1088-6842},
   MRCLASS = {65N30 (65N12)},
  MRNUMBER = {1677478},
MRREVIEWER = {Tsorng-Whay\ Pan},
       DOI = {10.1090/S0025-5718-99-01175-8},
       URL = {https://doi.org/10.1090/S0025-5718-99-01175-8},
}

@book {MR3136501,
    AUTHOR = {Jovanovi\'c, Bo\v sko S. and S\"uli, Endre},
     TITLE = {Analysis of finite difference schemes},
    SERIES = {Springer Series in Computational Mathematics},
    VOLUME = {46},
      NOTE = {For linear partial differential equations with generalized
              solutions},
 PUBLISHER = {Springer, London},
      YEAR = {2014},
     PAGES = {xiv+408},
      ISBN = {978-1-4471-5459-4},
   MRCLASS = {65-02 (35B30 35J40 35K35 35L35 65M06 65N06)},
  MRNUMBER = {3136501},
MRREVIEWER = {Enrico\ Bertolazzi},
       DOI = {10.1007/978-1-4471-5460-0},
       URL = {https://doi.org/10.1007/978-1-4471-5460-0},
}

@article {MR2461254,
    AUTHOR = {Guermond, Jean-Luc and Pasciak, Joseph E.},
     TITLE = {Stability of discrete {S}tokes operators in fractional
              {S}obolev spaces},
   JOURNAL = {J. Math. Fluid Mech.},
  FJOURNAL = {Journal of Mathematical Fluid Mechanics},
    VOLUME = {10},
      YEAR = {2008},
    NUMBER = {4},
     PAGES = {588--610},
      ISSN = {1422-6928,1422-6952},
   MRCLASS = {35Q30 (35B35 65N12 65N30 76M10)},
  MRNUMBER = {2461254},
MRREVIEWER = {Jordi\ Blasco},
       DOI = {10.1007/s00021-007-0244-z},
       URL = {https://doi.org/10.1007/s00021-007-0244-z},
}

@book {Oswald1994,
    AUTHOR = {Oswald, Peter},
     TITLE = {Multilevel finite element approximation},
    SERIES = {Teubner Skripten zur Numerik. [Teubner Scripts on Numerical
              Mathematics]},
      NOTE = {Theory and applications},
 PUBLISHER = {B. G. Teubner, Stuttgart},
      YEAR = {1994},
     PAGES = {160},
      ISBN = {3-519-02719-4},
   MRCLASS = {65N30 (65N55)},
  MRNUMBER = {1312165},
MRREVIEWER = {Weimin\ Han},
       DOI = {10.1007/978-3-322-91215-2},
       URL = {https://doi.org/10.1007/978-3-322-91215-2},
}

@article {MR4238777,
    AUTHOR = {Faustmann, Markus and Melenk, Jens Markus and Parvizi, Maryam},
     TITLE = {On the stability of {S}cott-{Z}hang type operators and
              application to multilevel preconditioning in fractional
              diffusion},
   JOURNAL = {ESAIM Math. Model. Numer. Anal.},
  FJOURNAL = {ESAIM. Mathematical Modelling and Numerical Analysis},
    VOLUME = {55},
      YEAR = {2021},
    NUMBER = {2},
     PAGES = {595--625},
      ISSN = {2822-7840,2804-7214},
   MRCLASS = {65N30 (35R11 65N22)},
  MRNUMBER = {4238777},
       DOI = {10.1051/m2an/2020079},
       URL = {https://doi.org/10.1051/m2an/2020079},
}

@incollection {MR3586796,
    AUTHOR = {Bourgain, Jean and Brezis, Haim and Mironescu, Petru},
     TITLE = {Another look at {S}obolev spaces},
 BOOKTITLE = {Optimal control and partial differential equations},
     PAGES = {439--455},
 PUBLISHER = {IOS, Amsterdam},
      YEAR = {2001},
      ISBN = {1-58603-096-5},
   MRCLASS = {46E30},
  MRNUMBER = {3586796},
}

@article {MR4654593,
    AUTHOR = {Du, Qiang and Mengesha, Tadele and Tian, Xiaochuan},
     TITLE = {{$L^p$} compactness criteria with an application to
              variational convergence of some nonlocal energy functionals},
   JOURNAL = {Math. Eng.},
  FJOURNAL = {Mathematics in Engineering},
    VOLUME = {5},
      YEAR = {2023},
    NUMBER = {6},
     PAGES = {Paper No. 097, 31},
      ISSN = {2640-3501},
   MRCLASS = {49J45 (74A70)},
  MRNUMBER = {4654593},
       DOI = {10.3934/mine.2023097},
       URL = {https://doi.org/10.3934/mine.2023097},
}

@article {MR2777246,
    AUTHOR = {Brenner, Susanne C. and Gudi, Thirupathi and Sung, Li-Yeng},
     TITLE = {A weakly over-penalized symmetric interior penalty method for
              the biharmonic problem},
   JOURNAL = {Electron. Trans. Numer. Anal.},
  FJOURNAL = {Electronic Transactions on Numerical Analysis},
    VOLUME = {37},
      YEAR = {2010},
     PAGES = {214--238},
      ISSN = {1068-9613},
   MRCLASS = {65N30 (65N15)},
  MRNUMBER = {2777246},
MRREVIEWER = {Manfred\ Dobrowolski},
}

@article {MR2755946,
    AUTHOR = {Georgoulis, Emmanuil H. and Houston, Paul and Virtanen, Juha},
     TITLE = {An {\it a~posteriori} error indicator for discontinuous
              {G}alerkin approximations of fourth-order elliptic problems},
   JOURNAL = {IMA J. Numer. Anal.},
  FJOURNAL = {IMA Journal of Numerical Analysis},
    VOLUME = {31},
      YEAR = {2011},
    NUMBER = {1},
     PAGES = {281--298},
      ISSN = {0272-4979,1464-3642},
   MRCLASS = {65N30 (65N15)},
  MRNUMBER = {2755946},
MRREVIEWER = {Manfred\ Dobrowolski},
       DOI = {10.1093/imanum/drp023},
       URL = {https://doi.org/10.1093/imanum/drp023},
}

@article {MR2034620,
    AUTHOR = {Karakashian, Ohannes A. and Pascal, Frederic},
     TITLE = {A posteriori error estimates for a discontinuous {G}alerkin
              approximation of second-order elliptic problems},
   JOURNAL = {SIAM J. Numer. Anal.},
  FJOURNAL = {SIAM Journal on Numerical Analysis},
    VOLUME = {41},
      YEAR = {2003},
    NUMBER = {6},
     PAGES = {2374--2399},
      ISSN = {0036-1429,1095-7170},
   MRCLASS = {65N15 (65N30 65N50 65N55)},
  MRNUMBER = {2034620},
MRREVIEWER = {M.\ M.\ Sussman},
       DOI = {10.1137/S0036142902405217},
       URL = {https://doi.org/10.1137/S0036142902405217},
}

@incollection {MR1115238,
    AUTHOR = {Wahlbin, Lars B.},
     TITLE = {Local behavior in finite element methods},
 BOOKTITLE = {Handbook of numerical analysis, {V}ol.\ {II}},
    SERIES = {Handb. Numer. Anal.},
    VOLUME = {II},
     PAGES = {353--522},
 PUBLISHER = {North-Holland, Amsterdam},
      YEAR = {1991},
      ISBN = {0-444-70365-9},
   MRCLASS = {65-02 (65M60 65N30)},
  MRNUMBER = {1115238},
}

@article {MR1052086,
    AUTHOR = {Bramble, James H. and Pasciak, Joseph E. and Xu, Jinchao},
     TITLE = {The analysis of multigrid algorithms with nonnested spaces or
              noninherited quadratic forms},
   JOURNAL = {Math. Comp.},
  FJOURNAL = {Mathematics of Computation},
    VOLUME = {56},
      YEAR = {1991},
    NUMBER = {193},
     PAGES = {1--34},
      ISSN = {0025-5718,1088-6842},
   MRCLASS = {65N22 (65N55)},
  MRNUMBER = {1052086},
MRREVIEWER = {Jan\ Mandel},
       DOI = {10.2307/2008527},
       URL = {https://doi.org/10.2307/2008527},
}

@book {MR1805196,
    AUTHOR = {Axler, Sheldon and Bourdon, Paul and Ramey, Wade},
     TITLE = {Harmonic function theory},
    SERIES = {Graduate Texts in Mathematics},
    VOLUME = {137},
   EDITION = {Second},
 PUBLISHER = {Springer-Verlag, New York},
      YEAR = {2001},
     PAGES = {xii+259},
      ISBN = {0-387-95218-7},
   MRCLASS = {31-01 (30-01 46Exx)},
  MRNUMBER = {1805196},
       DOI = {10.1007/978-1-4757-8137-3},
       URL = {https://doi.org/10.1007/978-1-4757-8137-3},
}

@article {MR290001,
    AUTHOR = {Brudny\u i, Ju.\ A.},
     TITLE = {Piecewise polynomial approximation and local approximations},
   JOURNAL = {Dokl. Akad. Nauk SSSR},
  FJOURNAL = {Doklady Akademii Nauk SSSR},
    VOLUME = {201},
      YEAR = {1971},
     PAGES = {16--18},
      ISSN = {0002-3264},
   MRCLASS = {41.15},
  MRNUMBER = {290001},
MRREVIEWER = {V.\ A.\ Popov},
}

@book {MR781540,
    AUTHOR = {Triebel, Hans},
     TITLE = {Theory of function spaces},
    SERIES = {Monographs in Mathematics},
    VOLUME = {78},
 PUBLISHER = {Birkh\"auser Verlag, Basel},
      YEAR = {1983},
     PAGES = {284},
      ISBN = {3-7643-1381-1},
   MRCLASS = {46Exx},
  MRNUMBER = {781540},
       DOI = {10.1007/978-3-0346-0416-1},
       URL = {https://doi.org/10.1007/978-3-0346-0416-1},
}

@book {TriebelFunctionSpacesII,
    AUTHOR = {Triebel, Hans},
     TITLE = {Theory of function spaces. {II}},
    SERIES = {Monographs in Mathematics},
    VOLUME = {84},
 PUBLISHER = {Birkh\"auser Verlag, Basel},
      YEAR = {1992},
     PAGES = {viii+370},
      ISBN = {3-7643-2639-5},
   MRCLASS = {46Exx (46-02)},
  MRNUMBER = {1163193},
MRREVIEWER = {P.\ Szeptycki},
       DOI = {10.1007/978-3-0346-0419-2},
       URL = {https://doi.org/10.1007/978-3-0346-0419-2},
}

@book {MR2250142,
    AUTHOR = {Triebel, Hans},
     TITLE = {Theory of function spaces. {III}},
    SERIES = {Monographs in Mathematics},
    VOLUME = {100},
 PUBLISHER = {Birkh\"auser Verlag, Basel},
      YEAR = {2006},
     PAGES = {xii+426},
      ISBN = {978-3-7643-7581-2},
   MRCLASS = {46E35 (46-02)},
  MRNUMBER = {2250142},
MRREVIEWER = {Lubo\v s\ Pick},
}

@article {MR1152321,
    AUTHOR = {DeVore, Ronald A. and Sharpley, Robert C.},
     TITLE = {Besov spaces on domains in {${\bf R}^d$}},
   JOURNAL = {Trans. Amer. Math. Soc.},
  FJOURNAL = {Transactions of the American Mathematical Society},
    VOLUME = {335},
      YEAR = {1993},
    NUMBER = {2},
     PAGES = {843--864},
      ISSN = {0002-9947,1088-6850},
   MRCLASS = {46E35},
  MRNUMBER = {1152321},
MRREVIEWER = {Vladimir\ D.\ Stepanov},
       DOI = {10.2307/2154408},
       URL = {https://doi.org/10.2307/2154408},
}

@book {MR2355272,
    AUTHOR = {Lai, Ming-Jun and Schumaker, Larry L.},
     TITLE = {Spline functions on triangulations},
    SERIES = {Encyclopedia of Mathematics and its Applications},
    VOLUME = {110},
 PUBLISHER = {Cambridge University Press, Cambridge},
      YEAR = {2007},
     PAGES = {xvi+592},
      ISBN = {978-0-521-87592-9},
   MRCLASS = {41-02 (41A15 41A63 65D05 65D07)},
  MRNUMBER = {2355272},
MRREVIEWER = {Peter\ Alfeld},
       DOI = {10.1017/CBO9780511721588},
       URL = {https://doi.org/10.1017/CBO9780511721588},
}

@book {MR473503,
    AUTHOR = {Simader, Christian G.},
     TITLE = {On {D}irichlet's boundary value problem},
    SERIES = {Lecture Notes in Mathematics},
    VOLUME = {Vol. 268},
 PUBLISHER = {Springer-Verlag, Berlin-New York},
      YEAR = {1972},
     PAGES = {iv+238},
   MRCLASS = {35J40},
  MRNUMBER = {473503},
MRREVIEWER = {M.\ Schechter},
}

@book {MR1962933,
    AUTHOR = {Giusti, Enrico},
     TITLE = {Direct methods in the calculus of variations},
 PUBLISHER = {World Scientific Publishing Co., Inc., River Edge, NJ},
      YEAR = {2003},
     PAGES = {viii+403},
      ISBN = {981-238-043-4},
   MRCLASS = {49-02 (35J50 49J10 49K10 49N60)},
  MRNUMBER = {1962933},
MRREVIEWER = {Giovanni\ Alberti},
       DOI = {10.1142/9789812795557},
       URL = {https://doi.org/10.1142/9789812795557},
}

@article {MR2643399,
    AUTHOR = {Diening, Lars and Růžička, Michael and
              Schumacher, Katrin},
     TITLE = {A decomposition technique for {J}ohn domains},
   JOURNAL = {Ann. Acad. Sci. Fenn. Math.},
  FJOURNAL = {Annales Academi\ae\ Scientiarum Fennic\ae. Mathematica},
    VOLUME = {35},
      YEAR = {2010},
    NUMBER = {1},
     PAGES = {87--114},
      ISSN = {1239-629X,1798-2383},
   MRCLASS = {26D10 (35Q35 46E30)},
  MRNUMBER = {2643399},
MRREVIEWER = {Shusen\ Ding},
       DOI = {10.5186/aasfm.2010.3506},
       URL = {https://doi.org/10.5186/aasfm.2010.3506},
}

@article {MR2375698,
    AUTHOR = {Lerner, Andrei K. and P\'erez, Carlos},
     TITLE = {A new characterization of the {M}uckenhoupt {$A_p$} weights
              through an extension of the {L}orentz-{S}himogaki theorem},
   JOURNAL = {Indiana Univ. Math. J.},
  FJOURNAL = {Indiana University Mathematics Journal},
    VOLUME = {56},
      YEAR = {2007},
    NUMBER = {6},
     PAGES = {2697--2722},
      ISSN = {0022-2518,1943-5258},
   MRCLASS = {42B25 (46A16)},
  MRNUMBER = {2375698},
MRREVIEWER = {Ale\v s\ Nekvinda},
       DOI = {10.1512/iumj.2007.56.3112},
       URL = {https://doi.org/10.1512/iumj.2007.56.3112},
}

@book {MR4242224,
    AUTHOR = {Ern, Alexandre and Guermond, Jean-Luc},
     TITLE = {Finite elements {I}---{A}pproximation and interpolation},
    SERIES = {Texts in Applied Mathematics},
    VOLUME = {72},
 PUBLISHER = {Springer, Cham},
      YEAR = {2021},
     PAGES = {xii+325},
      ISBN = {978-3-030-56340-0},
   MRCLASS = {65-01},
  MRNUMBER = {4242224},
       DOI = {10.1007/978-3-030-56341-7},
       URL = {https://doi.org/10.1007/978-3-030-56341-7},
}

@article {MR1156467,
    AUTHOR = {Fromm, Stephen J.},
     TITLE = {Potential space estimates for {G}reen potentials in convex
              domains},
   JOURNAL = {Proc. Amer. Math. Soc.},
  FJOURNAL = {Proceedings of the American Mathematical Society},
    VOLUME = {119},
      YEAR = {1993},
    NUMBER = {1},
     PAGES = {225--233},
      ISSN = {0002-9939,1088-6826},
   MRCLASS = {35J25 (35B65)},
  MRNUMBER = {1156467},
MRREVIEWER = {Yehuda\ Pinchover},
       DOI = {10.2307/2159846},
       URL = {https://doi.org/10.2307/2159846},
}

@article {MR1198129,
    AUTHOR = {Alkhutov, Yu.A. and Kondrat'ev, V. A.},
     TITLE = {Solvability of the {D}irichlet problem for second-order
              elliptic equations in a convex domain},
   JOURNAL = {Differentsial'nye Uravneniya},
  FJOURNAL = {Differentsial'nye Uravneniya},
    VOLUME = {28},
      YEAR = {1992},
    NUMBER = {5},
     PAGES = {806--818, 917},
      ISSN = {0374-0641},
   MRCLASS = {35J25},
  MRNUMBER = {1198129},
}

@book {MR3243734,
    AUTHOR = {Grafakos, Loukas},
     TITLE = {Classical {F}ourier analysis},
    SERIES = {Graduate Texts in Mathematics},
    VOLUME = {249},
   EDITION = {Third},
 PUBLISHER = {Springer, New York},
      YEAR = {2014},
     PAGES = {xviii+638},
      ISBN = {978-1-4939-1193-6},
   MRCLASS = {42-01 (42Bxx)},
  MRNUMBER = {3243734},
MRREVIEWER = {Atanas\ G.\ Stefanov},
       DOI = {10.1007/978-1-4939-1194-3},
       URL = {https://doi.org/10.1007/978-1-4939-1194-3},
}

@article {MR2763004,
    AUTHOR = {Lerner, Andrei K.},
     TITLE = {Some remarks on the {F}efferman-{S}tein inequality},
   JOURNAL = {J. Anal. Math.},
  FJOURNAL = {Journal d'Analyse Math\'ematique},
    VOLUME = {112},
      YEAR = {2010},
     PAGES = {329--349},
      ISSN = {0021-7670,1565-8538},
   MRCLASS = {42B25 (46E15)},
  MRNUMBER = {2763004},
MRREVIEWER = {Yoshihiro\ Sawano},
       DOI = {10.1007/s11854-010-0032-1},
       URL = {https://doi.org/10.1007/s11854-010-0032-1},
}

@misc{lerner2026dual,
  title={Dual boundedness of the maximal operator on concave Banach function spaces}, 
  author={Andrei K. Lerner},
  year={2026},
  eprint={2608.01140},
  archivePrefix={arXiv},
  primaryClass={math.CA},
  doi={10.48550/arXiv.2608.01140},
  url={https://arxiv.org/abs/2608.01140}
}

@book {MR961439,
    AUTHOR = {Dauge, Monique},
     TITLE = {Elliptic boundary value problems on corner domains},
    SERIES = {Lecture Notes in Mathematics},
    VOLUME = {1341},
      NOTE = {Smoothness and asymptotics of solutions},
 PUBLISHER = {Springer-Verlag, Berlin},
      YEAR = {1988},
     PAGES = {viii+259},
      ISBN = {3-540-50169-X},
   MRCLASS = {35J67 (35J40 47F05 58G20)},
  MRNUMBER = {961439},
MRREVIEWER = {H.\ Triebel},
       DOI = {10.1007/BFb0086682},
       URL = {https://doi.org/10.1007/BFb0086682},
}

@book {MR775683,
    AUTHOR = {Grisvard, P.},
     TITLE = {Elliptic problems in nonsmooth domains},
    SERIES = {Monographs and Studies in Mathematics},
    VOLUME = {24},
 PUBLISHER = {Pitman (Advanced Publishing Program), Boston, MA},
      YEAR = {1985},
     PAGES = {xiv+410},
      ISBN = {0-273-08647-2},
   MRCLASS = {35J25 (35-02)},
  MRNUMBER = {775683},
MRREVIEWER = {P.\ Szeptycki},
}

@incollection{Lerner2017,
  author    = {Lerner, Andrei K.},
  title     = {On a dual property of the maximal operator
               on weighted variable {$L^p$} spaces},
  booktitle = {Functional Analysis, Harmonic Analysis,
               and Image Processing:
               A Collection of Papers in Honor of
               Bj{\"o}rn Jawerth},
  editor    = {Cwikel, Michael and Milman, Mario},
  series    = {Contemporary Mathematics},
  volume    = {693},
  pages     = {283--300},
  publisher = {American Mathematical Society},
  year      = {2017},
  doi       = {10.1090/conm/693/13932},
  eprint    = {1509.07664},
  archivePrefix = {arXiv}
}

@book{CruzUribeMartellPerez2011,
  author    = {Cruz-Uribe, David V. and
               Martell, Jos{\'e} Mar{\'i}a and
               P{\'e}rez, Carlos},
  title     = {Weights, Extrapolation and the Theory
               of Rubio de Francia},
  series    = {Operator Theory: Advances and Applications},
  volume    = {215},
  publisher = {Birkh{\"a}user},
  year      = {2011},
  doi       = {10.1007/978-3-0348-0072-3}
}

@incollection{CruzUribe2017,
  author    = {Cruz-Uribe, David},
  title     = {Extrapolation and Factorization},
  booktitle = {Function Spaces, Embeddings and
               Extrapolation X},
  editor    = {Luke{\v{s}}, Jaroslav and Pick, Lubo{\v{s}}},
  pages     = {45--92},
  publisher = {Matfyzpress, Charles University},
  address   = {Prague},
  year      = {2017},
  note      = {Spring School on Analysis, Paseky 2017},
  eprint    = {1706.02620},
  archivePrefix = {arXiv}
}

@book {MR1156767,
    AUTHOR = {Kokilashvili, Vakhtang and Krbec, Miroslav},
     TITLE = {Weighted inequalities in {L}orentz and {O}rlicz spaces},
 PUBLISHER = {World Scientific Publishing Co., Inc., River Edge, NJ},
      YEAR = {1991},
     PAGES = {xii+233},
      ISBN = {981-02-0612-7},
   MRCLASS = {42B25 (26D10 42B10 46E30 47B38 47G10)},
  MRNUMBER = {1156767},
MRREVIEWER = {Vladimir\ D.\ Stepanov},
       DOI = {10.1142/9789814360302},
       URL = {https://doi.org/10.1142/9789814360302},
}

@article {MR2837636,
    AUTHOR = {Cruz-Uribe, David and Diening, Lars and H\"ast\"o, Peter},
     TITLE = {The maximal operator on weighted variable {L}ebesgue spaces},
   JOURNAL = {Fract. Calc. Appl. Anal.},
  FJOURNAL = {Fractional Calculus and Applied Analysis. An International
              Journal for Theory and Applications},
    VOLUME = {14},
      YEAR = {2011},
    NUMBER = {3},
     PAGES = {361--374},
      ISSN = {1311-0454,1314-2224},
   MRCLASS = {42B25 (42B35)},
  MRNUMBER = {2837636},
       DOI = {10.2478/s13540-011-0023-7},
       URL = {https://doi.org/10.2478/s13540-011-0023-7},
}

@book {MR2790542,
    AUTHOR = {Diening, Lars and Harjulehto, Petteri and H\"ast\"o, Peter and
              Růžička, Michael},
     TITLE = {Lebesgue and {S}obolev spaces with variable exponents},
    SERIES = {Lecture Notes in Mathematics},
    VOLUME = {2017},
 PUBLISHER = {Springer, Heidelberg},
      YEAR = {2011},
     PAGES = {x+509},
      ISBN = {978-3-642-18362-1},
   MRCLASS = {46-02 (26D10 31B15 35J60 35Q35 46E30 46E35 46N20)},
  MRNUMBER = {2790542},
MRREVIEWER = {Dorothee\ D.\ Haroske},
       DOI = {10.1007/978-3-642-18363-8},
       URL = {https://doi.org/10.1007/978-3-642-18363-8},
}

@book {MR3026953,
    AUTHOR = {Cruz-Uribe, David V. and Fiorenza, Alberto},
     TITLE = {Variable {L}ebesgue spaces},
    SERIES = {Applied and Numerical Harmonic Analysis},
      NOTE = {Foundations and harmonic analysis},
 PUBLISHER = {Birkh\"au\-ser\slash Springer, Heidelberg},
      YEAR = 2013,
     PAGES = {x+312},
      ISBN = {978-3-0348-0547-6},
   MRCLASS = {42-02 (42B20 42B25 42B35 46E30)},
  MRNUMBER = 3026953,
MRREVIEWER = {Paul\ Alton\ Hagelstein},
       DOI = {10.1007/978-3-0348-0548-3},
       URL = {https://doi.org/10.1007/978-3-0348-0548-3},
}

@article {MR4951202,
    AUTHOR = {Adamadze, Daviti and Diening, Lars and Kopaliani, Tengiz},
     TITLE = {Maximal operator on variable exponent spaces},
   JOURNAL = {J. Math. Anal. Appl.},
  FJOURNAL = {Journal of Mathematical Analysis and Applications},
    VOLUME = {554},
      YEAR = {2026},
    NUMBER = {2},
     PAGES = {Paper No. 129972, 21},
      ISSN = {0022-247X,1096-0813},
   MRCLASS = {42B25 (42B35 46E30)},
  MRNUMBER = {4951202},
MRREVIEWER = {Julian\ Weigt},
       DOI = {10.1016/j.jmaa.2025.129972},
       URL = {https://doi.org/10.1016/j.jmaa.2025.129972},
}

@article{Oswald1990,
  author  = {Oswald, Peter},
  title   = {On Function Spaces Related to Finite Element Approximation Theory},
  journal = {Z. Anal. Anwend.},
  volume  = {9},
  number  = {1},
  pages   = {43--64},
  year    = {1990},
  doi     = {10.4171/ZAA/380}
}

@article{MartinPerkowski2019,
  author  = {Martin, J{\"o}rg and Perkowski, Nicolas},
  title   = {Paracontrolled Distributions on {Bravais} Lattices
             and Weak Universality of the 2d Parabolic {Anderson} Model},
  journal = {Ann. Inst. Henri Poincar{\'e} Probab. Stat.},
  volume  = {55},
  number  = {4},
  pages   = {2058--2110},
  year    = {2019},
  doi     = {10.1214/18-AIHP942}
}

@book{Grafakos2014,
  author    = {Grafakos, Loukas},
  title     = {Modern Fourier Analysis},
  series    = {Graduate Texts in Mathematics},
  volume    = {250},
  edition   = {3},
  publisher = {Springer},
  address   = {New York},
  year      = {2014},
  isbn      = {978-1-4939-1229-2},
  doi       = {10.1007/978-1-4939-1230-8}
}

@article {KermanTorchinsky1982,
    AUTHOR = {Kerman, R. A. and Torchinsky, A.},
     TITLE = {Integral inequalities with weights for the {H}ardy maximal
              function},
   JOURNAL = {Studia Math.},
  FJOURNAL = {Polska Akademia Nauk. Instytut Matematyczny. Studia
              Mathematica},
    VOLUME = {71},
      YEAR = {1982},
    NUMBER = {3},
     PAGES = {277--284},
      ISSN = {0039-3223,1730-6337},
   MRCLASS = {42B25 (28A25)},
  MRNUMBER = {667316},
MRREVIEWER = {Roberto\ A.\ Mac\'ias},
       DOI = {10.4064/sm-71-3-277-284},
       URL = {https://doi.org/10.4064/sm-71-3-277-284},
}

@article{FengHenningsNeilan2017,
  author  = {Feng, Xiaobing and Hennings, Lauren and Neilan, Michael},
  title   = {Finite element methods for second order linear elliptic
             partial differential equations in non-divergence form},
  journal = {Mathematics of Computation},
  volume  = {86},
  number  = {307},
  pages   = {2025--2051},
  year    = {2017},
  doi     = {10.1090/mcom/3168}
}

@article{SmearsSuli2013,
  author  = {Smears, Iain and S{\"u}li, Endre},
  title   = {Discontinuous {Galerkin} finite element approximation of
             nondivergence form elliptic equations with {Cord{\`e}s}
             coefficients},
  journal = {SIAM Journal on Numerical Analysis},
  volume  = {51},
  number  = {4},
  pages   = {2088--2106},
  year    = {2013},
  doi     = {10.1137/120899613}
}

@article{FengNeilanSchnake2018,
  author  = {Feng, Xiaobing and Neilan, Michael and Schnake, Stefan},
  title   = {Interior penalty discontinuous {Galerkin} methods for
             second order linear non-divergence form elliptic {PDEs}},
  journal = {Journal of Scientific Computing},
  volume  = {74},
  number  = {3},
  pages   = {1651--1676},
  year    = {2018},
  doi     = {10.1007/s10915-017-0519-3}
}

@article{Arnold1982,
  author  = {Arnold, Douglas N.},
  title   = {An interior penalty finite element method with
             discontinuous elements},
  journal = {SIAM Journal on Numerical Analysis},
  volume  = {19},
  number  = {4},
  pages   = {742--760},
  year    = {1982},
  doi     = {10.1137/0719052}
}

@article{ArnoldBrezziCockburnMarini2002,
  author  = {Arnold, Douglas N. and Brezzi, Franco and
             Cockburn, Bernardo and Marini, L. Donatella},
  title   = {Unified analysis of discontinuous {Galerkin} methods
             for elliptic problems},
  journal = {SIAM Journal on Numerical Analysis},
  volume  = {39},
  number  = {5},
  pages   = {1749--1779},
  year    = {2002},
  doi     = {10.1137/S0036142901384162}
}

@article{BrennerSung2005,
  author  = {Brenner, Susanne C. and Sung, Li-Yeng},
  title   = {{$C^0$} interior penalty methods for fourth order
             elliptic boundary value problems on polygonal domains},
  journal = {Journal of Scientific Computing},
  volume  = {22/23},
  pages   = {83--118},
  year    = {2005},
  doi     = {10.1007/s10915-004-4135-7}
}

@article{BartelsDolzmannNochetto2010,
  author  = {Bartels, S{\"o}ren and Dolzmann, Georg and Nochetto, Ricardo H.},
  title   = {A finite element scheme for the evolution of orientational
             order in fluid membranes},
  journal = {ESAIM: Mathematical Modelling and Numerical Analysis},
  volume  = {44},
  number  = {1},
  pages   = {1--31},
  year    = {2010},
  doi     = {10.1051/m2an/2009040}
}

@article{DieningKroenerRuzickaToulopoulos2014,
  author  = {Diening, Lars and Kr{\"o}ner, Dietmar and
             R{\r u}{\v z}i{\v c}ka, Michael and Toulopoulos, Ioannis},
  title   = {A local discontinuous {Galerkin} approximation
             for systems with {$p$}-structure},
  journal = {IMA Journal of Numerical Analysis},
  volume  = {34},
  number  = {4},
  pages   = {1447--1488},
  year    = {2014},
  doi     = {10.1093/imanum/drt040}
}

@article{Gallistl2015,
  author  = {Gallistl, Dietmar},
  title   = {{Morley} finite element method for the eigenvalues
             of the biharmonic operator},
  journal = {IMA Journal of Numerical Analysis},
  volume  = {35},
  number  = {4},
  pages   = {1779--1811},
  year    = {2015},
  doi     = {10.1093/imanum/dru054}
}

@article{CarstensenPuttkammer2020,
  author  = {Carstensen, Carsten and Puttkammer, Sophie},
  title   = {How to prove the discrete reliability
             for nonconforming finite element methods},
  journal = {Journal of Computational Mathematics},
  volume  = {38},
  number  = {1},
  pages   = {142--175},
  year    = {2020},
  doi     = {10.4208/jcm.1908-m2018-0174}
}

@article{CarstensenNataraj2021,
  author  = {Carstensen, Carsten and Nataraj, Neela},
  title   = {A priori and a posteriori error analysis of the
             {Crouzeix--Raviart} and {Morley} {FEM}
             with original and modified right-hand sides},
  journal = {Computational Methods in Applied Mathematics},
  volume  = {21},
  number  = {2},
  pages   = {289--315},
  year    = {2021},
  doi     = {10.1515/cmam-2021-0029}
}

@article{KarkulikMelenk2015,
  author  = {Karkulik, Michael and Melenk, Jens Markus},
  title   = {Local high-order regularization and applications to
             {$hp$}-methods},
  journal = {Computers \& Mathematics with Applications},
  volume  = {70},
  number  = {7},
  pages   = {1606--1639},
  year    = {2015},
  doi     = {10.1016/j.camwa.2015.06.026}
}

@book{BennettSharpley1988,
  author    = {Bennett, Colin and Sharpley, Robert},
  title     = {Interpolation of Operators},
  series    = {Pure and Applied Mathematics},
  volume    = {129},
  publisher = {Academic Press},
  address   = {Boston, MA},
  year      = {1988},
  isbn      = {0-12-088730-4}
}

@article{BennettRudnick1980,
  author  = {Bennett, Colin and Rudnick, Karl},
  title   = {On {Lorentz--Zygmund} spaces},
  journal = {Dissertationes Mathematicae (Rozprawy Matematyczne)},
  volume  = {175},
  year    = {1980},
  pages   = {67}
}

@article{LeykekhmanVexler2016,
  author  = {Leykekhman, Dmitriy and Vexler, Boris},
  title   = {Finite element pointwise results on convex polyhedral domains},
  journal = {SIAM Journal on Numerical Analysis},
  volume  = {54},
  number  = {2},
  pages   = {561--587},
  year    = {2016},
  doi     = {10.1137/15M1013912}
}

@misc{DieningStorn2026,
  author     = {Diening, Lars and Storn, Johannes},
  title      = {Optimal {AFEM} convergence for any bulk parameter smaller than one},
  year       = 2026,
  eprint     = {2610.39719},
  eprinttype = {arxiv},
  doi        = {10.48550/arXiv.2610.39719}
}

@book {MR506343,
    AUTHOR = {Kre\u in, S. G. and Petunin, Ju.\ I. and Sem\"enov, E. M.},
     TITLE = {{Interpolyatsiya line\u{i}nykh operatorov}},
 PUBLISHER = {``Nauka'', Moscow},
      YEAR = {1978},
     PAGES = {400},
   MRCLASS = {46M35 (46E30 46E35)},
  MRNUMBER = {506343},
MRREVIEWER = {E.\ Gerlach},
}

@article {MR2316449,
    AUTHOR = {Kurtz, Douglas S.},
     TITLE = {Classical operators on mixed-normed spaces with product
              weights},
   JOURNAL = {Rocky Mountain J. Math.},
  FJOURNAL = {The Rocky Mountain Journal of Mathematics},
    VOLUME = {37},
      YEAR = {2007},
    NUMBER = {1},
     PAGES = {269--283},
      ISSN = {0035-7596,1945-3795},
   MRCLASS = {42B25 (42B30)},
  MRNUMBER = {2316449},
MRREVIEWER = {Paul\ Alton\ Hagelstein},
       DOI = {10.1216/rmjm/1181069331},
       URL = {https://doi.org/10.1216/rmjm/1181069331},
}

@article {MR909030,
    AUTHOR = {Fernandez, Dicesar Lass},
     TITLE = {Vector-valued singular integral operators on {$L^p$}-spaces
              with mixed norms and applications},
   JOURNAL = {Pacific J. Math.},
  FJOURNAL = {Pacific Journal of Mathematics},
    VOLUME = {129},
      YEAR = {1987},
    NUMBER = {2},
     PAGES = {257--275},
      ISSN = {0030-8730,1945-5844},
   MRCLASS = {42B20 (47B38 47G05)},
  MRNUMBER = {909030},
MRREVIEWER = {Norman\ J.\ Weiss},
       URL = {http://projecteuclid.org/euclid.pjm/1102690575},
}

\end{document}